\documentclass[]{amsart}

\usepackage[Symbol]{upgreek}

\usepackage[all,arc,curve,color,frame,line]{xypic}

\usepackage{amssymb}
\usepackage{epic}
\usepackage{mathrsfs}
\usepackage{graphicx}
\usepackage{accents}
\usepackage{listliketab}
\usepackage{stmaryrd}
\usepackage{hyperref}

\renewcommand{\marginpar}[1]{}

\numberwithin{equation}{section}

\newcommand{\bbold}{\mathbb}

\def\R { {\bbold R} }
\def\Q { {\bbold Q} }
\def\Z { {\bbold Z} }
\def\C { {\bbold C} }
\def\N { {\bbold N} }

\def\c {\operatorname{c}}
\def \T{\mathbb{T}}
\def \G {\operatorname{G}}
\def \I {\operatorname{I}}

\def \ex{\operatorname{e}}

\renewcommand\epsilon{\varepsilon}

\def \<{\langle}
\def \>{\rangle}

\def \tilde {\widetilde}

\def \inv{\operatorname{inv}}

\def \hat {\widehat}

\def \supp {\operatorname{supp}}

\def \RxE {\R [[  x^{\R}  ]] ^{\operatorname {E}}  }
\def \((  {(\!(}
\def \)) {)\!)}

\def \RlognE {\R [[  \ell_{n}^\R  ]] ^{\operatorname {E}}  }
\def \RlognplusE {\R [[  \ell_{n+1}^\R  ]] ^{\operatorname {E}}  }
\def \CM {C[[\fM]]}
\def \CN {C[[\fN]]}

\def \Hom{\operatorname{Hom}}

\def \res{\operatorname{res}}

\def \k {{{\boldsymbol{k}}}}
\def \l {{{\boldsymbol{l}}}}
\def \flatter{\mathrel{\prec\!\!\!\prec}}
\DeclareMathSymbol{\precequ}{\mathrel}{symbols}{"16}
\DeclareMathSymbol{\succequ}{\mathrel}{symbols}{"17}
\def \flattereq{\mathrel{\precequ\!\!\!\precequ}}
\def \steeper{\mathrel{\succ\!\!\!\succ}}
\def \steepereq{\mathrel{\succequ\!\!\!\succequ}}
\def \comp{\mathrel{-{\hskip0.06em\!\!\!\!\!\asymp}}}

\def \nasymp{\not\asymp}

\newcommand{\claim}[2][\!\!]{\medskip\noindent {\sc Claim #1.} {\it #2}\medskip}

\newtheorem{theorem}{Theorem}[section]
\newtheorem{lemma}[theorem]{Lemma}
\newtheorem{prop}[theorem]{Proposition}
\newtheorem{cor}[theorem]{Corollary}

\newtheorem*{theoremA}{Theorem~A}
\newtheorem*{theoremB}{Theorem~B}
\newtheorem*{corC}{Corollary~C}
\newtheorem*{corD}{Corollary~D}

\theoremstyle{definition}

\theoremstyle{remark}
\newtheorem*{example}{Example}
\newtheorem*{examples}{Examples}

\newtheorem*{notation}{Notation}

\newtheorem*{remarks}{Remarks}
\newtheorem*{remark}{Remark}

\newtheorem{remarkNumbered}[theorem]{Remark}

\newcommand{\abs}[1]{\lvert#1\rvert}

\def \fG {{\mathfrak G}}
\def \fH {{\mathfrak H}}
\def \fM {{\mathfrak M}}
\def \fN {{\mathfrak N}}
\def \fP {{\mathfrak P}}
\def \fR {{\mathfrak R}}
\def \fS {{\mathfrak S}}
\def \fd {{\mathfrak d}}
\def \fg {{\mathfrak g}}
\def \fm {{\mathfrak m}}
\def \fn {{\mathfrak n}}
\def \fr {{\mathfrak r}}
\def \fv {{\mathfrak v}}
\def \fu {{\mathfrak u}}
\def \fw {{\mathfrak w}}

\def \CM {{{C}}[[ \fM ]] }

\def \RM {{{\R}}[[ \fM ]] }

\def \id{\operatorname{id}}

\let\oldi\i
\let\oldj\j
\renewcommand\i{\relax\ifmmode{\boldsymbol{i}}\else\oldi\fi}
\renewcommand\j{\relax\ifmmode{\boldsymbol{j}}\else\oldj\fi}

\renewcommand\leq{\leqslant}
\renewcommand\geq{\geqslant}
\renewcommand\preceq{\preccurlyeq}
\renewcommand\succeq{\succcurlyeq}
\renewcommand\le{\leq}
\renewcommand\ge{\geq}
\renewcommand\frak{\mathfrak}

\DeclareMathAlphabet{\mathbf}{OML}{cmm}{b}{it}

\DeclareFontFamily{U}{fsy}{}
\DeclareFontShape{U}{fsy}{m}{n}{<->s*[.9]psyr}{}
\DeclareSymbolFont{der@m}{U}{fsy}{m}{n}
\DeclareMathSymbol{\der}{\mathord}{der@m}{182}

\DeclareSymbolFont{der@m}{U}{fsy}{m}{n}
\DeclareMathSymbol{\derdelta}{\mathord}{der@m}{100}

\DeclareSymbolFont{imag@m}{OT1}{cmr}{m}{ui}
\DeclareMathSymbol{\imag}{\mathord}{imag@m}{105}

\DeclareFontFamily{OMS}{smallo}{}
\DeclareFontShape{OMS}{smallo}{m}{n}{<->s*[.65]cmsy10}{}
\DeclareSymbolFont{smallo@m}{OMS}{smallo}{m}{n}
\DeclareMathSymbol{\smallo}{\mathord}{smallo@m}{79}

\DeclareFontFamily{OMS}{largerdot}{}
\DeclareFontShape{OMS}{largerdot}{m}{n}{<->s*[.8]cmsy10}{}
\DeclareSymbolFont{largerdot@m}{OMS}{largerdot}{m}{n}
\DeclareMathSymbol{\largerdot}{\mathord}{largerdot@m}{15}

\newcommand\Aut{\operatorname{Aut}}

\DeclareMathSymbol{\llambda}{\mathord}{der@m}{108}
\DeclareMathSymbol{\rrho}{\mathord}{der@m}{114}

\def \upo{\upomega}

\def \SigmaAut{\Sigma\!\Aut}

\title[Strong Automorphisms of Transseries]{The Group of Strong Automorphisms of the Differential Field of Transseries}

\begin{document}

\author[Aschenbrenner]{Matthias Aschenbrenner}
\address{Kurt G\"odel Research Center for Mathematical Logic\\
Universit\"at Wien\\
1090 Wien\\ Austria}
\email{matthias.aschenbrenner@univie.ac.at}

\author[van den Dries]{Lou van den Dries}
\address{Department of Mathematics\\
University of Illinois at Urbana-Cham\-paign\\
Urbana, IL 61801\\
U.S.A.}
\email{vddries@illinois.edu}

\author[van der Hoeven]{Joris van der Hoeven}
\address{CNRS, LIX (UMR 7161)\\ 
Campus de l'\'Ecole Polytechnique\\  91120 Palaiseau \\ France}
\email{vdhoeven@lix.polytechnique.fr}

\date{November 2026}

\begin{abstract}
We describe the group of strongly additive automorphisms of the differential field $\T$ of transseries as a non-trivial semidirect
product. We show that this group is torsion-free and has trivial center, and that many naturally arising differential subfields of $\T$ are definably closed in $\T$. 

\end{abstract}

\maketitle

\section*{Introduction}
 
\noindent
A {\it transseries}\/  
is a formal series in an indeterminate $x$ with real coefficients and monomials  involving exponential
and logarithmic terms, such as
\begin{align*}
   &   - 3 \ex^{\ex^x} + \ex^{\textstyle\frac{\ex^x}{\log x} +
  \frac{\ex^x}{\log^2 x} + \frac{\ex^x}{\log^3 x} + \cdots} - x^{11} + 7\\
  &   \hspace{2em} + \frac{\pi}{x} + \frac{1}{x \log x} + \frac{1}{x \log^2
  x} + \frac{1}{x \log^3 x} + \cdots \nonumber\\
  &   \hspace{2em} + \frac{2}{x^2} + \frac{6}{x^3} + \frac{24}{x^4} +
  \frac{120}{x^5} + \frac{720}{x^6} + \cdots \nonumber\\
  &   \hspace{2em} + \ex^{- x} + 2 \ex^{- x^2} + 3 \ex^{- x^3} + 4 \ex^{-
  x^4} + \cdots, \nonumber
\end{align*}
where $\log^2 x := (\log x)^2$, etc.~Transseries are 
also known as {\it logarithmic-exponential series}\/ (for brevity: LE-series).
They form a real closed ordered  subfield $\T$ of a certain ordered field $\R[[G^{\operatorname{LE}}]]$
of well-based  series (also known as Hahn series) with real coefficients. (Section~\ref{sec:prelims, well-based} below has basic facts about well-based series.) Here~$G^{\operatorname{LE}}$ is an  ordered subgroup of the multiplicative group of~$\T$ containing  $x$, whose elements are called
 {\it transmonomials.}\/ The {\em raison d'\^etre\/} for~$\T$ is that,
 unlike its subfield $\R(\!(t)\!)$ of Laurent series in $t=x^{-1}$ or the larger Hahn field $\R[[G^{\operatorname{LE}}]]$, $\T$ carries a natural 
 {\it exponential map,}\/  that is, an
isomorphism $\exp$ from the ordered  additive group of $\T$ onto the   
ordered multiplicative group of positive elements of~$\T$. We set $\ex^f:=\exp(f)$ for $f\in\T$ and denote the inverse of~$\exp$ by~$\log$. 

See \cite{DMM2} and Appendix~A of~\cite{ADH} for a detailed construction of $\T$;
we recall some salient points in Section~\ref{sec:prelims, trans}.
The main result of this paper requires the notion of  a family of elements of~$\T$ being {\it summable}\/ (in~$\T$);
for a summable  family 
$(f_\lambda)_{\lambda\in \Lambda}$ in~$\T$, its  sum $\sum_{\lambda} f_\lambda$ is defined as an element of~$\T$.
A  map $\Phi\colon\T\to\T$ is then said to be {\it strongly additive}\/
if for each summable family~$(f_\lambda)_{\lambda\in \Lambda}$ in~$\T$  
the family~$\big(\Phi(f_\lambda)\big){}_{\lambda\in \Lambda}$ is summable
 and $\Phi(\sum_\lambda f_\lambda) = \sum_\lambda \Phi(f_\lambda)$.
There is a unique strongly additive derivation~$\der$ on~$\T$ such that $\der(r)=0$ for all $r\in\R$, $\der(x)=1$, and 
$\der(\ex^f)=\der(f)\ex^f$ for all $f\in\T$.
We set $f':= \der(f)$, $f'':= \der^2(f)$, and so on, 
for $f\in \T$.  Unless explicitly stated otherwise, $\T$ will denote from now on
the {\it differential}\/ field~$\T$ with this derivation, and so by an automorphism of $\T$ we mean an automorphism of the differential field $\T$. Since $\T$ is real closed, its ordering is definable, hence an automorphism of $\T$
is also an automorphism of $\T$ as an ordered field. 
The constant field of 
$\T$ is $\R$, hence an automorphism of $\T$ restricts to the identity on $\R$, and hence is $\R$-linear. The convex hull of $\R$ in the ordered field $\T$ is a valuation ring of~$\T$, making $\T$ a valued field, so any automorphism of $\T$ is also an automorphism of this valued  field.  Here is the first main result, proved in Section~\ref{sec:proof of main thm} below:

\begin{theoremA}
Let $\alpha\colon\T\to\R$ be an additive map which vanishes on the convex hull of~$\R$ in $\T$,  and let $c\in\R$; then there is a unique strongly additive automorphism~$\sigma$ of~$\T$ such that $\sigma(x)=x+c$ and $\sigma(\ex^f)=\ex^{\alpha(f)+\sigma(f)}$ for all $f\in\T$.
Moreover, each strongly additive automorphism   of $\T$ arises in this way from a unique pair $(\alpha,c)$.
\end{theoremA}

\noindent
As a consequence the inverse of a strongly additive automorphism of $\T$
is also strongly additive, and hence the strongly additive  automorphisms of $\T$ form a group~$\SigmaAut_\der(\T)$ under composition.
We have the subgroups  
\begin{align*} \mathcal T\ &:=\ \big\{\sigma\in\SigmaAut_\der(\T):\ \sigma(\ex^f) = \ex^{\sigma(f)} \text{ for all }f\in \T \big\},\\
\mathcal E\  &:=\ \big\{ \sigma\in\SigmaAut_\der(\T):\ \sigma(x)=x \big\}
\end{align*}
of $\SigmaAut_\der(\T)$. 
The group $\mathcal T$ is abelian, indeed, isomorphic to the additive group of~$\R$, but is not a normal subgroup of $\SigmaAut_\der(\T)$. On the other hand,
$\mathcal E$ is a  non-abelian normal subgroup of $\SigmaAut_\der(\T)$, of size $2^{2^{\aleph_0}}$.  By Theorem~A,
$$\SigmaAut_\der(\T)\ =\ \mathcal E\rtimes\mathcal T\qquad \text{(internal semidirect product)},$$
see Corollary~\ref{cor:semidirect product}.  An element $\sigma\in \mathcal T$ is uniquely determined by the real constant~$c$ such that $\sigma(x)=x+c$, and accordingly we call $\mathcal T$ the {\it translation group\/} of $\T$, and we call its elements {\it translations}. 
 A group of differential automorphisms similar to~$\mathcal E$ in the context of meromorphic differential equations over~$\C$ was considered by  Martinet and Ramis~\cite{MR1,MR2}. Copying their terminology
we call $\mathcal E$ the {\it exponential torus}\/ of~$\T$. (See also \cite[Section~3.2]{vdPS}.)
Our other main result is the following: 

\begin{theoremB}
The group $\SigmaAut_\der(\T)$ is torsion-free and has trivial center.
\end{theoremB}


\noindent
Given a set  $\mathcal G$ of automorphisms of  $\T$, the fixed field
$$\T^{\mathcal G} := \big\{ f\in\T: \text{$\sigma(f)=f$ for all $\sigma\in\mathcal G$} \big\}$$
of $\mathcal G$ is   definably closed in the differential field~$\T$, and thus
a real closed differential subfield of $\T$. 
Conversely, given a differential subfield~$K$
of~$\T$, we can hope to show that $K$ is definably closed in $\T$ by identifying a subgroup $\mathcal G$ of $\SigmaAut_\der(\T)$ such that
$K=\T^{\mathcal G}$. 
For example, $\R=\T^{\mathcal T}$ by Corollary~\ref{cor:R definably closed}, so~$\R$ is definably closed in~$\T$. 
Using the same approach and exploiting the flexibility for constructing  automorphisms of $\T$ afforded by
Theorem~A,  we treat in Section~\ref{sec:filtrations}  some other naturally arising differential subfields of $\T$:

\begin{corC}
The following differential subfields of $\T$ are definably closed in $\T$:
$$\R[[\ex_n^\R\cdots\ex_0^\R]],\quad \R[[\ell_0^\R\cdots\ell_n^\R]],  \quad
\T_{\exp}, \quad \T_{\log}\quad \text{ and } \quad\T_{\log}^\Q.$$ 
\end{corC}

\noindent
Here $\ex_n$, $\ell_n$ are transmonomials defined recursively by $\ex_0=\ell_0=x$ and $\ex_{n+1}=\exp(\ex_n)$,
$\ell_{n+1}=\log(\ell_n)$; $\T_{\exp}$ is the differential field of {\it exponential trans\-series}\/
(see~\cite[Appen\-dix~A]{ADH} or \cite[Section~2.1]{DMM2});  
$\T_{\log} = \bigcup_n \R[[\ell_0^\R\cdots\ell_n^\R]]$
is the differential field of {\it logarithmic transseries}\/; and
$\T_{\log}^\Q$ its differential sub\-field~$\bigcup_n \R[[\ell_0^\Q\cdots\ell_n^\Q]]$.
Indeed, we'll see that $\T_{\log}^\Q=\T^{\mathcal L_{\infty}}$ where 
$$\mathcal L_{\infty}\ =\ 
\big\{\sigma\in\SigmaAut_\der(\T):\ \text{$\sigma(\ell_n)=\ell_n$ for all $n$}\big\}.$$
In  Section~\ref{sec:filtrations} we also give a description of those differential subfields of $\T$  
which arise as
 the fixed fields of subgroups
 of $\mathcal L_{\infty}$. One   interesting example   is the
differential field~$\T^\flat$ of {\it flat transseries}, introduced in Section~\ref{sec:prelims, trans}. It plays a key role in the proof of
Theorem~B. For a general study of ``definable closure" in closed $H$-fields, see \cite{ADHdc}. 

Finally, we establish a result (Corollary~\ref{cor:move exps}) which allows us to complete the proof of a fact 
 stated in \cite{ADHdim}:

\begin{corD}
The $0$-definable subset $S:=\big\{y\in\T:yy''=(y')^2\big\}$ of $\T$ is not internal to~$\R$: for no $n$ is there a definable $f\colon \R^n\to \T$ with $S\subseteq f(\R^n)$.
\end{corD}

\subsection*{Organization} Sections~\ref{sec:prelims, well-based} and \ref{sec:prelims, trans} contain preliminaries on well-based series and transseries, much of it summarizing what we need from the literature, but also with some new material. (The relevant literature includes \cite{DMM2, vdH, Sch}.) Section~\ref{sec:logtrans} studies the differential field  $\T_{\log}$ and Section~\ref{sec:Tlog}  determines its group of strongly additive automorphisms, which turns out to be abelian, divisible, and torsion-free. Section~\ref{sec:proof of main thm} proves Theorem~A. Sections~\ref{sec:group} and~\ref{sec:center} go into more detail about  the group~$\SigmaAut_\der(\T)$, and include the proof of Theorem~B.
In Section~\ref{sec:filtrations} we prove Corollaries~C and~D. Sections~\ref{sec:center} and ~\ref{sec:filtrations} are independent. 

\subsection*{Notations and terminology} Throughout, $m$ and $n$ range over the set $\N=\{0,1,2,\dots\}$
of natural numbers. ``Ordered'' means ``totally ordered'' in this paper. If $A$ is an abelian group written additively, then $A^{\ne}:=A\setminus\{0\}$, and if $A$ is in addition an ordered abelian group, then $A^>:=\{a\in A:\, a>0\}$,
and likewise with~$\ge$,~$<$, or~$\le$ in place of $>$. Rings are commutative with $1$, and if
$R$ is ring, then $R^\times$ denotes the multiplicative group of units of $R$.
Thus if $K$ is a field, then~$K^{\ne}=K^\times$, and if $K$ is an ordered field, then $K^{>}$ is its multiplicative group of positive elements.  For us a {\rm differential\/}  field $K$ contains $\Q$ as a subfield and is equipped with a single derivation $\der\colon K\to K$, with $C_K$ denoting its constant field. If $\der$ is clear from the context, then  for $a\in K$ we also write $a'$, $a''$ and so on in stead of $\der(a)$, $\der^2(a)$, and we set $a^\dagger:= a'/a$ for $a\ne 0$. Note that then  for
$a,b\in K^\times$ we have $(ab)^\dagger=a^\dagger+b^\dagger$,  and  $a^\dagger=b^\dagger\Leftrightarrow a/b\in C_K$.

\subsection*{Acknowledgements}
Work on this paper began about ten years ago in the wake of \cite{ADHdim}, and has been supported
by NSF   Grant DMS-1700439 (Aschenbrenner) and the Institut Henri Poincar\'e (Aschenbrenner, van den Dries).
Van der Hoeven has been supported by an ERC-2023-ADG grant for the
ODELIX project (number 101142171). Funded by the European Union. Views and opinions expressed are however those
of the authors only and do not necessarily reflect those of the European
Union or the European Research Council Executive Agency. Neither the European
Union nor the granting authority can be held responsible for them.

\subsection*{AI-Declaration} No AI was used in the research for this paper or its writing. 

\section{Well-Based Series}\label{sec:prelims, well-based}

\noindent
In this section we mainly fix terminology and notation concerning well-based series and state some basic facts about them, mostly well-known. We pay special attention to Hahn fields and strongly additive maps on them. An essential feature of $\T$ is that it is
not a Hahn field, but a directed union of Hahn fields, so we also include some generalities about such directed unions. 

\subsection*{Well-based series}
In this section $\frak{M}$ is a set (totally) ordered by $\preceq$.
We think of the elements of~$\mathfrak{M}$ as {\em monomials\/} and we let $\fm$,~$\fn$
range over $\fM$. A sub\-set~$\mathfrak{S}$ of $\frak{M}$ 
is said to be {\it well-based}\/
if it is well-ordered for the reverse ordering~$\succeq$, that is, 
there is no strictly increasing infinite sequence 
$\fm_0 \prec \fm_1 \prec \frak{m}_2 \prec \cdots$ in~$\mathfrak{S}$.
Next, let   $C$ be a field
 whose elements we think of as {\it coefficients.}\/ Write a function $f\colon \frak{M} \to C$ as  a series
$\sum_{\mathfrak{m}} f_{\mathfrak{m}} \mathfrak{m}$,
$f_{\mathfrak{m}}:=f(\mathfrak{m})$; its 
{\it support\/} is
$\supp f := \{ \mathfrak{m}   : f_{\mathfrak{m}}
\neq 0 \}$. Then 
\[ C [[ \mathfrak{M} ]]\ :=\ \left\{ f\colon  \frak{M} \to C :\,  
   \text{$\supp f \subseteq \mathfrak{M}$ is well-based} \right\} \]
is a subspace of the $C$-linear space  of all functions $\mathfrak M\to C$. 
 We identify $\mathfrak{M}$ with a subset of $\CM$ in the obvious way. Given any ordered subset $\fN$ of $\fM$ we have the $C$-linear subspace $\CN$ of $\CM$.

\medskip
\noindent
For nonzero $f\in \CM$ we let
$\mathfrak d(f) := \max \supp f$
be the {\it dominant monomial}\/ of~$f$; also $\mathfrak d(0):=0\in \CM$. We extend the ordering of $\mathfrak{M}$ to an ordering of
the subset~$\mathfrak{M}\cup \{0\}$ of $\CM$ by making $0$ its least element. The map~${\mathfrak d\colon \CM\to \mathfrak{M}\cup\{0\}}$
is then a valuation on the abelian group $\CM$ as defined in~\cite[Section~2.2]{ADH}, but with the reverse ordering on its value set~$\mathfrak{M}\cup\{0\}$. This makes the notation introduced in that section of \cite{ADH} available; for
example, given $f,g\in\CM$,  
$$f\preceq g\ \Leftrightarrow\ \mathfrak d(f)\preceq\mathfrak d(g),\qquad f\prec g\ \Leftrightarrow\  \mathfrak d(f)\prec\mathfrak d(g),$$
so the binary relations $\preceq$ and $\prec$ on $\CM$ restrict to the previously defined relations~$\preceq$ and~$\prec$ on
$\mathfrak{M}\cup\{0\}$. 
If $C$ comes equipped with an ordering making it an ordered field, then we equip $C[[\fM]]$ with the ordering~$\leq$
making it an ordered abelian group such that for all  $f\in C[[\fM]]$,
$$
0\leq f \quad\Longleftrightarrow\quad 0 \leq f_{\mathfrak d(f)}.
$$
This ordering  makes $C[[\fM]]$ an ordered vector space over $C$; cf.~\cite[p.~106]{ADH}. The restriction of this ordering to $\fM$ is the given ordering of $\fM$. 





\medskip 
\noindent
Let~$\square$ denote any of the symbols $\succ$, $\succeq$, $\prec$, $\preceq$. For $f\in\CM$  we set
$$f_{\square\fn}\ :=\ \sum_{\fm\square \fn} f_{\fm}\fm \in \CM\qquad (\text{so }\supp(f_{\square\fn})=\{\fm\in \supp f:\, \fm\square \fn\}).$$ 
Thus $f=f_{\succ\fn}+f_{\preceq\fn}$ and $f_{\preceq\fn}=f_{\fn}\fn+f_{\prec\fn}$. 
Given a subset $S$ of $C[[\fM]]$,  put 
$$S_{\square\fn}\ :=\ \{f\in S: f_{\square\fn}=f\}\ =\  \{f\in S:\supp f \,\square\, \fn\}.$$ 
If $S$ is an additive subgroup, respectively  a $C$-linear subspace, of $C[[\fM]]$, then so is~$S_{\square\fn}$. 
The map
$$f\mapsto f_{\square\fn}\ \colon\ \CM\to\CM_{\square\fn}$$ is $C$-linear, and $\CM$ is the following internal direct sums of $C$-linear subspaces:
$$C[[\fM]]\ =\  C[[\fM]]_{\succ\fn}\oplus  C\fn \oplus C[[\fM]]_{\prec\fn}.$$
We call $f_{\succ\fn}=\sum_{\fm \succ \fn} f_{\fm}\fm$ the {\it truncation $f_{\succ\fn}$ of $f$ at $\fn$}.
A subset $S$ of~$\CM$ is said to be {\it truncation closed\/} if $S_{\succ\fn}\subseteq S$ for all~$\fn$, that is, for all $f\in S$ and monomials~$\fn$ we have $f_{\succ\fn}\in S$. If $S$ is a
truncation closed $C$-linear subspace of~$C[[\mathfrak{M}]]$,
then~$\supp f\subseteq S$ for all $f\in S$.

\subsection*{Summing well-based series}
In the rest of this section $\Lambda$ denotes an index set and $\lambda$ ranges over $\Lambda$.
A family $(f_\lambda)$ in $\CM$ is said to be
{\it summable\/} if $\bigcup_\lambda \supp f_\lambda$ is well-based and
for each $\fm$ there are only finitely many $\lambda$ such that $\fm\in\supp f_\lambda$; in that case we define
its sum $\sum_\lambda f_\lambda$ to be the series $f\in \CM$ such that $f_{\fm}=\sum_\lambda f_{\lambda,\fm}$ for each $\fm$. (This agrees with the usual notation  for elements of $\CM$: for $f\colon \fM\to C$,  $f\in \CM$ iff the family~$(f_{\fm}\fm)$ is summable; in that case the just defined sum $\sum_{\fm} f_{\fm}\fm\in C[[\fM]]$ equals $f$.)
Here are some basic facts about this notion with routine proofs and often used tacitly below:

\begin{lemma}\label{lem:summ part}
Let $(f_\lambda)$ be a summable family in~$\CM$.
\begin{enumerate}
\item[\rm(i)] If $\pi\colon \textup M\to\Lambda$ is a bijection, then the family $(f_{\pi(\mu)})_{\mu\in\textup M}$
is summable and $\sum_{\mu} f_{\pi(\mu)} = \sum_\lambda f_\lambda$;
\item[\rm(ii)] if $(c_\lambda)$ is any family in~$C$ and $(g_\lambda)$ is a summable family in $\CM$  then the families $(c_\lambda f_\lambda)$ and $(f_\lambda+g_\lambda)$ are summable; we have  $\sum_\lambda (f_\lambda+g_\lambda)=\sum_\lambda f_\lambda + \sum_\lambda g_\lambda$, and if $c_\lambda=c$ for all $\lambda$ then $\sum_\lambda cf_\lambda=c\sum_\lambda f_\lambda$;
\item[\rm(iii)] if $(\Lambda_i)_{i\in I}$ is a family of pairwise disjoint subsets of $\Lambda$ with
$\Lambda=\bigcup_{i\in I} \Lambda_i$, then $(f_\lambda)_{\lambda\in\Lambda_i}$ is summable for each $i\in I$,
$\big(\sum_{\lambda\in\Lambda_i} f_\lambda\big){}_{i\in I}$ is summable, and~$\sum_{\lambda\in\Lambda} f_\lambda = \sum_{i\in I}
\sum_{\lambda\in\Lambda_i} f_\lambda$.
\end{enumerate}
\end{lemma}


\subsection*{Strongly linear maps}
Let $\fN$ also be an ordered set, allowing of course $\fN=\fM$. A map $\Phi\colon\CM \to \CN$ is said to be {\it strongly additive\/} if for each summable family $(f_\lambda)$ in $\CM$ the family~$\big(\Phi(f_\lambda)\big)$ in~$\CN$ is summable with $\Phi\big(\sum_\lambda f_\lambda\big)=\sum_\lambda \Phi(f_\lambda)$. 
If $\Phi$ is strongly additive, then it is additive. 
We say that $\Phi$ is {\it strongly linear}\/ if $\Phi$ is $C$-linear and strongly additive;
then for each summable family $(f_\lambda)$ in $\CM$ and each family $(c_\lambda)$ in $C$ we have
$\Phi(\sum_\lambda c_\lambda f_\lambda)=\sum_\lambda c_\lambda\Phi(f_\lambda)$.

\begin{example}
For $\fn\in\fM$ and $\square$ being $\succ$, $\succeq$, $\prec$, $\preceq$, the map 
$$f\mapsto f_{\square\fn}\ \colon\  \CM\to\CM$$ 
is strongly linear.
\end{example}

\noindent
The following is easy to verify:

\begin{lemma}\label{lem:add strongly additive}
Let $\Phi,\Phi^*\colon \CM\to \CN$  be strongly additive and $c\in C$. Then the maps $\Phi+\Phi^*\colon \CM\to \CN$ and $c\Phi\colon \CM\to \CN$  given by
$$(\Phi+\Phi^*)(f)\ :=\ \Phi(f)+\Phi^*(f), \qquad (c\Phi)(f)\ :\ =c\Phi(f) $$ 
are also strongly additive.
If $\fP$ is  an ordered  set and $\Psi\colon \CN\to C[[\fP]]$ is strongly additive, then  
 $\Psi\circ\Phi\colon \CM\to C[[\fP]]$ is strongly additive.
\end{lemma}

\noindent
A map $\phi\colon\fM \to \CN$ is called {\it well-based}\/ if for each well-based subset $\mathfrak S$ of~$\fM$,
the family $\big(\phi(\fm)\big){}_{\fm\in\mathfrak S}$ in $\CN$ is summable; equivalently,  for each 
sequence $(\fm_i,\fn_i)_{i\in\N}$ in $\fM\times\fN$ such that
$$\fm_0 \succ \fm_1 \succ\cdots\succ \fm_i\succ\fm_{i+1}\succ\cdots\quad \text{and}\quad\text{$\fn_i\in\supp \phi(\fm_i)$ for all $i\in\N$}$$ 
there are $i<j$ such that
$\fn_i\succ\fn_j$. 
For example, every map $\phi\colon\fM\to\CN$ such that~$\phi(\fM)\subseteq C\fN$ and  $\fm_1\prec \fm_2 \Rightarrow \phi(\fm_1) \prec  \phi(\fm_2)$ for all $\fm_1,\fm_2\in\fM$ with~${\phi(\fm_2)\neq 0}$ is well-based.
Every strongly linear map $\CM\to\CN$ restricts to a well-based map $\fM\to\CN$; conversely, by \cite[Proposition~3.5]{vdH:noeth}:

\begin{prop}\label{prop:extend to strongly linear}
Every well-based map $\fM\to\CN$ extends uniquely to a strongly linear map $\CM\to\CN$.
\end{prop}

\noindent
In the next proposition (a special case of \cite[Theorem~6.2]{vdH:noeth}) and its three corollaries we assume $\Phi\colon\CM\to\CM$ 
is strongly linear with $\Phi(\fm)\prec\fm$ for all $\fm$. 

\begin{prop}\label{prop:6.2}
Let $\fS\subseteq\fM$ be well-based. Then the family $\big(\Phi^n(\fm)\big){}_{(\fm,n)\in \fS \times\N}$ is summable.
\end{prop}
\begin{proof}
Let $i$, $j$, $k$ range over $\N$.
It suffices to show:
\begin{enumerate}
\item $\bigcup_n \supp \Phi^n(\fS)$ is well-based;
\item for each $\fn$ there are only finitely many   $(\fm,n)\in\fS\times\N$ with $\fn\in\supp\Phi^n(\fm)$.
\end{enumerate}
Call a tuple $\tau=(\fm_0,\dots,\fm_n)$ in $\fM^{1+n}$ a {\it path}\/ of length $n$ from $\fm$ to $\fn$ if~$\fm_0=\fm\in\fS$,
$\fm_{i+1}\in\supp\Phi(\fm_{i})$
for $i<n$, and $\fm_{n}=\fn$. For such $\tau$ we have ${\fm_0\succ\cdots\succ\fm_n}$, and we set $\fm(\tau):=\fm_n$.
For $\fm\in\fS$ and $\fn\in\supp\Phi^n(\fm)$ there is a path $\tau$  of length $n$
from $\fm$ to $\fn$. 
Moreover, if~$\tau=(\fm_0,\dots,\fm_n)$ is a path of length $n\ge 1$,
then~$\check{\tau}:=(\fm_0,\dots,\fm_{n-1})$ is a path of length $n-1$,
with $\fm(\check\tau)\succ\fm(\tau)$. We shall obtain (1) from the following a-priori stronger claim:

\begin{enumerate}
\item[(1)$^*$]
the set $\big\{\fm(\tau):\,  \text{$\tau$ is a path (of some length)}\big\}$ is well-based. 
\end{enumerate}
Suppose (1)$^*$ fails. Then
we have a  sequence $(\tau_i)$ 
of paths  $\tau_i=(\fm_{i,0},\dots,\fm_{i,n_i})$  
which is ``bad'', that is, $\fm(\tau_i)\prec\fm(\tau_j)$ for all $i<j$. 
We may assume that this bad sequence is minimal in the following sense:
for no $i$ is there a  path $\tau$   
of length~$<n_i$ such that $\tau_0,\dots,\tau_{i-1},\tau$
can be continued to a bad sequence. 

Since $\fS$ is well-based, we have $n_i\ge 1$ for all but finitely many $i\in\N$.
We claim that $\check{\fS}:=\big\{\fm(\check{\tau}_i):n_i\ge 1\big\}$ is well-based.
Suppose towards a contradiction that we have 
 a  sequence $(i_k)$ in $\N$ such that
  $n_{i_k}\ge 1$  and $\fm(\check{\tau}_{i_k}) \prec \fm(\check{\tau}_{i_{k+1}})$ for all $k$.
 Then $\tau_0,\dots,\tau_{i_0-1},\check\tau_{i_0},\check\tau_{i_1},  \dots$
is also a bad sequence of paths, which contradicts the minimality of our original bad sequence~$(\tau_i)$.
Hence $\check{\fS}$ is well-based.  As $\Phi$ is strongly linear,  $\supp \Phi(\check{\fS})$
is also well-based, with
$\fm(\tau_i)\in \supp\Phi\big(\fm(\check{\tau}_i)\big)\subseteq\supp\Phi(\check{\fS})$: a contradiction. This proves (1)$^*$ and thus (1). 

We shall obtain (2) from the following a-priori stronger claim:

\begin{enumerate}
\item[(2)$^*$]
for each $\fn$ there are only finitely many  $(\fm,n)\in\fS\times\N$ for which there is a path   of length $n$ from $\fm$ to  $\fn$. 
\end{enumerate}
Suppose towards a contradiction that (2)$^*$ fails.  Then (1)$^*$ yields a largest counterexample $\fn$.
Take distinct pairs $(\fm_i, n_i)\in \fS\times \N$, $i\in \N$, with a path $\tau_i$ of length~$n_i$ from~$\fm_i$
to $\fn$. If $n_i=0$, then $\fm_i=\fn$, hence $n_i = 0$ for  at most one $i$. 
Now the subset $\big\{\fm(\check{\tau}_i):n_i\ge 1\big\}$ of the well-based set in (1)$^*$
is well-based.  Let any $\fm$ be given. Then there are only finitely many $i$ with $n_i\ge 1$ and $\fm(\check{\tau}_i)=\fm$: otherwise, $\fm$ is a counterexample to (2)$^*$ and $\fm\succ \fn$, contradicting the maximality of~$\fn$. 
 Hence by strong linearity of $\Phi$  the family $\big(\Phi(\fm(\check{\tau}_i))\big){}_{n_i\ge 1}$ is summable,
so there are only finitely many $i$ such that $n_i\ge 1$ and $\fn\in\supp\Phi\big(\fm(\check{\tau}_i)\big)$, contradicting that~$\fn\in\supp\Phi(\fm(\check{\tau}_i))$ for all $i$ with $n_i\ge 1$. This proves (2)$^*$ and thus (2). 
\end{proof}

\begin{cor}\label{cor:cnPhin}
Let $(c_n)$ be a family in~$C$. Then for each $f\in C[[\fM]]$,  the family~$\big(c_n\Phi^n(f)\big)$
is summable, and the map
\begin{equation}\label{eq:sumcnPhin}
f\mapsto \sum_{n=0}^\infty c_n\Phi^n(f)\ \colon\  C[[\fM]]\to C[[\fM]]
\end{equation}
is strongly linear.
\end{cor}
\begin{proof}
Let $\fS\subseteq\fM$ be well-based. By Proposition~\ref{prop:6.2}, the family
 $\big(\Phi^n(\fm)\big){}_{(\fm,n)\in \fS \times\N}$ is summable, and hence
 so is  $\big(c_n\Phi^n(\fm)\big){}_{(\fm,n)\in \fS \times\N}$.
Put $\Psi(\fm):=\sum_n c_n\Phi^n(\fm)$. Then~$\big(\Psi(\fm)\big){}_{\fm\in\fS}$ is summable
by Lemma~\ref{lem:summ part}(iii). Hence   
$$\fm\mapsto\Psi(\fm)\ :\  \fM\to C[[\fM]]$$ is well-based and
thus extends uniquely, by Proposition~\ref{prop:extend to strongly linear}, to the strongly linear map \eqref{eq:sumcnPhin}.
\end{proof}

\noindent
For a formal power series $P=\sum_n c_n t^n\in C[[t]]$ ($c_n\in C$)
 we denote the strongly linear map
\eqref{eq:sumcnPhin} by $P(\Phi)$. 
The map $P\mapsto P(\Phi)$ is a $C$-linear ring morphism from~$C[[t]]$
to the ring of strongly linear endomorphisms of $\CM$; in fact:

\begin{cor}\label{cor:P(Phi)}
Let $P,Q\in C[[t]]$. Then 
$$(P+Q)(\Phi)\ =\ P(\Phi)+Q(\Phi)\quad\text{ and }\quad (P\cdot Q)(\Phi)\ =\ P(\Phi)\circ Q(\Phi).$$
If $Q\in tC[[t]]$, then $Q(\Phi)(\fm)\prec\fm$ for all $\fm\in \fM$ and $(P\circ Q)(\Phi)=P\big(Q(\Phi)\big)$.
\end{cor}

\begin{cor}\label{cor:Inversion-of-Strongly-Linear-Operators}
The strongly linear map $\operatorname{id}+\Phi$ on $C[[\frak{M}]]$
is bijective with strongly linear inverse given by
\[(\operatorname{id} + \Phi)^{-1}(g)\ =\ \sum_{n=0}^\infty (-1)^n \Phi^n(g).\]
\end{cor}

\noindent
 Corollary~\ref{cor:Inversion-of-Strongly-Linear-Operators} is  for infinite $C$ also contained in \cite[Corollary~1.4]{ADHdagap}, where $\fM$ is not assumed totally ordered. 
Together with Lemma~\ref{lem:add strongly additive} this corollary yields:

\begin{cor}\label{prop:inverse strongly linear}
Let $\Psi\colon\CM\to\CM$ be strongly linear such that 
$\Psi(\fm)\sim \fm$ for all $\fm$. Then $\Psi$ is bijective with strongly linear inverse $\Psi^{-1}\colon\CM\to\CM$.
\end{cor}

\noindent
The proof of the corollary below is similar to that of \cite[Lemma~1.7]{ADHdagap}.

\begin{cor}\label{cor:extend to strongly linear}
Suppose $\Phi\colon\CM\to\CN$ is a strongly linear bijection such that~$\fm_1\prec \fm_2\Rightarrow\Phi(\fm_1)\prec\Phi(\fm_2)$ for all $\fm_1,\fm_2 \in \fM$. Then $\Phi^{-1}$ is strongly linear.
\end{cor}
\begin{proof}  From the hypothesis we obtain that $f\prec g\Rightarrow
\Phi(f)\prec\Phi(g)$ for all $f,g$ in~$\CM^{\ne}$. It follows that for every $\fn\in \fN$ there are unique $c\in C^\times$ and $\fm\in \fM$ such that~$\Phi(c\fm)\sim \fn$. 
Define $\I\colon\fN\to\CM$ by $\I(\fn):=c\fm$ with $c\in C^\times$, $\fm$ as above. 
Then by Proposition~\ref{prop:extend to strongly linear}, $\I$ extends to a
strongly linear map $\CN\to\CM$, also denoted by $\I$.  Now the strongly linear
operator $\Psi=\Phi\circ\I$ on $\CN$ satisfies~${\Psi(\fn)\sim \fn}$ for all $\fn\in \fN$.
Hence by Corollary~\ref{prop:inverse strongly linear},  $\Psi$ has a strongly linear inverse~$\Psi^{-1}$, and so~$\Phi^{-1}=\I\circ\Psi^{-1}$
is strongly linear.
\end{proof}

\subsection*{Strong subgroups} A subgroup $A$ of the additive group of $\CM$  is said to be {\it strong}\/  if
for each fam\-i\-ly~$(f_\lambda)$ in $A$ which is summable in $\CM$ we have
$\sum_\lambda f_\lambda \in A$. 
For example, if $\mathfrak S$ is an ordered subset of $\fM$, then $C[[\mathfrak S]]$ is a strong $C$-linear subspace of $\CM$. Moreover:

\begin{lemma}\label{lem:fixsigma, 0}
Let $S$ be a strong $C$-linear subspace of $\CM$. Then the following conditions on $S$ are equivalent:
\begin{enumerate}
\item[\rm(i)] $\mathfrak d(S^{\neq})\subseteq S$;
\item[\rm(ii)] $S=C[[\mathfrak S]]$ where $\mathfrak S=\mathfrak d(S^{\neq})$;
\item[\rm(iii)] $S=C[[\mathfrak S]]$ for some $\mathfrak S\subseteq\fM$;
\item[\rm(iv)] $S$ is truncation closed;
\item[\rm(v)] $\supp f\subseteq S$ for all $f\in S$. 
\end{enumerate}
\end{lemma}
\begin{proof}
Suppose $\mathfrak S:=\mathfrak d(S^{\neq})\subseteq S$.
If $f\in C[[\mathfrak S]]$, then $f_\fm \fm\in S$ for all $\fm$, so~$f=\sum f_{\fm} \fm\in S$. Thus $C[[\mathfrak S]]\subseteq S$.
Towards a contradiction, assume $f\in S\setminus C[[\mathfrak S]]$. Let $\fm$ be the largest element of $\supp f \setminus \mathfrak S$. Then $g:=f_{\succ\fm}\in C[[\mathfrak S]]\subseteq S$,  so~$\fm=\mathfrak d(f-g)\in\mathfrak d(S^{\neq})=\mathfrak S$, a contradiction.
This shows (i)~$\Rightarrow$~(ii), and the im\-pli\-ca\-tions~(ii)~$\Rightarrow$~(iii)~$\Rightarrow$~(iv)~$\Rightarrow$~(v)~$\Rightarrow$~(i)  are clear. 
\end{proof}

\noindent
The strong truncation closed $C$-linear subspaces of $\CM$ form a distributive sublattice of the lattice
of all $C$-linear subspaces of $\CM$:

\begin{cor}\label{cor:distributive}
If the $C$-linear subspaces $U$, $V$, $W$ of $\CM$ are strong and truncation closed, then so are the $C$-linear subspaces $U+V$ and $U\cap V$ of $\CM$, and~$(U+V)\cap W=(U\cap W)+(V\cap W)$, $(U\cap V)+W=(U+W)\cap(V+W)$. 
\end{cor}
\begin{proof}
Use that if $\mathfrak U,\mathfrak V\subseteq\fM$ with $U=C[[\mathfrak U]]$ and 
$V=C[[\mathfrak V]]$, then   
$U+V=C[[\mathfrak U\cup\mathfrak V]]$ and
$U\cap V=C[[\mathfrak U\cap\mathfrak V]]$.
\end{proof}

\noindent
In the proof of Lemma~\ref{lem:Em-1,n} we need the following variant:

\begin{cor}\label{cor*:distributive} Let the $C$-linear subspaces $U, V$ of $\CM$ be strong and truncation closed, and 
$W= \bigcup_m W_m$, where $W_0\subseteq W_1\subseteq W_2\subseteq \cdots$ is an increasing sequence of strong and truncation closed $C$-linear subspaces of $\CM$. Then $$U\cap(V+ W)\ =\ (U\cap V)+(U\cap W).$$
\end{cor}
\begin{proof} We have $U\cap(V+W_m)=(U\cap V)+(U\cap W_m)$ for all $m$, by Corollary~\ref{cor:distributive}, and
taking unions over all $m$ then yields the desired result. 
\end{proof}

\subsection*{Multiplying  well-based series} In the rest of this section $\fM$ is equipped with a binary operation making it
 a multiplicatively written abelian group
(so with identity element~$1$),  such that for all $\fm_1$,~$\fm_2$,~$\fn$ in $\fM$ we have
$\fm_1\preceq\fm_2  \Rightarrow  \fm_1\fn\preceq \fm_2\fn$. So~$\fM$ is an ordered abelian group. 

\begin{notation}
Given subsets $\mathfrak S$, $\mathfrak T$  of $\fM$,
we put 
$$\mathfrak S\mathfrak T\  :=\ \{\mathfrak s\mathfrak t:\mathfrak s\in\mathfrak S,\mathfrak t\in
\mathfrak T\}$$
and
$$\mathfrak S^n\ :=\ \{\mathfrak s_1\cdots\mathfrak s_n:\ \mathfrak s_1,\dots,\mathfrak s_n\in\mathfrak S\}
\qquad\text{with $\mathfrak S^0=\{1\}$,}$$
and we let
$$\mathfrak S^{*}\ :=\ \bigcup_n \mathfrak S^n$$ 
be the submonoid of $\fM$ generated by $\mathfrak S$.
\end{notation}

\noindent
The following fact from \cite{Neumann} is an important source for well-based subsets of $\fM$:

\begin{lemma}[Neumann's Lemma]\label{lem:Neumann}
If $\mathfrak S,\mathfrak T\subseteq\fM$ are well-based, then
so is $\mathfrak S\mathfrak T$, and for each $\fm\in\mathfrak S\mathfrak T$
there are only finitely many pairs $(\mathfrak s,\mathfrak t)\in\mathfrak S\times\mathfrak T$
with~${\fm=\mathfrak s\mathfrak t}$.
If $\mathfrak S\subseteq\fM^{\prec 1}$ is well-based then
$\mathfrak S^{*}$ is well-based, and for each $\fm\in \mathfrak S^{*}$
there are only finitely many tuples $(n,\mathfrak s_1,\dots,\mathfrak s_n)$ with 
$\mathfrak s_1,\dots,\mathfrak s_n\in\mathfrak S$ such that $\fm=\mathfrak s_1\cdots\mathfrak s_n$.
\end{lemma}

\noindent
As a consequence,
for $f,g\in\CM$,
the family $(f_{\fm_1}g_{\fm_2} \fm_1\fm_2)_{(\fm_1,\fm_2)\in\fM\times\fM}$ is summable,
and denoting its sum by $f\cdot g$ we have
$$f\cdot g\ =\ \sum_{\mathfrak m\in\mathfrak M} \left(\sum_{\mathfrak m_1\cdot\mathfrak m_2=\mathfrak m} f_{\mathfrak m_1}g_{\mathfrak m_2}\right)\mathfrak m.$$
Equipped with the binary operation
$(f,g)\mapsto f\cdot g$,   $\CM$  is a $C$-algebra, and also a field~\cite[Lem\-ma~3.1.3]{ADH}. We identify $C$ with a subfield of  $C[[\mathfrak M]]$   via the field embedding $c\mapsto f$ with
$f_{1}=c$ and $f_{\mathfrak{m}}=0$ for all $\mathfrak{m}\ne 1$. Note that $\fM$ is subgroup of 
the multiplicative group $C[[\mathfrak M]]^{\times}$. If $\fN$ is a submonoid of~$\fM$, then $\CN$ is a subalgebra of the $C$-algebra $\CM$.
Multiplication in $\CM$ is compatible with infinite summation:

\begin{lemma}[{see \cite[Proposition~3.3]{vdH:noeth}}] \label{lem:summ mult}
If $(f_\lambda)_{\lambda\in\Lambda}$, $(g_{\mu})_{\mu\in\textup M}$ are summable families in $\CM$,
then the family $(f_\lambda g_\mu)_{(\lambda,\mu)\in\Lambda\times\textup M}$ in $\CM$ is
also summable, and
$$\left(\sum_\lambda f_\lambda\right) \cdot \left(\sum_{\mu} g_\mu\right) = 
\sum_{\lambda,\mu} f_\lambda g_\mu.$$
In particular, for $g\in\CM$ the operator $f\mapsto fg$ on $\CM$ is strongly linear: if~$(f_{\lambda})$ is a summable family in $\CM$, then so is
$(f_\lambda g)$ and $\sum_\lambda f_\lambda g  = \left( \sum_\lambda f_\lambda\right)g$.
\end{lemma}

\noindent
Using Lemmas~\ref{lem:summ part} and \ref{lem:summ mult} one shows easily:

\begin{cor}\label{cor:strongly linear field embedd}
If $\fN$ is also ordered abelian group,  then each strongly linear map~$\Phi\colon \CM\to\CN$ satisfying
$$\Phi(\fm_1\fm_2)=\Phi(\fm_1)\Phi(\fm_2)\quad\text{ for all $\fm_1,\fm_2\in\fM$}$$ 
is a ring morphism.
Moreover, each strongly linear map $\Delta\colon\CM\to\CM$ with
$$\Delta(\fm_1\fm_2)=\Delta(\fm_1)\fm_2+\fm_1\Delta(\fm_2)\quad\text{ for all $\fm_1,\fm_2\in\fM$}$$
is a derivation.
\end{cor}

\noindent
For $f\in C[[\fM]]$ and~$\square$ any of the symbols $\succ$, $\succeq$, $\prec$, $\preceq$,
we  set $f_\square:=f_{\square 1}$.  Similarly, given a subset $S$ of $\CM$ and $\square$ as before,  put 
$S_\square:=S_{\square 1}$. Then $\CM$ is the internal direct sums of $C$-linear subspaces
$$C[[\fM]]\ =\  C[[\fM]]_\succ\oplus  C \oplus C[[\fM]]_\prec.$$
This yields the valuation ring
$\mathcal O:=C[[\fM]]_{\preceq}=C\oplus C[[\fM]]_\prec$ of  $C[[\fM]]$, with maximal ideal $\smallo:=C[[\fM]]_{\prec}$ of $\mathcal O$.
Below we view $K=C[[\fM]]$ as a valued field with this valuation ring.
The group morphism $f\mapsto  \mathfrak d(f) \colon K^\times\to\fM$ has kernel~$\mathcal O^\times$, hence induces an isomorphism between
the (ordered) value group of $K$ and the group $\fM$ with the reversed ordering $\succeq$.
The surjective ring morphism $f\mapsto f_1\colon\mathcal O\to C$ has kernel~$\smallo$ and thus induces an isomorphism $\res(K)=\mathcal O/\smallo\to C$ between the residue field of~$K$ and the coefficient field $C$.
By Corollary~\ref{cor:extend to strongly linear}:

\begin{cor}\label{cor:strongly linear inverse}
The inverse of each strongly linear automorphism of the valued field $\CM$  is also strongly linear.
\end{cor}

\noindent
We say that a strongly linear operator $\Phi$ on $\CM$ is {\it diagonal}\/ if $\Phi(\fm)\in C\fm$ for all $\fm$.
For each map $\chi\colon \fM\to C$ the map $\fm\mapsto\chi(\fm)\fm$ is well-based, so we have a strongly linear   map $f\mapsto f_{\times\chi}\colon\CM\to\CM$ given by
$$f_{\times\chi}\ :=\ \sum_\fm f_\fm\chi(\fm)\fm.$$ 
The operator  $f\mapsto f_{\times\chi}$ on $\CM$ is diagonal;
every diagonal strongly linear operator on $\CM$ arises in this way from a unique map $\chi\colon\fM\to C$.
For $\chi_1, \chi_2\colon \fM\to C$, 
$$f_{\times (\chi_1+\chi_2)}\ =\ f_{\times \chi_1}+f_{\times \chi_2},\qquad f_{\times \chi_1\chi_2}\ =\ f_{\times \chi_1}\circ f_{\times\chi_2}\ =\ \big( f_{\times \chi_1}\big){}_{\times \chi_2}.$$
Using this and Corollary~\ref{cor:strongly linear field embedd} we obtain:

\begin{lemma}
If $\chi\colon\fM\to C^\times$ is a morphism from the group $\fM$ to the multiplicative group of $C$, then $f\mapsto f_{\times\chi}$ is an automorphism of the field $\CM$. If $\chi\colon\fM\to C$ is a morphism from the group $\fM$ to the additive group of $C$, then $f\mapsto f_{\times\chi}$ is a derivation on $\CM$.
\end{lemma}


\noindent
An automorphism of the field $\CM$ is said to be {\it strong}\/ if it and its inverse are strongly linear. Let $\SigmaAut(\CM)$ be the group of strong
automorphisms of the field~$\CM$.
Let $\operatorname{Hom}(\fM,C^\times)$ be the abelian group of all group morphisms~${\fM\to C^\times}$ with group operation given by pointwise multiplication.
Then the map which assigns to~${\chi\in \operatorname{Hom}(\fM,C^\times)}$  the automorphism $f\mapsto f_{\times\chi}$
of $\CM$ is a group embedding of $\operatorname{Hom}(\fM,C^\times)$ into~$\SigmaAut(\CM)$. 

\begin{remarkNumbered}\label{rem:S}
Suppose $C$ has characteristic zero. The formal power series
$$\exp(t):=\sum_{n}  \frac{1}{n!}t^n\in 1+tC[[t]],\qquad
\log(1+t):=\sum_{n\geq 1} \frac{(-1)^{n-1}}{n}t^n\in tC[[t]]$$
satisfy~${\log\circ\exp}=t$ and $\exp\circ\log(1+t)=1+t$.
Let $\der$ be a strongly linear derivation on $\CM$ such that~$\der(\fm)\prec\fm$ for each $\fm$.
Then we have the strongly linear map~$\exp(\der)\colon\CM\to\CM$ (see Corollary~\ref{cor:cnPhin}),
and one easily verifies (using the Product Rule for~$\der$) that $\sigma:=\exp(\der)$ is a field embedding
with $\sigma(\fm)\sim\fm$ for each~$\fm$,
 and hence~$\sigma\in\SigmaAut(\CM)$ by Corollary~\ref{prop:inverse strongly linear}. 
 (Theorem~3.13 in~\cite{BKKPS} shows conversely that each~$\sigma\in\SigmaAut(\CM)$ satisfying~$\sigma(\fm)\sim\fm$
for each $\fm$ arises this way;
 we will not use this later.) Next consider the formal power series
 $$S := \frac{t}{\exp(t)-1}\  =\   
1-\frac{t}{2}+\sum_{n=1}^\infty \frac{B_{2n}}{(2n)!} t^{2n}\ =\  
1-\frac{t}{2}+\frac{t^2}{12} + \text{higher powers of } t,$$
where the $B_{2n}$ are the Bernoulli numbers.
Then 
$$\big(\exp(\der)-\id\big)\circ S(\der)\  =\  
\big( (\exp(t)-1)\cdot S \big)(\der)\  =\  
\der$$ by Corollary~\ref{cor:P(Phi)}.
Likewise, for $c\in C\setminus\{1\}$ we have 
$$S_c\ :=\  \frac{1}{c\exp(t)-1}\  =\  \frac{1}{c-1}-\frac{c}{(c-1)^2}t+\frac{c(c+1)}{2(c-1)^3}t^2 + \text{higher powers of } t,$$
and
$$\big(c\exp(\der)-\id\big)\circ S_c(\der)\  =\  
\big( (c\exp(t)-1)\cdot S_c \big)(\der)\  =\  \id.$$
\end{remarkNumbered}

\subsection*{Ordering well-based series}
Suppose $C$ is an ordered field. Then $\CM$ with the induced ordering is an ordered field extension of $C$ with $\fM$ as an ordered subgroup of~$\CM^>$.
We call the ordered field $C[[\fM]]$ an {\it ordered Hahn field.}\/ The valuation ring $\CM_\preceq$ of $\CM$ is the convex hull of $C$ in $\CM$.
The  ordered Hahn field~$\CM$ is real closed iff the field $C$ is real closed and the group~$\fM$ is divisible (see for example
\cite[Section~3.5]{ADH}).

\subsection*{Decomposing $\fM$}
Suppose $\fG$ and $\fR$ are ordered subgroups of $\fM$ such that 
$$ \text{$\fG$ is convex in $\fM$}, \qquad \fG\cap \fR\ =\ \{1\}, \quad \fM\ =\ \fG\fR:=\{\fg\fr:\ \fg\in \fG,\ \fr\in \fR\}.$$
Let $f\in C[[\fM]]$. Then for each $\fr\in \fR$ the  family $(f_{\fg\fr}\fg)_{\fg\in \fG}$ is summable in~$C[[\fG]]$, and the family
$\left(\sum_{\fg\in \fG}f_{\fg\fr}\fg\right)_{\fr\in\fR}$ is summable
in $C[[\fG]][[\fR]]$.
The map
$$f\ =\ \sum_{\fm}f_\fm \fm\ \mapsto\ \sum_{\fr\in \fR}\left(\sum_{\fg\in \fG}f_{\fg\fr}\fg\right)\fr\ \colon\ \CM\to C[[\fG]][[\fR]]$$
is a $C[[\fG]]$-algebra isomorphism.
For $f\in C[[\fM]]$ we have in fact $$f\ =\ \sum_{(\fg,\fr)\in \fG\times \fR} f_{\fg\fr}\fg\fr\ =\ \sum_{\fr\in \fR}\left(\sum_{\fg\in \fG}f_{\fg\fr}\fg\right)\fr\ \text{ in }\ \CM$$ where the indicated sums exist in $C[[\fM]]$ according to the definition of summability, using
Lemma~\ref{lem:summ part}.
Whenever convenient we 
identify below~$C[[\fM]]$ and $C[[\fG]][[\fR]]$ via the above isomorphism. 
If in addition $C$ is an ordered field, then $C[[\fM]]$ and~$C[[\fG]]$ are ordered Hahn fields, and so is
$C[[\fG]][[\fR]]$, and the above isomorphism is also an isomorphism of ordered fields.

\subsection*{Substitution into power series}
Let $t_1,\dots, t_n$ be distinct indeterminates and $t:=(t_1,\dots, t_n)$. Let
$$c\ =\ c(t)\ =\ \sum_{\nu}c_{\nu}t^{\nu}\in C[[t]]\ :=\ C[[t_1,\dots, t_n]]$$
be a formal power series over $C$; here and below $\nu=(\nu_1,\dots, \nu_n)$ ranges over $\N^n$, and
$c_{\nu}\in C$, $t^{\nu}:= t_1^{\nu_1}\cdots t_n^{\nu_n}$. For any tuple $\varepsilon=(\varepsilon_1,\dots, \varepsilon_n)$ of elements of $C[[\fM]]_{\prec}$ the family $(c_{\nu}\epsilon^{\nu})$,  where $\varepsilon^{\nu}:= \varepsilon_1^{\nu_1}\cdots \varepsilon_n^{\nu_n}$, is summable, by Lemma~\ref{lem:Neumann}. 
Put
$$c(\varepsilon):= \ \sum_{\nu}c_{\nu}\varepsilon^{\nu}\in C[[\fM]]_{\preceq}.$$
Equip $\Z^n$ with an ordering $\le$ making it an ordered abelian group with $0\le \nu$ for all $\nu$, for example, the lexicographic ordering. It has $\N^n$ as a well-ordered subset, by Dickson's Lemma~\cite{Dickson}. Equip the multiplicative monoid $t^*:= t_1^{\N}\cdots t_n^{\N}$ with
the ordering $\preceq$ such that $t^{\mu}\preceq t^{\nu}\Leftrightarrow \mu \ge \nu$, for all $\mu,\nu\in \N^n$. Then $t^*$ is well-based with largest element $1=t^0$. 
Thus $C[[t^{*}]]=C[[t]]$, and in this way the notion of summability make sense
for families of formal power series in~$C[[t]]$. Indeed,  a family $(c_i)_{i\in I}$ in $C[[t]]$ is summable iff for each~$\nu$ there
are only finitely $i\in I$ such that $c_{i\nu}\neq 0$, and then
$\sum_{i} c_i=\sum_\nu \left(\sum_i c_{i\nu}\right) t^\nu$. (So this notion of summability does not depend on the chosen ordering of $\Z^n$.) 
Similarly, we have the notion of strongly additive maps $C[[t]]\to C[[\fM]]$, and so on.
Given an $n$-tuple $\varepsilon$ of elements of $C[[\fM]]_{\prec}$, the map $t^\nu\mapsto \varepsilon^{\nu}\colon t^*\to \CM$ is well-based, so the map
\begin{equation}\label{eq:substitution into power series}
c(t)\mapsto c(\varepsilon)\ \colon\ C[[t]]\to C[[\fM]]
\end{equation}
defined above is a strongly linear $C$-algebra morphism.

\medskip
\noindent
In the rest of this subsection $n=1$, $t=t_1$, and $\nu=\nu_1$ ranges over $\N$.

\begin{lemma}
Let $c=\sum_\nu c_\nu t^\nu\in C[[t]]$  
and $d\in tC[[t]]$, and set $c\circ d:=c(d)=\sum_{\nu} c_\nu d^\nu$.
Let $\varepsilon\in C[[\fM]]_\prec$. Then $d(\varepsilon)\preceq \varepsilon\prec 1$
and $(c\circ d)(\varepsilon)=c\big(d(\varepsilon)\big)$.
\end{lemma}
\begin{proof}
Using strong linearity of the substitution map \eqref{eq:substitution into power series} we get
$$(c\circ d)(\varepsilon)\ =\  \left(\sum_\nu c_\nu d^\nu\right)(\varepsilon)\ =\ \sum_\nu c_\nu d^\nu(\varepsilon)\ 
=\ \sum_\nu c_\nu d(\varepsilon)^n\ =\ c\big(d(\varepsilon)\big)$$
as claimed.
\end{proof}

\begin{examples}
Suppose $C$ has characteristic zero. Then the power series 
$\exp(t):=\sum_{n}  \frac{1}{n!}t^n$ and $\log(1+t):=\sum_{n\geq 1} \frac{(-1)^{n-1}}{n}t^n$ in $C[[t]]$ yield partial
functions 
\begin{align*} \varepsilon\ \mapsto\ \exp(\varepsilon)\ &=\ \sum_{n} \frac{1}{n!}\varepsilon^n\ :\ C[[\fM]]_{\prec} \to 1+C[[\fM]]_{\prec},\\
\delta\ \mapsto 
\log(1+\delta)\ &=\ \sum_{n\geq 1} \frac{(-1)^{n-1}}{n}\delta^n\  :\ C[[\fM]]_{\prec} \to C[[\fM]]_{\prec},
\end{align*}
the first an isomorphism of the additive subgroup $C[[\fM]]_{\prec}$ of $C[[\fM]]$
onto the multiplicative subgroup $1+C[[\fM]]_{\prec}$ of $C[[\fM]]^\times$, with inverse
$$1+\delta\ \mapsto\ \log (1+\delta)\  :\ 1+C[[\fM]]_{\prec}\to C[[\fM]]_{\prec}.$$
For $c\in C$ the binomial series $(1+t)^c:=\sum_n {c\choose n} t^n\in C[[t]]$ yields
a map
$$1+\delta\ \mapsto\ (1+\delta)^c\ :=\ \sum_n {c\choose n} \delta^n\ :\ 1+C[[\fM]]_{\prec}\to 1+C[[\fM]]_{\prec}.$$
The identity $(1+t)^c = \exp\big(c\log (1+t)\big)$ in $C[[t]]$ and the lemma above yield
$$(1+\delta)^c\ =\ \exp\big(c\log(1+\delta)\big)\qquad (\delta\in C[[\fM]]_{\prec}),$$
so $1+\delta\ \mapsto\ (1+\delta)^c$ is an endomorphism of the multiplicative group $1+C[[\fM]]_{\prec}$, and for
$\delta\in \CM_{\prec}$ and $c_1, c_2\in C$ we have
$$(1+\delta)^{c_1+c_2}\ =\ (1+\delta)^{c_1}(1+\delta)^{c_2},\qquad (1+\delta)^{c_1c_1}\ =\  \big((1+\delta)^{c_1}\big){}^{c_2}.$$ 
\end{examples}

\subsection*{Powers of well-based series}
This subsection considers a scalar multiplication $$(r,\fm)\mapsto \fm^r\ \colon\ \R\times\fM\to\fM$$ that
makes $\fM$ an ordered vector space over $\R$ (written multiplicatively). For $f$ in~$\RM^>$ we have 
$f=a\fm(1+\varepsilon)\in\RM$ where $a\in \R^>$, $\fm$, and $\varepsilon\in \RM_{\prec}$ are uniquely determined by $f$. Given such $f$ and any real number $r$, set
$$f^r\ :=\ a^r\fm^r(1+\varepsilon)^r\in \RM^{>} .$$
For $f,g\in\RM^>$ and $r,s\in\R$ this yields
$$f^{r+s}\ =\ f^rf^s,\quad f^{rs}\ =\ (f^r)^s,\quad (fg)^r\ =\ f^rg^r.$$

\subsection*{Directed unions of well-based series fields}
Let $(\fM_i)_{i\in I}$ with $I\ne \emptyset$ be a family of ordered subgroups
of the ordered multiplicative group~$\fM$ with $\fM=\bigcup_i \fM_i$. We assume that $(\fM_i)$ is {\it directed}\/: for all
$i,j\in I$ there exists~$k\in I$ with $\fM_i, \fM_j\subseteq \fM_k$. 
This yields the valued subfield  
$$K\ :=\ \bigcup_i C[[\fM_i]]$$ 
of $C[[\fM]]$. Note that $K$ is truncation closed in $C[[\fM]]$. A {\it $K$-subgroup\/} of $\fM$ is an ordered subgroup $\fG$ of $\fM$ such that~${C[[\fG]]\subseteq K}$, inside the ambient $C[[\fM]]$;  each~$\fM_i$ is a $K$-subgroup of $\fM$.  
The family $(\fM_i)$ is said to be {\it healthy} if each 
$K$-subgroup of $\fM$ is contained in some~$\fM_i$; this might depend on $C$. If $I$ is countable (the only
relevant case in this paper), then $(\fM_i)$ is healthy. In the next section we shall construe the field $\T$ of transseries as such a union $\bigcup_{m,n} \R[[G_{m,n}]]$.

 
\medskip
\noindent
{\em In the rest of this section we assume that $(\fM_i)$ is healthy}. 

\medskip
\noindent
A family
$(f_\lambda)$ in $K$ is said to be {\it summable\/} if
there exists a $K$-subgroup $\fG$ of $\fM$ such that all $f_{\lambda}\in C[[\fG]]$
and $(f_{\lambda})$ is summable as a family in~$C[[\fG]]$; note that then~$\sum_\lambda f_{\lambda}$ is defined as an element of $K$ (lying in~$C[[\fG]]$ for~$\fG$ as above). 

\begin{example}
For  $F = \sum_{\nu}c_{\nu}t^{\nu}\in C[[t]]$ with $t=(t_1,\dots, t_n)$ and $\epsilon\in (K_\prec)^n$, the family~$(c_{\nu}\epsilon^{\nu})_{\nu\in \N^n}$ is summable, with
$F(\epsilon):=\sum_{\nu}c_{\nu}\epsilon^{\nu}\in   K_\preceq$. 
\end{example}

\noindent
A subgroup $A$ of the additive group of $K$  is said to be {\it strong}\/  if
for each family~$(f_\lambda)$ in $A$ which is summable in $K$ we have
$\sum_\lambda f_\lambda \in A$. 
The intersection of a family of strong subgroups of $K$ is  a strong subgroup of $K$.
Clearly~$K$ itself is a strong $C$-linear subspace of $K$, as is    $K\cap C[[\mathfrak G]]$ for each subset~$\mathfrak G$ of~$\fM$. 

Note that $\CM$ is a valued vector space over $C$, and so any $C$-linear subspace of $\CM$ is a valued subspace of $\CM$, 
cf.~\cite[Section~2.3]{ADH}. 

\begin{lemma}\label{lem:fixsigma, 1}
For a strong $C$-linear subspace $A$ of $K$, the following are equivalent:
\begin{enumerate}
\item[\rm(i)] $\mathfrak d(A^{\neq})\subseteq A$;
\item[\rm(ii)] $A=K\cap C[[\mathfrak G]]$ where $\mathfrak G=\mathfrak d(A^{\neq})$;
\item[\rm(iii)] $A=K\cap C[[\mathfrak G]]$ for some set $\mathfrak G\subseteq\fM$;
\item[\rm(iv)] $A$ is truncation closed in $\CM$;
\item[\rm(v)] $\supp f\subseteq A$ for all $f\in A$.
\end{enumerate}
\end{lemma}
\begin{proof} Let $i$ range over $I$. Since
$A\cap C[[\fM_i]]$ is a strong subgroup of $C[[\fM_i]]$, the implications (ii)~$\Rightarrow$~(iii)~$\Rightarrow$~(iv)~$\Rightarrow$~(v)~$\Rightarrow$~(i) are obvious or clear from 
Lemma~\ref{lem:fixsigma, 0}. For~(i)~$\Rightarrow$~(ii), assume (i). Then
$\mathfrak d\big(A^{\neq}\cap C[[\fM_i]])\subseteq A^{\ne}\cap C[[\fM_i]]$. Hence by Lem\-ma~\ref{lem:fixsigma, 0} we have
$A\cap C[[\fM_i]]=C[[\mathfrak G_i]]$ where $\mathfrak G_i=\mathfrak d(A^{\ne}\cap C[[\fM_i]])$. Set 
$$\mathfrak G\ :=\ \mathfrak d(A^{\neq})\ =\ \bigcup_i \mathfrak G_i.$$ 
From (i) we get $\fM_i\cap\mathfrak G = \fM_i\cap \fd(A^{\ne}) \subseteq \fM_i\cap A \subseteq \mathfrak G_i$
and hence $\fM_i\cap \mathfrak{G}=\mathfrak G_i$,  and thus (ii). 
 \end{proof}
 
 \noindent
 Note that $K$ is not only truncation closed in $\CM$, but even closed under subseries: for any set $\fS\subseteq \fM$ and
 $f\in K$ we have $\sum_{\fm\in \fS} f_{\fm}\fm\in K$. By Lemma~\ref{lem:fixsigma, 1}, being closed under subseries remains true for any strong $C$-linear subspace of $K$ that is truncation closed in $\CM$. Moreover:
 
 \begin{cor}\label{corbestapp}  Let $A$ be a strong $C$-linear subspace of $K$ that is truncation closed in $\CM$. Then every $f\in K$ has a best approximation in $A$, that is, the subset~$\fd(f-A)$ of $\mathfrak{M}\cup \{0\}$ has a least element.
 \end{cor} 
 \begin{proof}  This is clear for $f\in A$, so assume $f\in K\setminus A$ and suppose towards a contradiction that $\fd(f-A)$ has no least element.  This yields a divergent pc-sequence $(a_{\rho})$ in $A$ such that $a_{\rho}\leadsto f$. Hence $A+Cf$ is an immediate extension of~$A$ by  \cite[Lemma~2.3.1]{ADH}.   Take $i\in I$ with $f\in C[[\mathfrak{M}_i]]$.   Set $\fn_{\rho}:=\fd(f-a_{\rho})$.
Then $f$ and~$a_{\rho}$ have the same truncation $b_{\rho}$ at $\fn_{\rho}$, so with $\mathfrak{G}_i$ as in the proof of Lemma~\ref{lem:fixsigma, 1} we have $b_{\rho}\in A\cap C[[\mathfrak{M}_i]]= C[[\mathfrak{G}_i]]$. Moreover, $f-b_{\rho}\preceq f-a_{\rho}$. Hence the $\fd(f-b_{\rho})$ are coinitial in $\fd(f-A)$, contradicting the spherical completeness of $C[[\mathfrak{G}_i]]$. 
\end{proof}


\noindent
If $A$ is a strong subring of $K$, $F\in C[[t]]$,
$t=(t_1,\dots, t_n)$, and $\epsilon\in (A_\prec)^n$, then~${F(\epsilon)\in A_\preceq}$.
Using this, the proof of the following is like that of~\cite[Lem\-ma~2.5]{DMM1}:

\begin{lemma}
Each strong valued subfield of $K$ containing $C$ is henselian.
\end{lemma}

\noindent
Suppose now that $C$ is equipped with an ordering making $C$ an ordered field. We view $K$ as an ordered subfield
of the ordered Hahn field $\CM$. The previous lemma and \cite[Corollary~3.5.19]{ADH} yield:

\begin{cor}
Let $F\supseteq C$ be a strong ordered subfield of $K$. 
Then $F$ is real closed iff $C$ is real closed and each $f\in F^{>}$ has for every $n\geq 1$ an $n$th root.\end{cor}

\noindent
Let~$\fN$ also be a multiplicatively  written ordered abelian group with
$\fN=\bigcup_{j}\fN_j$ and $(\fN_j)_{j\in J}$ a directed family of ordered subgroups of $\fN$, $J\ne \emptyset$. Let
$L=\bigcup_j C[[\fN_j]]\subseteq C[[\fN]]$, and assume~$(\fN_j)$ is healthy. Then 
a map $\Phi\colon K\to L$ is said to be {\it healthy\/}  if for each $K$-subgroup~$\fG$ of~$\fM$ there exists an $L$-subgroup $\fH$ of~$\fN$ such that $\Phi\big(C[[\fG]]\big)\subseteq C[[\fH]]$.
A map $\Phi\colon K \to L$ is said to be {\it strongly additive\/} if for every summable family 
$(f_{\lambda})$ in $K$ the family $\big(\Phi(f_{\lambda})\big)$ is summable in $L$, and~$\Phi\big(\sum_{\lambda}f_{\lambda}\big)=\sum_{\lambda}\Phi(f_{\lambda})$.   Every strongly additive map $K\to L$ is additive. 
We use {\it strongly linear}\/ as shorthand for ``strongly additive and $C$-linear''.

\begin{lemma}\label{lem:comp of strongly add}
Let $\Phi,\Phi^*\colon K\to L$ be strongly additive and $c\in C$. Then the maps~$\Phi+\Phi^*\colon K\to L$ and $c\Phi\colon K\to L$  given by
$$(\Phi+\Phi^*)(f)\ =\ \Phi(f)+\Phi^*(f), \qquad (c\Phi)(f)\ =\ c\Phi(f)\qquad (f\in K)$$ 
are strongly additive. Let $\Psi\colon L\to M$ be a strongly additive map, 
$M=\bigcup_\lambda C[[\fP_\lambda]]$ with $(\fP_\lambda)_{\lambda\in \Lambda}$ $(\Lambda\ne \emptyset)$ a healthy directed family of ordered subgroups of an ordered commutative group
 $\fP=\bigcup_\lambda\fP_\lambda$. Then
 $\Psi\circ\Phi\colon K\to M$ is strongly additive.  
\end{lemma}

\noindent
This is routine to verify.  
Lemma~\ref{lem:comp of strongly add} continues to hold if ``strongly additive'' is replaced everywhere by ``healthy''.
If $\Phi\colon K\to L$ is strongly additive,
then $\ker \Phi$ is a strong subgroup of $K$.
Hence:

\begin{cor}\label{cor:strong fixed gp}
If $\Phi\colon K\to K$ is strongly additive, then the additive subgroup
$$\ker(\Phi-\id_K)=\big\{f\in K:\Phi(f)=f\big\}$$ 
of $K$ is strong.
\end{cor} 

\noindent
We say that $\Phi\colon K\to K$ is {\it very healthy}\/ if  every $K$-subgroup $\fG$ of $\fM$ is contained in some $\fM_i$ with $\Phi(C[[\fM_i]])\subseteq C[[\fM_i]]$. If $\Phi$ is very healthy, then it is healthy.

Suppose $\Phi\colon K\to K$ is strongly linear and very healthy, and $\Phi(\fm)\prec \fm$ for all~$\fm\in \fM$. 
Let $P\in C[[t]]$.  For $i\in I$ with $\Phi(C[[\fM_i]])\subseteq C[[\fM_i]]$ we have the strongly linear map
$\Phi|_{C[[\fM_i]]}\colon C[[\fM_i]]\to C[[\fM_i]]$ and thus the strongly linear map $$P(\Phi|_{C[[\fM_i]]})\colon C[[\fM_i]]\to C[[\fM_i]].$$ We let $P(\Phi)\colon K\to K$ denote the common extension of all these $P(\Phi|_{C[[\fM_i]]})$.  Then~$P(\Phi)$ is strongly linear and very healthy.
Corollary~\ref{cor:P(Phi)} goes through for such~$\Phi$.

\medskip
\noindent
By Corollary~\ref{cor:strongly linear inverse}, certain strongly linear field automorphisms of~$\CM$
have a strongly linear inverse. Even with our standing assumption that $(\fM_i)$ is healthy
we do not know whether the inverse of every strongly linear field automorphism of our directed union $K=\bigcup_i C[[\fM_i]]$
is strongly linear. Therefore we define
a field isomorphism~$K\to L$ to be {\it strong}\/ 
if it and its inverse are strongly linear. The composition of strong isomorphisms~$K\to L$
and~$L\to M$ (where $M$ is as in Lemma~\ref{lem:comp of strongly add} above) is a strong isomorphism~${K\to M}$.
Thus the strong automorphisms of the field $K$ form a subgroup $\SigmaAut(K)$ of the group 
$\Aut(K|C)$ of field automorphisms of $K$ fixing $C$ pointwise.
If $\der$ is a derivation on $K$, then
$\SigmaAut_\der(K)$ denotes the subgroup of $\SigmaAut(K)$ consisting of the
strong automorphisms of the {\it differential}\/ field $(K,\der)$. 
Note that for ${\chi\in \operatorname{Hom}(\fM,C^\times)}$, the strong automorphism~$f\mapsto f_{\times\chi}$ of~$\CM$ restricts to a strong automorphism of~$K$.
If $\der$ is a strongly linear derivation on $K$, then this restriction
commutes with $\der$ iff $\der(\fm)_{\times\chi}=\chi(\fm)\der(\fm)$ for
all~$\fm$.
 

\subsection*{A Galois connection}
{\em Assume that $K$ is equipped with a derivation $\der$}. (As is the case for $\T$ in the next section.) 
For a subset $F$ of $K$ we define the subgroup
$$\SigmaAut_\der(K|F)\ :=\ \big\{\sigma\in\SigmaAut_\der(K):\ \text{$\sigma(f)=f$ for all $f\in F$}\big\}$$
of $\SigmaAut_\der(K)$. We have
$$F_1\subseteq F_2\subseteq K \quad\Rightarrow\quad \SigmaAut_\der(K|F_1)\supseteq\SigmaAut_\der(K|F_2),$$
and 
$$\SigmaAut_\der(K|\emptyset)\ =\ \SigmaAut_\der(K|C)\ =\ \SigmaAut_\der(K).$$
Given
a subset $\mathcal G$ of $\Aut_\der(K)$ the subset
$$K^{\mathcal G}\ :=\ \big\{f\in K:\  \text{$\sigma(f)=f$ for all $\sigma\in\mathcal G$} \big\}$$
of $K$ is definably closed in the differential field $K$.
Clearly $K^{\mathcal G}$  is a differential subfield of~$K$.
We call $K^{\mathcal G}$ the
{\it fixed field}\/ of $\mathcal G$. 
We have $$\mathcal G_1\subseteq\mathcal G_2\subseteq\Aut_\der(K) \quad\Rightarrow\quad K^{\mathcal G_1} \supseteq K^{\mathcal G_2}.$$
By Corollary~\ref{cor:strong fixed gp}, the fixed field of a subset of $\SigmaAut_\der(K)$  is a strong differential subfield of $K$. 
For $f\in K$ we set $\SigmaAut_\der(K|f):=\SigmaAut_\der(K|\{f\})$, and for $\sigma$ in~$\Aut_{\der}(K)$ we set $K^\sigma:=K^{\{\sigma\}}$.

\medskip
\noindent
Let now $F$, $F_1$, $F_2$ range over differential subfields of $K$ and $\mathcal G$, $\mathcal G_1$, $\mathcal G_2$ over subgroups of $\SigmaAut_\der(K)$. 
From Kaplansky~\cite{Kaplansky1,Kaplansky2} we borrow the following convenient terminology:
We  call $\c(F):=K^{\SigmaAut_\der(K|F)}$ the {\it closure}\/ of $F$ (in $K$); then
$$F_1\subseteq F_2\ \Rightarrow\ \c(F_1)\subseteq\c(F_2),\quad F\subseteq\c(F), \quad \c(\c(F))=\c(F),$$
where the last equality comes from $\SigmaAut_\der(K|\c(F))=\SigmaAut_\der(K|F)$. 
We say that~$F$ is  {\it closed}\/ if $F = \c(F)$.
Similarly we define $\c(\mathcal G):=\SigmaAut_\der(K|K^{\mathcal G})$ and  we say that~$\mathcal G$ is {\it closed}\/ if
$\mathcal G = \c(\mathcal G)$.  So $\{\id_K\}$ is closed, and 
$$\mathcal G_1\subseteq \mathcal G_2\ \Rightarrow\ \c(\mathcal G_1)\subseteq\c(\mathcal G_2),\quad \mathcal G\subseteq\c(\mathcal G), \quad \c(\c(\mathcal G))=\c(\mathcal G),$$
the last equality coming from $K^{\c(\mathcal G)}=K^{\mathcal G}$.
Also $\c(K^{\mathcal G})=K^{\c(\mathcal G)}$, so each differential subfield $K^{\mathcal G}$ of $K$ is closed, and likewise, each subgroup~$\SigmaAut_\der(K|F)$ of~$\SigmaAut_\der(K)$ is closed.
Thus the map $F\mapsto \SigmaAut_\der(K|F)$ is an order-reversing bijection 
$$\big\{ \text{closed differential subfields of $K$} \big\} \longrightarrow
\big\{ \text{closed subgroups of $\SigmaAut_\der(K)$} \big\},$$ 
with inverse~$\mathcal G\mapsto K^{\mathcal G}$. Of particular interest are the normal subgroups of
$\SigmaAut_\der(K)$:

\begin{lemma}\label{lemgalois}
Let $F$ be a differential subfield of $K$ and $\mathcal G:=\SigmaAut_\der(K|F)$. 
\begin{enumerate}
\item[\rm(i)] If $\sigma\in\SigmaAut_\der(K)$ satisfies $\sigma(F)=F$, then $\sigma^{-1}\mathcal G\sigma=\mathcal G$.
Hence if $\sigma(F)\subseteq F$ for all
$\sigma\in \SigmaAut_\der(K)$, then $\mathcal G$ is a normal subgroup
of $\SigmaAut_\der(K)$, and the group morphism
$$\sigma\mapsto\sigma|_F\colon\SigmaAut_\der(K)\to\Aut_\der(F)$$ has kernel $\mathcal G$.
\item[\rm(ii)] Suppose $F$ is closed. Then for $\sigma\in\SigmaAut_\der(K)$ we have
$\sigma^{-1}\mathcal G\sigma=\mathcal G$ iff $\sigma(F)=F$. In particular, $\mathcal G$ is a normal subgroup of
$\SigmaAut_\der(K)$ iff $\sigma(F)\subseteq F$ for all~$\sigma\in  \SigmaAut_\der(K)$.
\end{enumerate} 
\end{lemma}

\begin{proof}
Part (i) is obvious. For (ii), suppose $F$ is closed, so $F=K^{\mathcal G}$.
Let $\sigma$ in~$\SigmaAut_\der(K)$ be such that $\sigma^{-1}\mathcal G\sigma=\mathcal G$, and let $f\in F$. Then $\sigma(f)\in F=K^{\mathcal G}$, that is, $\tau\sigma(f)=\sigma(f)$ for all $\tau\in\mathcal G$: indeed, for $\tau\in \mathcal G$ we have
$\sigma^{-1}\tau\sigma\in \sigma^{-1}\mathcal G\sigma=\mathcal G$, so~$\sigma^{-1}\tau\sigma(f)=f$ and so $\tau\sigma(f)=\sigma(f)$.
\end{proof}

\noindent
The rest of this subsection is not used later in this paper.
The {\it Galois topology}\/ on~$\mathcal G:=\SigmaAut_\der(K)$  is the coarsest topology on $\mathcal G$ such that
\begin{enumerate}
\item for each $\sigma\in\mathcal G$, left multiplication $\tau\mapsto \sigma\tau$ by $\sigma$
is continuous;
\item for each $f\in K$, the subgroup $\SigmaAut_\der(K|f)$ of $\mathcal G$ is closed.
\end{enumerate}
By
Soundararajan~\cite{Sound}, the   subgroups of $\mathcal G$ that are closed in \marginpar{didn't check yet that topologically closed
subgroups are closed subgroups} the Galois topology  are exactly the  closed subgroups of $\mathcal G$ defined above. 
Is there a group topology on $\mathcal G$ 
such that the topologically closed subgroups of $\mathcal G$ are exactly the   closed subgroups of $\mathcal G$?
Any such topology is hausdorff, but we face another strong constraint:

\begin{lemma}\label{lem:dense subgp}
Suppose $\mathcal G$ is equipped with a group topology such that the closed
subgroups of $\mathcal G$ are exactly the topologically closed subgroups of $\mathcal G$.
Let $\sigma\in\mathcal G$ be such that $K^\sigma=K^{\mathcal G}$. Then the subgroup
$\langle\sigma\rangle$ of $\mathcal G$ generated by
$\sigma$ is dense in $\mathcal G$, so~$\mathcal G$ is abelian.
\end{lemma}
\begin{proof}
We have
$\operatorname{c}(\langle\sigma\rangle)=\SigmaAut_\der(K|K^\sigma)=\SigmaAut_\der(K|K^{\mathcal G})=\mathcal G$, so 
$\langle\sigma\rangle$ is dense in $\mathcal G$.
Since $\mathcal G$ is a hausdorff, this forces $\mathcal G$ to be abelian.
\end{proof}

\noindent
By Lemma~\ref{lem:dense subgp} and Corollary~\ref{cor:R definably closed} below, there is no
group topology on  $\SigmaAut_\der(\T)$ 
such that the closed subgroups of  $\SigmaAut_\der(\T)$
are the topologically closed subgroups.

\section{Transseries}\label{sec:prelims, trans}

\noindent
An {\it exponential ordered field\/} is an ordered field $E$ equipped with a strictly increasing 
group morphism $E\to E^>$ from the ordered additive group of $E$ to the ordered multiplicative group 
of positive elements of~$E$; this group morphism is called the exponentiation of $E$ and  always denoted here 
by~$\exp$;
the inverse of the ordered group isomorphism $\exp\colon E \to \exp(E)$
is   an ordered group isomorphism 
$\log\colon \exp(E) \to E$, and
 for $f\in \exp(E)$, $g\in E$ we put $f^g:=\exp(g\log f)$.
A {\it logarithmic-exponential ordered field\/} is an exponential ordered field~$E$ with $\exp(E)=E^{>}$.
We consider the ordered field~$\R$ of real numbers as a logarithmic-exponential ordered field with exponentiation $r\mapsto \ex^r$. 
We can extend this exponentiation in a natural way to the valuation ring $\R[[t]]$ of the field $\R(\!(t)\!)$ of Laurent series, but $\exp(t^{-1})$ has no natural sense as an element of that field. 
To remedy this, we extend~$\R(\!(t)\!)$ in two steps to
the  logarithmic-exponential ordered field $\T$ of transseries.
We now briefly recall this construction and refer the reader to \cite{DMM2} and \cite[Appendix~A]{ADH} for proofs.

\subsection*{Exponential transseries}
The first step is to construct the exponential ordered field $\T_{\exp}$ of
{\it exponential transseries}\/ (denoted by $\RxE$ in \cite{DMM2}).
As an ordered field, $E:=\T_{\exp}$ is the union of an increasing sequence
$$E_{0}=\R[[G_{0}]]\ \subseteq\ E_{1}=\R[[G_{1}]]\ \subseteq\ \cdots\subseteq\ E_{m}=\R[[G_{m}]]\ \subseteq\cdots$$
of ordered  subfields of $E$, where each inclusion $\R[[G_{m}]]\subseteq\R[[G_{m+1}]]$
comes from a corresponding inclusion $G_{m}\subseteq G_{m+1}$ of multiplicatively written ordered abelian groups, so~$E=\bigcup_m \R[[G_m]]$ is an ordered  subfield of $\R[[G^{\operatorname{E}}]]$ where $G^{\operatorname{E}}=\bigcup_m G_m$.
The ordered group $G_0:=x^\R$ is a multiplicative copy of the ordered abelian (additive) group $\R$,
with isomorphism $r\mapsto x^r\colon \R\to x^\R$. (The field $\R[[G_0]]=\R[[x^{\R}]]$  contains the field of Laurent series: setting $t^k:= x^{-k}$ for $k\in \Z$ we have $\R[[x^{\Z}]]= \R(\!(t)\!)$.) 
 We have the 
$\R$-linear subspace 
$$A_m\ :=\  \big\{ f\in E_m:\,   \supp f\succ G_{m-1}\big\}\qquad (G_{-1}:=\{1\} \text{ by convention})$$
of $E_m=\R[[G_m]]$; so $A_0=(E_0)_\succ$.
The exponentiation $\exp$ of $E$ maps $E_m$ into~$E_{m+1}$; in more detail, it maps the ordered subgroup $A_m$ of $E_m$ onto the ordered subgroup $\exp(A_m)$ of $G_{m+1}$, and indeed   
  $G_{m+1}=G_m\exp(A_m)$ where $G_m$ is convex in~$G_{m+1}$ and $G_m\cap\exp(A_m)=\{1\}$. 
By  \cite[Lem\-ma~A.1]{ADH} we have 
\begin{align*} (E_m)_\succ\ &=\ A_0\oplus\cdots\oplus A_m,\quad
G_{m+1}\ =\ x^\R \exp(A_0\oplus\cdots\oplus A_m), \text{ where}\\ 
x^\R &\text{ is convex in }G_{m+1},\quad x^\R\cap \exp(A_0\oplus\cdots\oplus A_m)\ =\ \{1\}.
\end{align*} 
 In particular, 
 $E_\succ=\bigoplus_m A_m$ (internal direct sum of $\R$-linear subspaces) and $G^{\operatorname{E}}=x^\R\exp(E_\succ)$ where
$x^\R$ is convex in $G^{\operatorname{E}}$ and $x^\R\cap\exp(E_\succ)=\{1\}$. Also
$$E^{>}\ =\ x^{\R}\cdot \exp(E),\quad \exp(r+\varepsilon)\ =\    \ex^r\sum_n \frac{1}{n!} \varepsilon^n\in E_{m}\quad \text {for }r\in \R,\  \varepsilon\in (E_m)_\prec).$$
In particular, as an exponential ordered field, $\T_{\exp}$  extends the exponential ordered field $\R$. 
We recursively define $\ex_0:=x$ and $\ex_{m+1}:=\exp\ex_m$. Induction then yields
$$\ex_{m}\in A_m^{>},\qquad   \ex_m\in G_m,\quad \ex_m>G_{m-1}.$$

\subsection*{From exponential transseries to general  transseries} Although $x:=x^1>0$, it is easy to see that $x\notin\exp(E)$. This defect is remedied in the second step, in which the exponential ordered field
$E=\T_{\exp}$ is extended to the
logarithmic-exponential ordered field $\T$. The idea is to use
distinct symbols $\ell_0,\ \ell_1,\ \ell_2,\ \dots$ in the role of~$x,\ \log x,\ \log \log x,\ \dots$.
Replacing for any given $n$ the formal variable~$x$ in $\RxE:=E$ by $\ell_n$ changes informally any series~$f(x)\in \RxE$ into
a series~$f(\ell_n)\in \RlognE$. This is done so that each $\RlognE$ is an exponential ordered field with~$\RlognplusE$ extending~$\RlognE$. Then $\T$ is defined to be the union of the~$\RlognE$. The details of the construction are in \cite[Appendix~A]{ADH}. 
Here we just record the facts about $\T$ from \cite[Appendix~A]{ADH} that we shall need. 

\medskip
\noindent
Indeed, $\T$ is a logarithmic-exponential ordered field extension of the exponential ordered field $E=\T_{\exp}$. As an ordered field, $\T$ is an ordered subfield of a Hahn field~$\R[[G^{\operatorname{LE}}]]$ where $G^{\operatorname{LE}}$ is a multiplicatively written ordered abelian group containing $G^{\operatorname{E}}$ as an ordered subgroup; the elements of $G^{\operatorname{LE}}$ are called {\em transmonomials}. The corresponding inclusion map 
$\R[[G^{\operatorname{E}}]]\to \R[[G^{\operatorname{LE}}]]$ maps $\T_{\exp}$ into $\T$: $\T_{\exp}\subseteq \T$. Throughout we view $\R[[G^{\operatorname{LE}}]]$ as our ambient ordered Hahn field. 

\medskip
\noindent
The ordered group $G^{\operatorname{LE}}$ comes equipped with an automorphism $\fm \mapsto \fm{\downarrow}$. The unique extension of this automorphism to a strongly linear automorphism~${f\mapsto f{\downarrow}}$ of the Hahn field $\R[[G^{\operatorname{LE}}]]$ is called the {\em downward shift}. The $n$th iterate of the downward shift operator $\downarrow$ is denoted by $\downarrow_n$. The inverse of $\downarrow$, the {\em upward shift}, is denoted by $\uparrow$ and its $n$th iterate by $\uparrow^n$; they are also strongly linear. We have~${G^{\operatorname{E}}\subseteq G^{\operatorname{E}}{\downarrow}}$, hence an increasing sequence 
$$G^{\operatorname{E}}= G^{\operatorname{E}}{\downarrow_0}\ \subseteq\ G^{\operatorname{E}}{\downarrow_1}\ \subseteq\  G^{\operatorname{E}}{\downarrow_2}\ \subseteq\ \cdots$$
of ordered subgroups of $G^{\operatorname{LE}}$, and in fact, $\bigcup_n G^{\operatorname{E}}{\downarrow_n}=G^{\operatorname{LE}}$.  The downward shift restricts to an automorphism of the {\em exponential ordered field\/} $\T$, and $\T$
contains elements $\ell_0,\ \ell_1,\ \ell_2,\ \dots>\R$  with $\ell_0=x$ and $\ell_{n+1}=\ell_n{\downarrow}=\log \ell_n$ for all $n$, in particular, $x{\downarrow}=\ell_1=\log x$. Also, $x^r=\exp(r\log x)$ for $r\in \R$, justifying our earlier choice of power notation. 
With ``$\subseteq$'' in the model-theoretic sense of ``substructure of'', we have $E\subseteq E{\downarrow}$, and hence an increasing sequence
$$E = E{\downarrow_0}\  \subseteq\  E{\downarrow_1}\ \subseteq \cdots \subseteq\ E{\downarrow_n}\ \subseteq\cdots$$
 of exponential ordered subfields of $\T$. In fact,  $\T=\bigcup_n E{\downarrow_n}$. We obtain the $E{\downarrow_n}$ and $\T$ as directed unions of Hahn fields as follows.  Set 
 $$G_{m,n}\ :=\ G_m{\downarrow_n}, \qquad E_{m,n}\ :=\ E_m{\downarrow_n}\ =\ \R[[G_m]]{\downarrow_n}\ =\ \R[[G_{m,n}]].$$
  Recalling that $G^{\operatorname{E}}= \bigcup_m G_m$  with $(G_m)$ an increasing sequence of convex subgroups of the ordered group $G^{\operatorname{E}}$, we have $G^{\operatorname{E}}{\downarrow_n}=\bigcup_m G_{m,n}$ with $(G_{m,n})_m$ an increasing sequence of convex subgroups of the ordered group $G^{\operatorname{E}}{\downarrow_n}$. Likewise, $E=\bigcup_m E_m$ yields the ordered field $E{\downarrow_n}=\bigcup_m E_{m,n}=\bigcup_m \R[[G_{m,n}]]$ as an ordered subfield of the ordered Hahn field $\R[[G^{\operatorname{E}}{\downarrow_n}]]$.
 Hence $G^{\operatorname{LE}}=\bigcup_n G^{\operatorname{E}}{\downarrow_n}=\bigcup_{m,n} G_{m,n}$. Also~$G_{m,n}\subseteq G_{m+1,n+1}$ (not obvious) and $G_{m,n}\subseteq G_{m+1,n}$ for all~$m$,~$n$, so $(G_{m,n}){}_{(m,n)}$ is a (healthy) directed family of ordered subgroups of $G^{\operatorname{LE}}$. Thus,
 $$\T\ =\ \bigcup_n E{\downarrow_n}\ =\ \bigcup_{m,n} \R[[G_{m,n}]],$$
 so $G^{\operatorname{LE}}\subseteq \T^{>}$, and the notions {\em strongly additive map $\T\to \T$}, {\em strong subgroup of~$\T$}, and the like from the subsection on directed unions in Section~\ref{sec:prelims, well-based} are available. Towards the end of this section we show that $E=\T_{\exp}$ is a strong subgroup of $\T$. 

\medskip
\noindent
Set $A_{m,n}:=A_m{\downarrow_n}$. Then $(E_{m,n})_\succ=A_{0,n}\oplus\cdots\oplus A_{m,n}$, and using $x{\downarrow_n}=\ell_n$:
$$G_{0,n}\ =\ \ell_n^\R,\qquad  G_{m+1,n}\ =\ \ell_n^\R\cdot\exp\big((E_{m,n})_\succ\big),\qquad\ell_n^\R\cap \exp\big((E_{m,n})_\succ\big)\ =\ \{1\}.$$
Taking the union over all $m$ while fixing $n$ then yields
\begin{equation}\label{eq:GEn}
G^{\operatorname{E}}{\downarrow_n}\ =\ \ell_n^\R\exp\big((E{\downarrow_n})_\succ\big), \qquad \ell_n^\R\cap\exp\big((E{\downarrow_n})_\succ\big)\ =\ \{1\}.
\end{equation}
Hence $\log(G^{\operatorname{E}}{\downarrow_n})=\R\ell_{n+1}+(E{\downarrow_n})_\succ\subseteq (E{\downarrow_{n+1}})_\succ$, and thus $\log(G^{\operatorname{LE}})=\T_{\succ}$. Note that $\ex^g=\exp(g)$ for $g\in \T$ in view of 
$f^g:=\exp(g\log f)$ for $f\in \T^{>}$.  An easy induction gives $\ex_n{\downarrow_n}=x$, so $\ex_n=x{\uparrow^n}$. The strictly increasing sequence
$(\ex_n)$ is cofinal in $\T$, and the strictly decreasing sequence $(\ell_n)$ is coinitial in $\T^{>\R}$. 
It follows from $\log(G^{\operatorname{LE}})=\T_{\succ}$ that $\fm^r\in G^{\operatorname{LE}}$ for $\fm\in G^{\operatorname{LE}}$ and $r\in \R$. Thus
$$(r,\fm)\mapsto \fm^r\colon \R\times G^{\operatorname{LE}}\to G^{\operatorname{LE}}$$ makes $G^{\operatorname{LE}}$ a (multiplicatively written) ordered vector space over $\R$, and so as in the subsection ``Powers of well-based series'' in Section~\ref{sec:prelims, well-based}, this map extends to $$(r,f)\mapsto f^r\colon \R\times \R[[G^{\operatorname{LE}}]]^{>}\to \R[[G^{\operatorname{LE}}]]^{>};$$ 
for~$f\in \T^{>}$, $r\in \R$, this $f^r$ agrees with the already defined $f^r:=\exp(r\log f)$. 
Note that for~$r_0,\dots, r_n, s_0, \dots, s_n\in \R$ and the lexicographic ordering on $\R^{1+n}$ we have
$$\ell_0^{r_0}\cdots \ell_n^{r_n}\ \prec\ \ell_0^{s_0}\cdots \ell_n^{s_n}\ \Longleftrightarrow\ (r_0,\dots, r_n)\ <\  (s_0,\dots, s_n). $$
Going back to the $G_m$, induction on $m$ gives  $\ex_0,\dots, \ex_m\in G_m$ and $\fm^\R\subseteq G_m$ for~$\fm\in G_m$.
Applying ${\downarrow}_n$ this yields $\ell_{n}^{\R},\dots, \ell_{n-m}^\R\subseteq G_{m,n}$ for $n\ge m$ (to be used in the proof of Lemma~\ref{lem:intersect with Tlog}). Since $G_{m,n}$ is divisible, $E_{m,n}$ is real closed.

\medskip\noindent
We recall from \cite{DMM2} that $\T$ as an exponential ordered field is an elementary extension of
the ordered field of real numbers with its usual exponentiation $t\mapsto \ex^t$. Various familiar inequalities involving real exponentiation therefore also hold in $\T$, and that is how we sometimes tacitly use this elementary equivalence.

\subsection*{Composition of transseries} We begin with the following observation:

\begin{lemma}\label{lematmostone}  Given any $g\in \T$ there is at most one strongly linear endomorphism of the exponential ordered field $\T$ that maps $x$ to $g$. 
\end{lemma}
\begin{proof} Any strongly linear endomorphism $\eta$ of the exponential ordered field $\T$ is the identity on $\R$ and satisfies for $f\in \T$:
$\eta(\ex^f)=\ex^{\eta(f)}$, and when $f>0$, also~$\eta(\log f)=\log \eta(f)$ (and thus $\eta(x^r)=\eta(x)^r$ for $r\in \R$).
 Now use that by the very construction of $\T$, any logarithmic-exponential ordered subfield of $\T$ that contains $\R\cup \{x^\R\}$ and is closed under sums of summable families must be all of $\T$. 
 \end{proof}
 
 \noindent
 Since $\T$ is a field, any endomorphism $\eta$ as in this lemma is of course also injective, and satisfies $\eta(x)>\R$. (In fact, such $\eta$ is even bijective, as we shall see.) 

\medskip\noindent
{\em In the rest of this subsection $f$ ranges over $\T$ and $g$,~$h$ over $\T^{>\R}$}. There is a unique operation
$$ (f,g)\mapsto f\circ g\ :\ \T\times \T^{>\R}\to \T$$ 
such that for each $g$, the map 
$f\mapsto f\circ g\colon \T\to \T$ is a strongly linear endomorphism of the exponential ordered field $\T$ with  $x\circ g=g$.
The ``at most one'' part of this fact follows from Lemma~\ref{lematmostone}; for existence and more details, see  \cite[Proposition~A.7]{ADH} and the discussion following it.  Here are some properties of this operation: 
\begin{enumerate}
\item[(i)] $f\circ x=f$; 
\item[(ii)] $f{\uparrow}=f\circ \ex^x$, $f{\downarrow}=f\circ \ell_1$;
\item[(iii)] $\ex^x\circ g= \exp g$, $\ell_1\circ g= \log g$, $x^r\circ g=g^r$ ($r\in\R$); 
\item[(iv)] $g\circ h\in \T^{>\R}$ and $(f\circ g)\circ h= f\circ(g\circ h)$.
\end{enumerate}
In particular,
$$f{\uparrow^n}=f\circ \ex_n, \quad f{\downarrow_n}=f\circ \ell_n, \quad
(f\circ g){\uparrow^n} = f\circ (g{\uparrow^n}),\quad
(f\circ g){\downarrow_n} = f\circ (g{\downarrow_n}).$$
We often use a more general version of (iii):

\begin{lemma}\label{lem:(iii) general} We have $\ex^f\circ g= \ex^{f\circ g}$, and if $f>0$, then 
$$(\log f)\circ g\ =\ \log(f\circ g), \qquad f^r\circ g\ =\ (f\circ g)^r \quad (r\in \R).$$ 
\end{lemma} 
\begin{proof} If $f>\R$, then $\ex^f\circ g=(\ex^x\circ f)\circ g=\ex^x\circ (f\circ g)=\ex^{f\circ g}$ by (iii) and (iv), and likewise we obtain
$(\log f)\circ g= \log(f\circ g)$.  For exponentiation the case $f<\R$ reduces to the case $f>\R$ by considering $-f$ instead.  If $f\preceq 1$, then $f=c+\varepsilon$ with $c\in \R$ and~$\varepsilon\prec 1$, so we can use
 $\ex^f=\ex^c\sum_n \varepsilon^n/n!$ and strong linearity to obtain ${\ex^f}\circ g= \ex^{f\circ g}$.  For taking logarithms, the case
 $0<f\prec 1$ reduces to the case $f>\R$ by considering $1/f$ instead. 
  If $f=c(1+\varepsilon)$, $c\in \R^{>}$, $\varepsilon\prec 1$, use
$\log(1+\varepsilon)=\sum_{m=1}^\infty (-1)^{m-1} \varepsilon^m/m$ and strong linearity.  The case of powers $f^r$ ($f>0$, $r\in \R$) reduces to exponentiation and taking logarithms by $f^r=\ex^{r\log f}$. 
\end{proof} 

\noindent
By \cite[Proposition~A.8]{ADH} we have:

\begin{prop}\label{prop:comp}  If
$f\in E_{m,n}$ and $g\in E_{p,q}$ \textup{(}$p,q\in \N$\textup{)}, then
$f\circ g\in E_{m+n+p,n+q}$. Hence given $g$, 
the map $f\mapsto f\circ g$ is healthy, and if $f\in E$, $g\in E{\downarrow_n}$ then $f\circ g\in E{\downarrow_n}$.
\end{prop}

\noindent
The binary operation $(g,h)\mapsto g\circ h$ on $\T^{>\R}$ is a group operation and
 makes $\T^{>\R}$ a group with identity element $x$; see~\cite[Proposition~6.3]{DMM2}. 
 Let $g^{\inv}$ denote the inverse of~$g$ in this group.

\begin{cor}\label{cor:comp}
Suppose $g\in E{\downarrow_n}$. Then the embedding $f\mapsto f\circ g\colon E\to E{\downarrow_n}$  
of   exponential ordered fields is surjective
 iff $\ell_n\circ g^{\inv}\in E$.
\end{cor}

 \noindent 
In the proof of Lemma~\ref{faithful} we shall use the following elementary estimate:

\begin{lemma}\label{lemest}  Let $g>h>\R$ in $\T$. Then for $n\ge 1$ we have
$$\frac{{\ex_n}\circ g}{{\ex_n}\circ h}\  \ge\  1+ \Big(\prod_{i=1}^{n-1}{\ex_i}\circ h\Big)\cdot (g-h).$$
\end{lemma} 
\begin{proof} By induction on $n$. Set $\phi:= g-h$. The case $n=1$ follows from $\ex^\phi\ge 1+\phi$. Suppose the inequality holds for a certain $n\ge 1$. Then 
\begin{align*} {\ex_n}\circ g-{\ex_n}\circ h\ &\ge\ ({\ex_n}\circ h) \cdot\Big(\prod_{i=1}^{n-1}{\ex_i}\circ h\Big)\cdot \phi\ =\ 
 \Big(\prod_{i=1}^{n}{\ex_i}\circ h\Big)\cdot \phi,\ \text{ so}\\
\frac{{\ex_{n+1}}\circ g}{{\ex_{n+1}}\circ h}\ &=\  \ex^{{\ex_n} \circ g -{\ex_n}\circ h}\ \ge\  1+ \Big(\prod_{i=1}^{n}{\ex_i}\circ h\Big)\cdot \phi.\qedhere
\end{align*}
\end{proof} 

\noindent
The way this lemma gets used is that no matter how small $g-h\in \T^{>}$, the quotient~${{\ex_n}\circ g}/{{\ex_n}\circ h}$ gets arbitrarily large in $\T$ as $n\to \infty$. 

\medskip
\noindent
The next result is a special case of \cite[Proposition 4.9]{Edgar}. We shall use it to prove
Corollary 2.10, which in turn is needed to obtain Proposition 7.5.

\begin{prop}[Edgar]\label{prop:Edgar}
Let $f\in \T_{\succ}^{>}$. Then the map 
\begin{equation}\label{eq:Edgar}
g\mapsto f\circ g\colon \T^{>\R}\to\T^{>\R}
\end{equation}
is strictly increasing.
\end{prop}

\noindent
We include a proof of this proposition following  \cite[proof of Proposition~3.1]{Mantova}. We first define when $f\in\T$ has {\it exponential rank at most $k\in\N$}\/  (in
symbols:  $\operatorname{ER}(f)\le k$) as follows:
\begin{enumerate}
\item $\operatorname{ER}(f)\le 0$ iff $f\in\R\ell_n$ for some $n$,   and
\item
$\operatorname{ER}(f)\le k+1$ iff
$\operatorname{ER}(\log \fm)\le k$ for all $\fm\in\supp f$.
\end{enumerate}
Let $k\in \N$. It is clear than  $\big\{f\in \T:\, \operatorname{ER}(f)\le k+1\big\}$ is an $\R$-linear subspace of $\T$, and induction on $k$ shows for $f\in \T$: $\operatorname{ER}(f)\le k\ \Rightarrow\ \operatorname{ER}(f)\le k+1$. 
Using these facts we show that for each $f\in\T$ there exists $k\in\N$ such that $\operatorname{ER}(f)\le k$; more precisely:

\begin{lemma}
If $f\in E_{m,n}$, then $\operatorname{ER}(f)\le m+1$.
\end{lemma}
\begin{proof}
By induction on $m$. Let $f\in E_{0,n}$, so $f=\sum_r f_r \ell_n^r$ ($f_r\in\R$).
Since~$\log \ell_n^r=r\ell_{n+1}$ has exponential rank~$\le 0$, we  have $\operatorname{ER}(f)\le 1$.
Now suppose the lemma holds for a certain $m$ and $f\in E_{m+1,n}$, so $f=\sum_\fm f_\fm\fm$ where
$\fm\in G_{m+1,n}$. Then~$\log\fm \in\R\ell_{n+1}+(E_{m,n})_\succ$,
hence $\operatorname{ER}(\log\fm)\le m+1$ by the inductive hypothesis. Thus~$\operatorname{ER}(f)\le m+2$.
\end{proof}

\begin{proof}[Proof of Proposition~\ref{prop:Edgar}]  
By induction on $k\in\N$ we show:
if $f\in\T_{\succ}^{>}$, $\operatorname{ER}(f)\le k$, then the map \eqref{eq:Edgar} is strictly increasing.
This is clear for~$k=0$, since the restriction of $\log$ to a map $\T^{>\R}\to\T^{>\R}$ is strictly increasing,
and $\ell_n\circ g=\log^n( g )$ for all $g\in\T^{>\R}$ and all $n$; see Lemma~\ref{lem:(iii) general}.
 
 Suppose the claim holds for a certain $k$, and 
 $f\in\T_{\succ}^{>}$, $\operatorname{ER}(\log \fm)\le k$ for all~${\fm\in\supp f}$,
 and
 $g<h$ in $\T^{>\R}$; it is enough to show that then $f\circ g<f\circ h$.
 Let $f=\sum_{\fm} f_{\fm}\fm$ ($f_{\fm}\in\R$).
 Set $\fd:=\fd(f)$; so $f_{\fd}>0$ since~$f>0$, and replacing $f$ by $f_{\fd}^{-1}f$ we arrange $f_{\fd}=1$. 
 For $\fm\in\supp f$ we have $\log\fm\in\T_\succ$, and
 from $f\in\T_\succ$  we  obtain $\fm\succ 1$
 and thus $\log\fm>0$,
 hence by inductive hypothesis~$(\log \fm)\circ g<(\log \fm)\circ h$, and 
 using that $\exp$ is strictly increasing:
 $\fm\circ g<\fm\circ h$.
In particular, $\fd\circ g-\fd\circ h<0$;
thus it suffices to prove:
 
 \claim{$\fd \circ g - \fd\circ h \sim f\circ g-f\circ h$.}
 
 \noindent
 Let $\fm\in (\supp f)\setminus\{\fd\}$, so $\fm\prec\fd$. Then $\log(\fd/\fm)\in\T_\succ$, $\log (\fd/\fm) > 0$,  
 so  by the inductive hypothesis again, 
 $(\log (\fd/\fm)) \circ g< (\log (\fd/\fm))\circ h$ 
 and hence  $(\fd/\fm) \circ g < (\fd/\fm) \circ h$. This yields
 $$(\fd\circ g)/(\fd\circ h) - 1 < (\fm\circ g)/(\fm\circ h) - 1<0$$
 and hence
 $$(\fd\circ g)/(\fd\circ h) - 1 \succeq (\fm\circ g)/(\fm\circ h) - 1.$$
 From $0<\fm\prec \fd$ we get  
 $\fm\circ h\prec \fd\circ h$, therefore
 \begin{align*}
 \fd\circ g-\fd\circ h &= (\fd\circ h)\cdot \big((\fd\circ g)/(\fd\circ h) - 1\big)  \\ &\succ
 (\fm\circ h) \cdot \big((\fm\circ g)/(\fm\circ h) - 1\big) = \fm\circ g-\fm\circ h.
 \end{align*}
 This yields the claim and   finishes the inductive step.
 \end{proof}

\noindent
Just like Edgar~\cite[proof of 4.10 from 4.9]{Edgar} we deduce from Proposition~\ref{prop:Edgar}:

\begin{cor}\label{cor:Edgar, 1}
Let   $f_1\prec f_2$ in $\T_\succ$ and $g<h$ in $\T^{>\R}$.
Then $$f_1\circ h - f_1\circ g \prec f_2\circ h - f_2\circ g.$$
\end{cor}
\begin{proof}
Replace $f_2$ by $-f_2$ if necessary to arrange $f_2>0$.
Let $c\in\R$ and put~$f:=f_2-cf_1\in\T_\succ$; then $f>0$. Thus
$f\circ g<f\circ h$  by Proposition~\ref{prop:Edgar}, so~$c(f_1\circ h - f_1\circ g ) < f_2\circ h - f_2\circ g$.
Now use that this holds for all $c\in\R$.
\end{proof}

\begin{cor}\label{cor:Edgar, 2}
Let $c\in\R^>$ and $f\in\T_\succ$ with $f\succ x$. Then $f\circ (x+c)-f \succ 1$.
\end{cor}
\begin{proof}
Take $f_1:=x$, $f_2:=f$, $g:=x$, $h:=x+c$ in the previous corollary.
\end{proof}

\begin{cor}\label{cor:steep fm}
Let $c\in\R^>$ and $g\in\T^>$ with $\log g\succ x$.
Then
$$g\succ 1\ \Rightarrow\ g\circ(x+c) \succ g,\qquad g\prec 1\ \Rightarrow\ g\circ(x+c) \prec g.$$
\end{cor}
\begin{proof}
Set  $\fm=\fd(g)$. Then $g \asymp \fm$ and thus $g\circ(x+c) \asymp \fm\circ(x+c)$. As~$g\nasymp 1$, we have~$\fm\ne 1$, so
$\log \fm \in \T_{\succ}^{\ne}$, and thus $\log g\asymp\log\fm$. Thus
 we may arrange~${g=\fm\in G^{\operatorname{LE}}}$.
Then $f:=\log\fm\in \T_\succ$. If~$\fm \succ 1$, then $f>0$ and therefore~${f\circ (x+c) - f > \R}$ by Proposition~\ref{prop:Edgar} and Corollary~\ref{cor:Edgar, 2}, 
 so $\fm\circ(x+c)=\ex^{f\circ(x+c)}\succ \ex^f=\fm$.
 If $\fm\prec 1$, pass from $\fm$ to~$\fm^{-1}$.
\end{proof}

\subsection*{The derivation of $\T$}
There is a unique strongly linear derivation~$\der$ on~$\T$ such that  $\der(x)=1$ and 
$\der(\exp f)=\der(f)\exp f$ for   $f\in \T$;
see \cite[Proposition~A.5]{ADH}.

\medskip
\noindent
{\it Below we view $\T$ as a differential field with this derivation, and  $\T$
 denotes the \emph{differential}\/ field of transseries, 
unless indicated otherwise.}\/

\medskip
\noindent
We write $f'=\der(f)$, $f''=\der^2(f)$, etc., for $f\in\T$.
The constant field of~$\der$ is~$C_{\mathbb T}=\{f\in\T:f'=0\}=\R$, see~\cite[Theorem~3.9]{DMM2}. 
For all $f\in\T$ and all $n$ we have 
$$(f{\uparrow^n})'\ =\ (f'{\uparrow^n})\ex_1\cdots\ex_n\quad\text{and}\qquad
(f{\downarrow_n})'\ =\ (f'{\downarrow_n})\frac{1}{\ell_0\cdots\ell_{n-1}},$$
in particular
$$\ex_n' \ =\  \ex_1\cdots\ex_n\quad\text{and}\quad \ell_n' \ =\  \frac{1}{\ell_0\cdots\ell_{n-1}}.$$
The Chain Rule holds: if $f\in\T$, $g\in \T^{>\R}$ then 
$$(f\circ g)'\  =\  (f'\circ g)\cdot g';$$ 
see~\cite[Proposition~6.3]{DMM2}.
For each $m$, $E_{m}$ is a differential subfield of $\T$, as is~$E_{m,n}$ for~$m\geq n$;
so $E{\downarrow_n}$ is a differential subfield
of $\T$, and   $\der\colon\T\to\T$ is very healthy. 

By \cite[Section 5]{DMM2} we have $\der(\T)=\T$. We shall need the following refinement:

\begin{lemma} \label{lemint} The map $\int\colon \T\to \T$ that assigns to each $g\in \T$ the unique $f\in \T$ with $f'=g$ and $f_1=0$ is strongly linear.   
\end{lemma} 
\begin{proof} 
By \cite[5.3]{DMM2} we have $\der(E_m)=\der\big(\R[[G_m]]\big)=\der\big(\R[[G_m^{\ne}]]\big)=\R[[G_m\setminus \{x^{-1}\}]]$. Let $m\ge n$ and recall that $E_{m,n}=E_m{\downarrow_n}= \R[[G_{m,n}]]$. Applying the formula above for $(f{\downarrow_n})' $ to $f\in E_m$ gives
\begin{align*} \der(E_{m,n})\ &=\ \der\big(\R[[G_{m,n}^{\ne}\big)\ =\ \textstyle\frac{1}{\ell_0\cdots\ell_{n-1}}\left(\R[[G_m\setminus \{x^{-1}\}]]{\downarrow_n}\right)\\
&=\ \R\big[\big[G_{m,n}\setminus \big\{\textstyle\frac{1}{\ell_0\cdots \ell_n}\big\}\big]\big]\ =\ \R[[G_{m,n}\setminus\{\ell_n^\dagger\}]].
\end{align*} 
Moreover, $\der$ is injective on $\R[[G_{m,n}^{\ne}]]$, and $\der(\fm)\prec \der(\fn)$ whenever $\fm\prec \fn$ in $G_{m,n}^{\ne}$.
Thus by Corollary~\ref{cor:extend to strongly linear} the restriction of $\der$ to a bijective and strongly linear map $\R[[G_{m,n}^{\ne}]]\to \R[[G_{m,n}\setminus\{\ell_n^\dagger\}]]$ has a strongly linear inverse $$g\mapsto \textstyle\int g\ :\  \R[[G_{m,n}\setminus\{\ell_n^\dagger\}]]\to \R[[G_{m,n}^{\ne}]].$$ 
If $n\ge 1$, then $G_{m-1,n-1}\subseteq G_{m,n}\setminus\{\ell_n^\dagger\}$, so $g\mapsto \int g\colon \R[[G_{m-1,n-1}]]\to \R[[G_{m,n}]]$ is strongly linear. Thus $\int\colon \T\to \T$ is strongly linear. 
\end{proof}

\begin{notation}
Let $f, g \in \T^\times$.
We put
\begin{eqnarray*}
  f \flattereq g & : \Longleftrightarrow & f^{\dag} \preccurlyeq
  g^{\dag},\\
  f \flatter g & : \Longleftrightarrow & f^{\dag} \prec g^{\dag},\\
  f \comp g & : \Longleftrightarrow & f^{\dag} \asymp g^{\dag}.
\end{eqnarray*}
By
\cite[Proposition~7.3]{AD1}, if $f,g\succ 1$:
\begin{eqnarray*}
  f \flattereq g \quad & \Longleftrightarrow & \quad |f| \le |g|^n 
 \text{ for some $n\ge 1$},\\
  f \flatter g \quad & \Longleftrightarrow & \quad
|f|^n < |g| \text{ for all $n\ge 1$}.
\end{eqnarray*}
Note that $f\flatter\exp f$ for $f>\R$, so  $\ell_{n+1}\flatter\ell_n$ and $\ex_n\flatter\ex_{n+1}$. We have
$$1\prec f\preceq g\ \Rightarrow f\flattereq g, \qquad f\preceq g\prec 1\ \Rightarrow f\steepereq g,$$
and hence
$$f\asymp g\nasymp 1\ \Rightarrow\ f\comp g.$$
We also note:
$$f\flattereq g\ \Rightarrow\ fg \flattereq g,\qquad
f\flatter g\ \Rightarrow\ fg\comp g.$$
\end{notation}

\subsection*{Exponentiation and differentiation}
Let $K$ be an exponential ordered field. We say that a derivation $\der$ on $K$ is {\it compatible}\/ (tacitly: with exponentiation)
if~$\der(\exp f)=\der(f)\exp f$ for $f\in K$. 
Thus the above derivation  on the logarithmic-exponential ordered field $\T$ is compatible.
Let $\der$ be a compatible derivation on $K$.
Then $C_K$ is an exponential ordered subfield of $K$.
Note that   $\der(\log g)=\der(g)/g$ for~$g\in \exp(K)$;
in particular, if
$K$ is a logarithmic-exponential ordered field, then~$C_K$ is also a logarithmic-exponential ordered field.

\begin{lemma}\label{lem:sigma alpha}
Let $K$, $L$ be exponential ordered fields  equipped with compatible derivations, and let $\sigma\colon K\to L$ be  an  embedding of differential fields with~$\sigma(\exp K)\subseteq\exp L$. 
For $f\in K$, set
$$\alpha(f)\ :=\ \log\big(\sigma(\exp f)/\exp \sigma(f)\big)\  =\  \log\sigma(\exp f) - \sigma(f)\in L,$$ 
so $\sigma(\exp f)=\exp\big(\alpha(f)+\sigma(f)\big)$.
Then:
\begin{enumerate}
\item[\rm(i)] $\alpha(f)\in C_L$ for all $f\in  K$, and $\alpha\colon K\to C_L$ is a  morphism of additive groups with kernel 
$$\ker\alpha\ =\ \big\{f\in K:\ \sigma(\exp f)=\exp {\sigma(f)}\big\};$$
\item[\rm(ii)] for $f\in \exp(K)$, $g\in K$ we have $\sigma(\log f)=\log\sigma(f)-\alpha(\log f)$ and
$$ 
\sigma(f^g)\ =\ \exp\big(\alpha(g\log f)-\sigma(g)\alpha(\log f)\big)\cdot \sigma(f)^{\sigma(g)}.$$
\end{enumerate}
\end{lemma}
\begin{proof}
Let $f\in K$. Then 
$$\big(\!\log \sigma(\exp f)\big)' = \frac{\big(\sigma(\exp f)\big)'}{\sigma(\exp f)} = 
\frac{\sigma\big((\exp f)'\big)}{\sigma(\exp f)} =
 \frac{\sigma(f'\exp f)}{\sigma(\exp f)} = \sigma(f') = \big(\sigma(f)\big)'$$
and hence $\alpha(f)\in C_F$. 
The rest follows easily from the definitions.
\end{proof}


\begin{cor}
Let $\sigma$ be a differential field endomorphism of $\T$. Then: 
\begin{enumerate}
\item[\rm(i)]  $\sigma(\T_\succeq)\subseteq \T_\succeq\ \Longleftrightarrow\
\sigma(\R^>\cdot G^{\operatorname{LE}})\subseteq \R^>\cdot G^{\operatorname{LE}}$.
\item[\rm(ii)] $\sigma(\T_\succeq)=\T_\succeq\ \Longleftrightarrow\
\sigma(\R^>\cdot G^{\operatorname{LE}})=\R^>\cdot G^{\operatorname{LE}}$.
\end{enumerate} 
\end{cor}
\begin{proof} Since $\T$ is real closed, $\sigma$ is also an endomorphism of the ordered field $\T$, so~$\sigma(\exp \T)\subseteq \exp \T$. This makes Lemma~\ref{lem:sigma alpha} available. Item (i) is now an easy consequence of the identity $\sigma(\exp f)=\exp\big(\alpha(f)+\sigma(f)\big)$ from that lemma and the fact that $\R^>\cdot G^{\operatorname{LE}}=\exp(\T_{\succeq})$. For (ii), note first that
$\sigma$ is the identity on the constant field $\R$ of $\T$, and so $\R\subseteq \ker (\alpha)$ with $\alpha$ as in that lemma. 
Suppose~${\sigma(\T_\succeq)=\T_\succeq}$. To obtain $\sigma(\R^>\cdot G^{\operatorname{LE}})=\R^>\cdot G^{\operatorname{LE}}$, that is, $\sigma\big(\!\exp(\T_{\succeq})\big)=\exp(\T_{\succeq})$, let~$g\in \T_{\succeq}$; it is enough to find $f\in \T_{\succeq}$ such that $\sigma(\exp f)=\exp g$. Take $f_0\in \T_{\succeq}$ such that~$\sigma(f_0)=g$, and let $f=f_0+r$, $r\in \R$. Then
$\sigma(\exp f)=\exp\big(\alpha(f)+\sigma(f)\big)=\exp g$ if~$\alpha(f)+\sigma(f)=g$, that is, $\alpha(f_0)+r=0$, so $r:=-\alpha(f_0)$ works. A similar argument yields the other direction of (ii).  
\end{proof} 

\subsection*{Newtonianity} In this subsection we use the concepts of {\em $H$-field\/} and  {\em asymptotic couple of an $H$-field}, and deal with various properties that an $H$-field can have: (un)grounded, $\upo$-free, newtonian; see \cite[Chapters~6, 9, 10, 11, 14, 15]{ADH}. By \cite[Corollary~15.0.2]{ADH}, the $H$-field $\T$ is $\upo$-free and newtonian. In Corollary~\ref{cor:newt} below we extend this to a large class of $H$-subfields of $\T$. 
The following is implicit in \cite[proof of Theorem~6.4]{Kaplan}. 

\begin{lemma}\label{lem:newt} Let $F$ be a differential subfield of $\T$ that strictly contains $\R$ $($so $F$ is an $H$-subfield of $\T)$. Suppose $F$ is ungrounded. Then $F$ is $\upo$-free. 
\end{lemma} 
\begin{proof} Take a logarithmic sequence $(l_{\rho})_{\rho < \kappa}$ for $F$; see \cite[Section~11.5]{ADH}. Here $\kappa$ is an infinite limit ordinal. We claim that $\kappa=\omega$ and that $(l_{n})$ is a logarithmic sequence for $\T$. (The lemma follows from this claim since $\T$ is $\upo$-free.) Let $\chi\colon \Gamma_{\T}^{<} \to \Gamma_{\T}^{<}$ be the (increasing) contraction map of the asymptotic couple $(\Gamma_{\T}, \Psi_{\T})$ of $\T$. It restricts to the contraction map of the asymptotic couple $(\Gamma, \psi)$ of $F$, and for~$\gamma:= v(l_0)\in \Gamma^{<}$ we have $v(l_n)=\chi^n(\gamma)$. Take $m$ such that $\ex_m\succeq l_0$ (possible since $\{\ex_n: n\ge 0\}$ is cofinal in $\T$).
Then $\gamma^*:= v(\ex_m)\le \gamma$, so $\chi^n(\gamma^*)\le \chi^n(\gamma)$ for all $n$. Now~$\chi^{m+n}(\gamma^*)=v\big(\ell_n)$ for all $n$, so the sequence $\big(\chi^n(\gamma^*)\big)$ is cofinal in $\Gamma_{\T}^{<}$. Hence the sequence
$\big(\chi^n(\gamma)\big)$ is cofinal in $\Gamma_{\T}^{<}$ and thus in $\Gamma^{<}$. This proves the claim. 
\end{proof} 

\noindent
We need facts about the newtonization of an $\upo$-free $H$-field that were only proved in \cite[Chapter~14]{ADH}  under the assumption of the value group of the $H$-field being divisible. This assumption was removed in \cite{Pynn-Coates}, and we shall use this now.

\begin{cor} \label{cor:newt} Let $F$ be as in the lemma above, and suppose also that $F$ is strong and truncation closed. Then $F$ is newtonian.
\end{cor}
\begin{proof} Since $F$ is $\upo$-free, it has a newtonization $F^{\operatorname{n}}$, an immediate extension of $F$. We embed $F^{\operatorname{n}}$ over $F$ into $\T$ and identify it with its image in $\T$. 
This makes it an immediate extension of $F$ inside $\T$, so $F=F^{\operatorname{n}}$ by Corollary~\ref{corbestapp}. 
\end{proof}



\subsection*{Some strong subgroups of $\T$}
For later use  we record here some facts about intersections of subgroups of $G^{\operatorname{LE}}$. In passing we show that (the additive group of)~$E{\downarrow_n}$  
is a strong subgroup of $\T$.

\begin{lemma}\label{lem:Em-1,n}
$G_{m,n+1}\cap G^{\operatorname{E}}{\downarrow_n}=G_{m-1,n}$ \textup{(}with $G_{-1,n}:=\{1\}$ by convention\textup{)}.
\end{lemma}
\begin{proof}
It suffices to show this for $n=0$, that is, to show 
$G_{m}{\downarrow}\cap G^{\operatorname{E}}= G_{m-1}$;  equivalently (with $E_{-1}:=\R$ by convention): 
\begin{equation}\label{eq:Em-1,n}
E_m{\downarrow}\cap \R[[G^{\operatorname{E}}]]\  =\  E_{m-1}.
\end{equation}
We proceed by induction on $m$.
For $m=0$ we have $$G_{0}{\downarrow}\cap G^{\operatorname{E}}=\ell_1^\R\cap  G^{\operatorname{E}}=\{1\}=G_{-1}.$$
Suppose \eqref{eq:Em-1,n} holds for a certain $m$.
We have  $G_{m+1}{\downarrow} = \ell_{1}^\R\exp{\!\big((E_m{\downarrow})_\succ}\big)$. 
Moreover, $\ell_{1}^\R\cap\exp{\!\big((E{\downarrow})_\succ\big)}=\{1\}$ and 
$$G^{\operatorname{E}}\ =\ \ell_0^\R\exp{\big(E_\succ\big)}\ =\ \exp (\R\ell_{1}+E_\succ)\ \subseteq\ \exp{\big((E{\downarrow})_\succ\big)},$$ 
 hence  $G_{m+1}{\downarrow}\cap G^{\operatorname{E}}\subseteq \exp\big((E_{m}{\downarrow})_\succ\big)$ and thus
\begin{align*}
G_{m+1}{\downarrow}\cap G^{\operatorname{E}}\ 	&=\ \exp\big((E_{m}{\downarrow})_\succ\cap (\R\ell_{1}+E_\succ)\big)\\
													&=\ \exp\big(\R\ell_{1}+ ((E_{m}{\downarrow})_\succ\cap E_\succ)\big)\\
													&=\ \exp\big(\R\ell_{1}+(E_{m-1})_\succ\big)\ =\ G_{m},
\end{align*}
where we used Corollary~\ref{cor*:distributive} for the second equality and \eqref{eq:Em-1,n} for the third.
\end{proof}

\noindent
More generally, with $G_{-m,n}:=\{1\}$ for $m\ge 1$:

\begin{lemma}
For all $k\in\N$ and $m$, $n$  we have $G_{m,n+k}\cap G^{\operatorname{E}}{\downarrow_n}=G_{m-k,n}$.
\end{lemma}
\begin{proof}
By induction on $(k,m)$ where $\N^2$ is ordered lexicographically.
The case~$k=0$ being trivial, suppose the equality in the lemma holds for a certain $k$. Then $$G_{m,n+k+1} \cap G^{\operatorname{E}}{\downarrow_n} \ =\ 
\big(G_{m,n+k+1}\cap G^{\operatorname{E}}{\downarrow_{n+k}}\big) \cap G^{\operatorname{E}}{\downarrow_n} \ =\ G_{m-1,n+k}
\cap  G^{\operatorname{E}}{\downarrow_n}$$
by Lemma~\ref{lem:Em-1,n}. 
If $m=0$, then $G_{m-1,n+k}
\cap  G^{\operatorname{E}}{\downarrow_n}=\{1\}=G_{m-(k+1),n}$, and
if~$m\geq 1$, then $G_{m-1,n+k}
\cap  G^{\operatorname{E}}{\downarrow_n}=G_{m-(k+1),n}$ by inductive hypothesis.
\end{proof}

\begin{cor} For each $n$ the additive 
subgroup  $E{\downarrow_n}$  of $\T$ is strong. 
\end{cor}
\begin{proof}
Let $(f_\lambda)$ be a family of elements of $E{\downarrow_n}$ which is summable in $\T$.
Take~$p\in\N$ and $m$  such that $f_\lambda\in\R[[G_{m,p}]]$ for all $\lambda$
and $(f_\lambda)$ is summable in $\R[[G_{m,p}]]$.
Then~$f_\lambda\in  \R[[G_{m+n-p,n}]]$ for all $\lambda$:
if $p\leq n$, this follows from  $G_{m,p}\subseteq  G_{m+n-p,n}$, and if $p>n$
from the previous lemma applied with $k=p-n$.
Hence $\sum_\lambda f_\lambda \in \R[[G_{m+n-p,n}]] \subseteq E{\downarrow_n}$.
\end{proof}

\noindent
In view of Lemma~\ref{lem:fixsigma, 1} this yields $E{\downarrow_n}=\T\cap \R[[G^{\operatorname{E}}{\downarrow_n}]]$ (inside $\R[[G^{\operatorname{LE}}]]$). Thus~$E=\T_{\exp}$ is strong in $\T$ with $\T_{\exp}=\T\cap \R[[G^{\operatorname{E}}]]$. 
We finish this subsection with a technical fact to be used in the proof of Lemma~\ref{lem:intersect n+1 to n}.  

\begin{cor}\label{cor+:distributive}
$\big(E_\succ + (E_{m-1}{\downarrow})_\succ\big)\cap A_m{\downarrow}\ =\ (E_{m-1})_\succ \cap A_m{\downarrow}\ =\  E_\succ\cap A_m{\downarrow}$.
\end{cor}
\begin{proof}
By Lemma~\ref{lem:Em-1,n}, 
 $(E_m{\downarrow})_\succ \cap E_\succ= (E_{m-1})_\succ$.  Now $(E_m)_{\succ}=(E_{m-1})_{\succ}\oplus A_m$, so
$(E_m{\downarrow})_\succ = (E_{m-1}{\downarrow})_\succ + A_m{\downarrow}$, and thus
$$\big(E_\succ + (E_{m-1}{\downarrow})_\succ\big)\cap A_m{\downarrow}=
\big( (E_{m-1})_\succ + (E_{m-1}{\downarrow})_\succ\big)\cap A_m{\downarrow}.$$
From $(E_{m-1})_{\succ}\cap A_m=\{0\}$ we obtain $(E_{m-1}{\downarrow})_\succ\cap A_m{\downarrow}=\{0\}$.
The $\R$-linear subspaces
$(E_{m-1})_\succ$, $(E_{m-1}{\downarrow})_\succ$, $A_m{\downarrow}$  of $\R[[G^{\operatorname{LE}}]]$ are strong and truncation closed, so we can appeal to Corollary~\ref{cor:distributive}. 
\end{proof}

\subsection*{Decomposing $\T$} 
In this subsection we adapt some material from  Section~4 of \cite{ADHdagap} to our $\T$, which differs from $\T$ in \cite{ADHdagap}.
{\it In the rest of this section we let~$\fm$,~$\fn$ range over~$G:=G^{\operatorname{LE}}$.}\/
Let~$\fN$ be a subgroup of~$G$. Then $\fN$ is 
convex iff for
for all $\mathfrak{m}, \mathfrak{n}$:
\[ \mathfrak{m} \flattereq \mathfrak{n} \in
   \mathfrak{N} \quad\Longrightarrow\quad \mathfrak{m} \in \mathfrak{N} . \]
If $\fN$ is convex, then $\fN$ is also {\it closed under $\R$-powers}\/, that is:
\[ \mathfrak{n} \in
   \mathfrak{N},\ r\in\R \quad\Longrightarrow\quad \mathfrak{n}^r \in \mathfrak{N} . \]
{\it In the rest of this subsection $\mathfrak{N}$ is a convex subgroup of $G$.}\/
Note that if $\fN\ne\{1\}$, then there exists $m$ such that $\ell_n\in\fN$ for all $n\ge m$.
To $\mathfrak{N}$ we associate  
\[ I\  :=\   \{ \mathfrak{m}   :  
\text{$\fm\succ 1$ and $\exp \mathfrak{m} 
\flattereq \frak{n}$ for some $\frak{n}\in\mathfrak{N}$}   \}\ 
\subseteq\ \mathfrak{N}  . \]
Then $I$ is downward closed in  $G^{\succ}:=G^{\succ 1}$ (with
$I=\emptyset$ if $\mathfrak{N}=\{1\}$).
Consequently, the complement $F = G^{\succ} \setminus I$ of $I$
is upward closed in  $G^{\succ}$, and
\[ \    \mathfrak{R}\  :=\  \big\{ \mathfrak{r} \in G : \log \mathfrak{r}
   \in \R [[ F ]] \big\} \]
is also a subgroup of $G$ closed under $\R$-powers.

\begin{lemma} \label{lem:logN}
For all $\frak{m}$ we have
$\frak{m}\in\frak{N}$ iff $\log\frak{m} \in \R[[I]]$.
\end{lemma}
\begin{proof} Let $\fm$ be given.  Then
$$\log \fm \in \T_{\succ}\ \subseteq\  \R[[G^{\succ}]]\ =\ \R[[F]]\oplus \R[[I]].$$
Take $\varphi\in\R[[I]]$, $\psi\in\R[[F]]$ with  $\log\fm=\varphi+\psi$.
 Then $\ex^\varphi,\ex^\psi\in G$. Indeed, $\ex^\varphi\in\fN$, because if $\varphi\ne 0$, then
 $\fd(\varphi)\in I$, so $\ex^\varphi\comp\ex^{\fd(\varphi)}\flattereq\fn$ for some $\fn\in\fN$.
 Similarly, if~$\psi\ne 0$, then $\ex^\psi\notin\fN$. The desired equivalence now follows from $\fm=\ex^\varphi\cdot\ex^\psi$.
\end{proof}

\begin{cor}
$\mathfrak{N}\cap \mathfrak{R}=\{1\}$ and $G  = \mathfrak{N}\cdot \mathfrak{R}$.
For  $\mathfrak{n}\in \mathfrak{N}$ and  
$\mathfrak{r}\in \mathfrak{R}$ we have: 
$\mathfrak{n}\mathfrak{r}
\succ 1\ \Longleftrightarrow\ \mathfrak{r} \succ 1$, or $\mathfrak{r} =
1$ and $\mathfrak{n} \succ 1$. 
\end{cor}

\noindent
 We think of the monomials in the convex subgroup
$\frak N$ as being {\em flat}. Accordingly
we call~$\mathfrak{R}$ the
{\it steep supplement}\/ of $\mathfrak{N}$.
Then $\T_{\fN}:=\T\cap \R[[\fN]]$ is a strong truncation closed and real closed
subfield of $\T$ and $\T_{\fN}\cap G = \fN$, $\fd(\T_{\fN}^\times)=\fN$. Next we put
$$\fN_{m,n}:=G_{m,n}\cap\fN,\quad \fR_{m,n}:=G_{m,n}\cap\fR. 
$$
so $\fN=\bigcup_{m,n} \fN_{m,n}$ and $\T_{\fN}\  =\  \bigcup_{m,n} \R[[\fN_{m,n}]]$. 

\begin{cor} $G_{m,n} = \mathfrak{N}_{m,n}\cdot \mathfrak{R}_{m,n}$. Moreover, 
 $(\fN_{m,n})_{(m,n)}$ is a healthy directed family of ordered subgroups of $\fN$. 
\end{cor}
\begin{proof} Let $\fm\in G_{m,n}$. Then $\log\fm\in \R\ell_{n+1}+(E_{m-1,n})_\succ$ where~$(E_{-1,n})_\succ:=\{0\}$. Now 
$\fm=\fn\fr$ with $\fn\in \fN, \fr\in \fR$ and
$\supp(\log \fr)\subseteq \supp(\log \fm)\subseteq \{\ell_{n+1}\}\cup G_{m-1,n}^{\succ}$ where the first inclusion holds by the proof of Lemma~\ref{lem:logN}. Hence 
$$\log \fr\in \R\ell_{n+1}+(E_{m-1,n})_\succ,$$
 so $\fr\in G_{m,n}$, and thus $\fn\in G_{m,n}$. 
The rest is obvious.  
\end{proof}

\noindent
If $x\in \fN$, then $\ell_n\in \fN$ for all $n$, so  $\mathfrak L\subseteq \fN$, where  $\mathfrak L:=\bigcup_n \ell_0^\R\cdots\ell_{n}^\R$ is the ordered subgroup
of~$G$ consisting of the {\it logarithmic trans\-monomials}. 

\begin{cor}\label{cor:T_N}
Suppose   $x\in\fN$.
Then  $\T_{\fN}$ is a newtonian $\upo$-free $H$-subfield of~$\T$.
\end{cor}
\begin{proof}
By induction on $m$ we show $\der(\R[[\fN_{m,n}]])\subseteq  \R[[\fN]]$.
We have $\R[[\fN_{0,n}]]=\R[[\ell_n^{\R}]]$
and $\ell_n'=1/(\ell_0\cdots\ell_{n-1})\in\fN$, so  $\der(\R[[\fN_{0,n}]])\subseteq  \R[[\fN]]$.
Let $m\ge 1$; by strong linearity of $\der$ it suffices to show $\der(\fN_{m,n}) \subseteq  \R[[\fN]]$.
Now if $\fm\in \fN_{m,n}$, then~$\log\fm\in (\R\ell_{n+1}+ E_{m-1,n}) \cap\R[[\fN]]\subseteq \R[[\fN_{m-1,n+1}]]$, so
by the inductive hypothesis
$\fm'=\fm\cdot (\log\fm)'\in\R[[\fN]]$ as required.
This proves $\der(\T_{\fN})\subseteq\T_{\fN}$, and the rest now follows from Lemma~\ref{lem:newt} and Corollary~\ref{cor:newt}.
\end{proof}

\noindent
Thus if $x\in\fN$, then $\T_{\fN}$ is linearly surjective by \cite[Corollary~14.2.2]{ADH}, so 
$\T_{\fN}$ is closed under logarithms: $\log\big(\T_{\fN}^>\big)\subseteq\T_\fN$.  Note also that if $x\in \fN$ and $m\ge n$, then~$\der\big(\R[[\fN_{m,n}]]\big)\subseteq \R[[\fN_{m,n}]]$.

\begin{example}
Let   $\frak{g} \in G^{\succ}$ and $x\flatter \fg$. Then  $\mathfrak{N} :=  \{ \mathfrak{n}  : \mathfrak{n}
     \flatter \frak{g} \}$ is a convex subgroup of~$G$ with $x\in \fN$, and for the corresponding $I$ and $F$ we have
$$I\ =\  \{\frak{m} :\  \text{$\fm\succ 1$  and $\exp\frak{m} 
\flatter \frak g$} \}, \qquad F\ =\ \{\fm:\ \fm \succ 1 \text{ and } \fm \succeq 
\log \fg\},$$ and thus 
  \[ \mathfrak{R}\  =\  \bigl\{ \mathfrak{r} \in G :\  \supp \log
     \mathfrak{r} \succcurlyeq \frak{d}(\log \mathfrak{g}) \bigr\}. \]     
  With this $\fN$ we set  $\T_{\flatter \fg}:=\T_{\fN}$.   (In the next subsection we take~${\mathfrak{g}=\ex^x}$.)
\end{example}

\subsection*{The differential field of flat transseries}
We have 
$$\fm\flatter\ex^x \quad\Longleftrightarrow\quad \fm^\dagger=(\log \fm)'\prec 1  \quad\Longleftrightarrow\quad \log \fm \prec x.$$
If one of these equivalent conditions holds, we say that $\fm$ is {\it flat.}\/ 
We call the convex subgroup
$G^\flat:=\{\fm:\fm\flatter \ex^x\}$ of $G$  the
{\em flat part}\/ of $G$. Its steep supplement
$$G^{\sharp}\  :=\  \{ \fm \in G : 
  \supp\log \fm  \succcurlyeq x \}$$
is called the  {\it steep part}\/ of $G$.
We have $x\in G^\flat$, and we put
$$\T^\flat\ :=\  \T_{\flatter \ex^x}\  =\  \T_{G^\flat} \ 
=\  \{f\in \T:\text{all $\fm\in\supp f$ are flat}\}\ =\ \bigcup_{m,n}\R[[G^\flat_{m,n}]].$$
So $\T^\flat$ is a strong truncation closed subfield of $\T$, and $\T^\flat$ is
an $\upo$-free newtonian $H$-subfield of $\T$ with $\fd\big( (\T^\flat)^\times \big)=G^\flat$. Transseries in $\T^\flat$ are called {\it flat.}\/ Note that~$f^\dagger\prec 1$ for nonzero $f\in \T^\flat$. In particular, $\der(\fm)\prec\fm$ for $\fm\in G^\flat$.
 Since
the derivation $\der$ of~$\T$ restricts to a very healthy derivation
on $\T^\flat$ (also denoted by~$\der$), we have for each~${c\in\R}$ 
the strongly linear field automorphism $\exp(c\der)$ of $\T^\flat$, with inverse~$\exp(-c\der)$; see also remarks after Corollary~\ref{cor:strong fixed gp}. 
For~$f\in \T^\flat$ we have
$$\exp(c\der)(f)\ =\ \sum_{n=0}^\infty \frac{f^{(n)}}{n!}c^n,\qquad f\ne 0\ \Rightarrow\  \exp(c\der)(f) 
\sim f. $$
We show that $\exp(c\der)$ commutes with $\exp$ and $\log$:

\begin{lemma}\label{lem:5.11}
Let $f\in\log \big((\T^\flat)^>\big)$, $g\in(\T^\flat)^>$. Then
$$\ex^{\exp(c\der)(f)}\  =\  \exp(c\der)(\ex^f ),\qquad
  \log\!\big(\!\exp(c\der)(g)\big)\ =\ \exp(c\der)(\log g).$$
\end{lemma}
\begin{proof}
First note that $\ex^f\in\T^\flat$ gives $f'=(\ex^f)^\dagger\prec 1$,
so $f^{(n)}\prec 1$ for ${n\ge 1}$,   hence~$\exp(c\der)(f)=f+\varepsilon$
where $\varepsilon\prec 1$. 
Put $A:=\ex^{\exp(c\der)(f)}$, $B:=\exp(c\der)(\ex^f)$.
Then $A\sim \ex^f\sim B$, and  since $A$ and $B$ satisfy the linear differential equation~$y'=\exp(c\der)(f')\cdot y$, we get $A=B$.
This proves the first equality.  
Applying~$\log$ on both sides and taking $f:=\log g$ yields the other one.
\end{proof}

\noindent
As a matter of suggestive notation, we often write~$f(x)$  and $f(g(x))$ instead of
$f$ and $f\circ g$ for $f\in \T$, $g\in\T^{>\R}$. In particular, ${f(x+c)=f\circ (x+c)}$ for $f\in \T$, $c\in \R$. 
Next a special case of the Taylor formula from~\cite[Pro\-po\-si\-tion~5.11(c)]{JvdH}.
(For Taylor expansions in the well-based setting, see also  \cite[Proposition~2.6, Theorem~A.4]{vdH}, with a correction in \cite[Theorem~6.1.4]{Sch}, and for 
generalizations see~\cite{BM}.)

\begin{prop}\label{prop:5.11}
Let $c\in\R$.  Then for $f\in\T^\flat$ we have
\begin{equation}\label{eq:5.11}
f(x+c)\   =\  \exp(c\der)(f)\in\T^\flat.
\end{equation}
\end{prop}

\begin{proof} Clearly, \eqref{eq:5.11} holds for $f\in \R$ and for $f=x$. 
If \eqref{eq:5.11} holds for $g,h\in\T^\flat$ in place of $f$, then it holds for $f=gh$.
If
$(f_\lambda)$ is a summable family in $\T^\flat$ such that  
\eqref{eq:5.11} holds for each $f_\lambda$ in place of~$f$, then
\eqref{eq:5.11} holds for $f=\sum_\lambda f_\lambda$.
Let now~$f\in \log\big((\T^\flat)^>\big)$ and suppose \eqref{eq:5.11} holds.
Then
$$(\ex^f)  (x+c)\  =\  \ex^{f  (x+c)}\  =\ 
\ex^{\exp(c\der)(f)}\  =\  \exp(c\der)(\ex^f)$$
where we used Lemma~\ref{lem:5.11} for the last equality. 
Likewise one shows that if  \eqref{eq:5.11} holds for a given $f\in (\T^\flat)^>$,
then  \eqref{eq:5.11} holds with $\log f$ in place of $f$.  Hence for all~$n$ and all $r\in \R$, \eqref{eq:5.11} holds for $f=\ell_n^r=\ex^{r\log \ell_n}$. Now $G^\flat_{0,n}=\ell_n^\R$, so \eqref{eq:5.11} holds for all $f\in \R[[G^\flat_{0,n}]]$. Suppose for a certain~$m$ that \eqref{eq:5.11} holds for all $n$ and all $f\in \R[[G^\flat_{m,n}]]$, and let
$\fm\in G^\flat_{m+1,n}$. Then $\log \fm=r\ell_{n+1}+ g$ with $r\in \R$ and~$g\in \big(E_{m,n})_{\succ}\cap \T^\flat\subseteq \R[[G^\flat_{m,n}]]$.  Hence \eqref{eq:5.11} holds for $f=\fm=\ell_n^r \ex^g$, and thus for all $f\in \R[[G^\flat_{m+1,n}]]$. 
This inductive argument shows that \eqref{eq:5.11} holds for all~$f\in \T^\flat$.
\end{proof}

\begin{cor}
Let $c\in\R$ and $f\in\T^\flat\setminus (\R+\R x+\cdots+\R x^{m-1})$. Then 
$$f  (x+c) - \displaystyle\sum_{n=0}^m \frac{f^{(n)}}{n!}c^n\  \prec\  f^{(m)}.$$
\end{cor}
\begin{proof}
For   $n\le m$ we have $f^{(n)}\ne 0$, and
since $f^{(n)}\in\T^\flat$ we also have $f^{(n+1)}\prec f^{(n)}$.
Now use Proposition~\ref{prop:5.11}.
\end{proof}

\noindent
Proposition~\ref{prop:5.11} gives for $c\in \R$ the strong automorphism $f(x)\mapsto f(x+c)$ of the field~$\T^\flat$. 
Consider the  strongly linear map 
$$\Delta^\flat\ \colon\ \T^\flat\to\T^\flat,\qquad \Delta^\flat(f)\ :=\ f(x+1)-f(x).$$
We show that $\Delta^\flat$ is surjective; in fact:

\begin{prop}\label{prop:Delta flat}
There is a  strongly linear map $\Sigma^\flat\colon\T^\flat\to\T^\flat$  with $\Delta^\flat\circ\Sigma^\flat=\id_{\T^\flat}$.
\end{prop}
\begin{proof} As before, we denote the restriction of the derivation $\der$ on $\T$ to a (very healthy) derivation on $\T^\flat$ also by
$\der$. Then $\Delta^\flat=\big(\exp(t)-1\big)(\der)$ by Proposition~\ref{prop:5.11}.  
The power series $S\in\Q[[t]]$  from Remark~\ref{rem:S} satisfies $\big(\exp(t)-1\big)\cdot S=t$, so~${\Delta^\flat\circ S(\der)}=\der$ on
$\T^\flat$.  For $m\ge n\ge 1$, the integration map $\int\colon \T\to \T$ of Lem\-ma~\ref{lemint} maps  $\R[[G^\flat_{m-1,n-1}]]=\T^\flat\cap \R[[G_{m-1,n-1}]]$ into $\R[[G^\flat_{m,n}]]=\T^\flat\cap \R[[G_{m,n}]]$, and so restricts to a strongly linear map $\int^\flat\colon \T^\flat\to \T^\flat$. Thus $\Sigma^\flat:=S(\der)\circ \textstyle\int^\flat$ has the required properties. 
\end{proof}

\begin{remark} 
One can show more generally that every linear difference operator \marginpar{I (Lou) did not check this remark} 
$$f(x)\mapsto f(x+n) + a_{1} f(x+n-1) + \cdots + a_n f(x)\colon\T^\flat\to\T^\flat
\quad (a_1,\dots,a_{n}\in\T^\flat)$$
is surjective. See \cite[Section~8]{DO}, which studies the model theoretic properties of the field $\T$ equipped with
its automorphism~$f(x)\mapsto f(x+1)$ and
the flattened  valuation~$\T^\times\xrightarrow{\fd} G^{\operatorname{LE}}=G^\flat \cdot G^\sharp\to  G^\sharp$ 
(the second arrow being the natural projection). 
\end{remark}

\section{Logarithmic Transseries}\label{sec:logtrans} 

\noindent
The elements of the  subfield $\T_{\log}=\bigcup_n \R[[\ell_0^\R\cdots\ell_{n}^\R]]$ of~$\T$ are called
{\it logarithmic transseries.}\/
More generally, let $\k$ be a subring of $\R$ (so $\Z\subseteq \k$). 
Then~$\mathfrak L^{\k}:=\bigcup_n \mathfrak L_n^\k$ is the increasing union of the ordered subgroups 
$\mathfrak L^\k_n:=\ell_0^\k\cdots\ell_{n}^\k$ of $\mathfrak L$,
and
 $\T_{\log}^{\k}:=\bigcup_n\R[[\mathfrak L_n^\k]]$ is a subfield of $\T_{\log}$, truncation closed in $\R[[G^{\operatorname{LE}}]]$. 
For~$f\in \R[[\mathfrak L_n^\k]]^>$, $r\in\k$  we have $f^r\in \R[[\mathfrak L_n^\k]]$.
Also $\log\big((\T_{\log}^{\k})^>\big) \subseteq \T_{\log}^{\k}$; indeed,
$$\log \big(\R[[\mathfrak L_n^\k]] ^>\big)  \ =\  \left(\bigoplus_{i=1}^{n+1} \k\ell_i\right) \oplus \R[[\mathfrak L_n^\k]] _{\preceq} \qquad\text{(internal direct sum)}$$
and thus
 $$\log\big((\T_{\log}^\k)^>\big)\  =\  \left(\bigoplus_{i=1}^\infty \k\ell_i\right)  \oplus (\T_{\log}^\k)_{\preceq}\qquad\text{(internal direct sum).}$$
 We now collect further facts about~$\T_{\log}$ and its subfields~$\T_{\log}^{\k}$ as needed later.
 
\begin{lemma}\label{lem:intersect with Tlog}
If $m\leq n$, then
\begin{equation}\label{eq:intersect with Tlog, 1}
G_{m,n}\cap\mathfrak L^{\k}\  =\  \ell_{n-m}^\k\cdots\ell_{n}^\k\quad\text{ and }\quad
E_{m,n}\cap\T_{\log}^\k\ =\  \R[[\ell_{n-m}^\k\cdots\ell_{n}^\k]]
\end{equation}
and if $m\geq n$, then
\begin{equation}\label{eq:intersect with Tlog, 2}
G_{m,n}\cap\mathfrak L^\k\ =\ \ell_{0}^\k\cdots\ell_{n}^\k\quad\text{ and }\quad
E_{m,n}\cap\T_{\log}^\k\ =\ \R[[\ell_0^\k\cdots\ell_{n}^\k]].
\end{equation}
In particular,  $E{\downarrow_n} \cap \T_{\log}^{\k}\  =\  E_{n,n}\cap\T_{\log}^{\k}\ =\ \R[[\ell_0^\k\cdots\ell^\k_n]]$.
\end{lemma}
\begin{proof}
We show \eqref{eq:intersect with Tlog, 1} by induction on $m$. The case $m=0$ holds, since
$G_{0,n}\cap\mathfrak L^\k=\ell_n^\R\cap\mathfrak L^\k=\ell_n^\k$ and
thus $E_{0,n}\cap \T_{\log}^\k =\R[[G_{0,n}]]\cap\T_{\log}^{\k}=\R[[\ell_n^\k]]$.
Suppose  \eqref{eq:intersect with Tlog, 1} holds for a certain  $m\in\{0,\dots,n-1\}$. Then
by Corollary~\ref{cor:distributive} and inductive hypothesis,
$$(E_{m,n})_\succ \cap \bigcup_{i\in\N} \R\ell_1+\cdots+\R\ell_{i+1}\  =\  
\R\ell_{n-m}+\cdots+\R\ell_n.$$
It follows that this remains true with $\k$ instead of $\R$. Thus
\begin{align*}
G_{m+1,n}\cap \mathfrak L^\k \	&=\  \ell_n^\R \exp\big((E_{m,n})_\succ\big) \cap\mathfrak L^{\k} \\
				&=\ \ell_n^\k\left(\exp\big((E_{m,n})_\succ\big)\cap \bigcup_{i\in\N} \exp(\k\ell_1+\cdots+\k\ell_{i+1})\right)  \\
							&=\  \ell_n^{\k} \exp(\k\ell_{n-m}+\cdots+\k\ell_n)\ 
							 =\ \ell_{n-(m+1)}^\k\cdots\ell_{n-1}^\k\ell_n^\k,
\end{align*}
hence $E_{m+1,n}\cap\T^\k_{\log} = \R[[G_{m+1,n}]]\cap\T^\k_{\log}=\R[[\ell_{n-(m+1)}^\k\cdots\ell_{n}^\k]]$, as promised.

Next we prove \eqref{eq:intersect with Tlog, 2} by induction on $m\geq n$.
We established the case~$m=n$, so assume
\eqref{eq:intersect with Tlog, 2} for a certain $m\geq n$. Then
by Corollary~\ref{cor:distributive} and inductive hypothesis,
$$(E_{m,n})_\succ \cap \bigcup_{i\in\N} \R\ell_1+\cdots+\R\ell_{i+1}\  =\ 
\R\ell_{1}+\cdots+\R\ell_n.$$
As before we obtain
$G_{m+1,n}\cap\mathfrak L^{\k} 	= \ell_{0}^\k\cdots\ell_n^\k$
and so  
\[\qquad E_{m+1,n}\cap\T_{\log}^\k\  =\  \R[[G_{m+1,n}]]\cap\T_{\log}^{\k}\ =\ \R[[\ell_0^\k\cdots\ell_{n}^\k]].\qquad \qquad \qedhere \]
\end{proof}

\begin{cor} 
The additive subgroup $\T_{\log}^\k$ of $\T$ is strong.
\end{cor}
\begin{proof}
Let $(f_\lambda)$ be a family of elements of $\T_{\log}^\k$ which is summable in $\T$.
Take~$m$,~$n$  such that $f_\lambda\in\R[[G_{m,n}]]$ for all $\lambda$
and $(f_\lambda)$ is summable in $\R[[G_{m,n}]]$. Then $f_\lambda\in \R[[\ell_0^\k\cdots\ell^\k_n]]$ for all $\lambda$, by
Lemma~\ref{lem:intersect with Tlog}, so $\sum_\lambda f_\lambda\in\T_{\log}^\k$.
\end{proof}

\noindent
Thus $\T_{\log}^\k=\T\cap \R[[\mathfrak{L}^{\k}]]$ by Lemma~\ref{lem:fixsigma, 1}. Recall that $\ex_0=x$, $\ex_1=\ex^x$, $\ex_2=\ex^{\ex^x}$, etc. 

\begin{cor}\label{cor:Em it exp} Let $m\le n$. Then
$E_m\cap\R[[\ex_n^\k\cdots\ex_0^\k]]=\R[[\ex_m^\k\cdots\ex_0^\k]]$ if $m\leq n$.
\end{cor}
\begin{proof}
Applying ${\downarrow_n}$ on both sides, this is equivalent to
$E_{m,n}\cap \R[[\ell_0^\k\cdots\ell_n^\k]]=\R[[\ell_{n-m}^\k\cdots\ell_n^\k]]$, which 
 is part of \eqref{eq:intersect with Tlog, 1}.
\end{proof}

\subsection*{Substitution into logarithmic transseries}
This subsection will be applied in the next section.  Let $\k$ be a subring of $\R$ and $K:=\R[[\ell_0^\k\cdots\ell_n^\k]]$.
Fix $g=(g_0,\dots,g_n)\in K^{1+n}$ with $g_i \prec \ell_i$ for $i=0,\dots,n$,
and let $r=(r_0,\dots,r_n)$ range over~$\k^{1+n}$. Given $f=(f_0,\dots,f_n)\in (K^>)^{1+n}$,    set $f^r:=f_0^{r_0}\cdots f_n^{r_n}\in K$. We also put $\ell:=(\ell_0,\dots,\ell_n)$.
 
\begin{lemma}\label{lem:subst into ell}
The map $\ell^r\mapsto (\ell+g)^r\colon \ell_0^\k\cdots\ell_n^\k \to K$ is well-based.
\end{lemma}

\begin{proof}
Setting $\epsilon_i:=g_i/\ell_i\prec 1$ for $i=0,\dots,n$ we have
\begin{align*} (\ell+g)^r\ &=\ \ell^r \cdot (1+\epsilon_0)^{r_0} \cdots (1+\epsilon_n)^{r_n}, \text{ with}\\
(1+\epsilon_i)^{r_i}\ &=\ \sum_{j\in\N} {r_i \choose j} \epsilon_i^{j}.
\end{align*}
So $\supp (\ell+g)^r \subseteq \ell^r S_0^{*}\cdots S_n^{*}$ where $S_i:=\supp\epsilon_i\prec 1$. 
Thus if $\mathfrak S$ is a well-based subset of $\ell_0^\k\cdots\ell_n^\k$, 
then $\big((\ell+g)^r\big){}_{\ell^r\in\mathfrak S}$ is summable by Neumann's Lemma~\ref{lem:Neumann}.
\end{proof}
 
\noindent
Proposition~\ref{prop:extend to strongly linear} and Lemma~\ref{lem:subst into ell} yield a 
strongly linear map~$K\to K$ which sends~$f=\sum_r f_r\ell^r\in K$ ($f_r\in\R$)
to $f_{+g}:=\sum_r f_r (\ell+g)^r$. By Corollary~\ref{cor:strongly linear field embedd}, this map is a field embedding. Moreover,
$\big(f_{+g})^s=(f^s)_{+g}$ for $f\in K^{>}$, $s\in \k$. 
The proof above also gives $(\ell+g)^r\sim \ell^r$. Thus by Corollary~\ref{prop:inverse strongly linear}:


\begin{cor}\label{cor:subst into ell}
The strongly linear field embedding $f\mapsto f_{+g}\colon  K\to K$
is bijective with strongly linear inverse. For $g\in \R^{1+n}$ this inverse is $f\mapsto f_{-g}$.  
\end{cor}

\subsection*{Modifying the logarithmic sequence}
Let now a sequence $\alpha=(\alpha_n)$ of real numbers  be given. We
recursively define the elements $\ell_n^*$   of~$\T_{\log}^{>\R}$  by
$$\ell_0^*\ :=\ \ell_0+\alpha_0,\qquad \ell_{n+1}^*\ :=\ \log\ell_n^* +  \alpha_{n+1}.$$
Then $(\ell_{n+1}^*)'=(\ell_n^*)^\dagger$ and $\ell_n^*\le 2\ell_n$ for all $n$, so $(\ell_n^*)$ is a logarithmic sequence for~$\T$ in the sense of \cite[Section~11.5]{ADH}. 
An easy induction on $n$ yields $(\ell_{n}^*)'=\frac{1}{\ell_0^*\cdots\ell_{n-1}^*}$. 
In the next lemma we let~$t_0,t_1,\dots$ be distinct indeterminates over $\R$.

\begin{lemma}\label{lem:ellnstar}
Recursion on $n$ yields power series $G_n=G_n^\alpha\in\R[[t_0,\dots,t_{n-1}]]$ by
$$G_0\ :=\ \alpha_0,\qquad G_{n+1}\ :=\ \alpha_{n+1}+\log(1+G_n t_n),$$
and then
$$\ell_n^*\  =\  \ell_n + G_n(1/\ell_0,\dots,1/\ell_{n-1})\quad\text{and}\quad G_n(0)\ =\ \alpha_n.$$
\end{lemma}
\begin{proof}
This is clear for $n=0$.
Suppose the claim holds for a certain $n$. Then
$$g:=\ell_n^*-\ell_n=G_n(1/\ell_0,\dots,1/\ell_{n-1})\preceq 1\prec\ell_n,$$ so
$\ell_{n+1}^*= \log(\ell_n+g)+\alpha_{n+1}$, and thus
$$\ell_{n+1}^*\ =\ \log\left(\ell_n\left(1+\frac{g}{\ell_n}\right)\right)+\alpha_{n+1}\ =\ \ell_{n+1} + \alpha_{n+1} + \log\left(1+\frac{g}{\ell_n}\right)$$
where $\log\big(1+\frac{g}{\ell_n}\big)=\sum_{m=1}^\infty \frac{(-1)^{m+1}}{m} (g/\ell_n)^m$.
Hence  $G_{n+1}\in \R[[t_0,\dots,t_n]]$ has the required properties.
\end{proof}

\begin{remark}
If $\alpha_1=\cdots=\alpha_{n}=0$, then $G_{n+1} \in \alpha_{n+1}+  t_0 \R[[t_0,\dots,t_{n}]]$.
\end{remark}

\begin{lemma}\label{lem:ellnstarinverse}
$\ell_n\circ (\ell_n^*)^{\inv}\in \T_{\exp}$.
\end{lemma}
\begin{proof}
Recursion on $n$ defines elements $\ex_n^*$ of $\T_{\exp}^{>\R}$ by 
$\ex_0^* := x-\alpha_0$ and $\ex_{n+1}^*  := \ex_n^* \circ \ex^{x-\alpha_{n+1}}$.
An easy induction on $n$ shows that $\ex_n^*=(\ell_n^*)^{\inv}$.
We now prove by induction on~$n$ that $\ell_n\circ\ex_n^* \in \T_{\exp}^{>\R}$. This is clear for $n=0$, so suppose it holds for a certain $n$. Then
$$\ell_{n+1}\circ \ex_{n+1}^*\  =\   \ell_{n+1}\circ (\ex_n^* \circ \ex^{x-\alpha_{n+1}})\ =\ 
\ell_1\circ (\ell_n\circ\ex_n^*) \circ \ex^{x-\alpha_{n+1}}$$
and $f:= \ell_n\circ\ex_n^* \in \T_{\exp}^{>\R}$ by inductive hypothesis.
Take $c\in\R^>$, $r\in \R$, $g\in (\T_{\exp})_\succ$, and $\varepsilon\in (\T_{\exp})_\prec$  such that
$f=cx^r\ex^g(1+\varepsilon)$. Then 
$$\ell_1\circ f\  =\  \log f = \log c + r\log x + g + \log(1+\varepsilon)\qquad\text{where $\log(1+\varepsilon)\in \T_{\exp}$}$$
and hence $\ell_{n+1}\circ \ex_{n+1}^*= \ell_1\circ f\circ  \ex^{x-\alpha_{n+1}}$ yields
$$\ell_{n+1}\circ \ex_{n+1}^*\  =\  \log c + r\cdot(x-\alpha_{n+1}) + (g\circ \ex^{x-\alpha_{n+1}}) +
\log(1+\varepsilon\circ \ex^{x-\alpha_{n+1}})$$
with $g\circ \ex^{x-\alpha_{n+1}}, \log(1+\varepsilon\circ \ex^{x-\alpha_{n+1}})\in \T_{\exp}$, so $\ell_{n+1}\circ \ex_{n+1}^*\in \T_{\exp}$.
\end{proof}

\begin{lemma}\label{lem:id on elln}
The following are equivalent:
\begin{enumerate}
\item[\rm(i)] $\ell_n^*=\ell_n$;
\item[\rm(ii)] $\ell_i^*=\ell_i$ for $i=0,\dots,n$;
\item[\rm(iii)] $\alpha_0=\cdots=\alpha_n=0$.
\end{enumerate}
\end{lemma}
\begin{proof}  By induction on $n$. The case $n=0$ is trivial. 
Suppose the equivalences hold for a certain $n$, and assume $\ell_{n+1}^*=\ell_{n+1}$.
Then
$$\ell_n\ =\ \ex^{\ell_{n+1}}\ =\ \ex^{\ell_{n+1}^*}\ =\ \ex^{\alpha_{n+1}+\log\ell_n^*}\in \R^> \cdot \ell_n^*.$$
Also $\ell_n^*-\ell_n \preceq 1$ by Lemma~\ref{lem:ellnstar}, so  $\ell_n^*=\ell_n$,
and thus $\ell_i^*=\ell_i$ for $i=0,\dots,n$ by inductive hypothesis.
This proves (i)~$\Rightarrow$~(ii) for $n+1$ in place of $n$. The proofs of~(ii)~$\Rightarrow$~(iii)~$\Rightarrow$~(i) with $n$ replaced by $n+1$ are obvious. 
\end{proof}

\begin{remark}
Suppose $m\geq 1$ and $\ell_{m-1}^*=\ell_{m-1}$, and let $G_n=G_n^\alpha$ be the power series from Lemma~\ref{lem:ellnstar}. Then $G_0=\cdots=G_{m-1}=0$, and for $n\geq m$ we have~$G_{n}\in \R[[t_{m},\dots,t_{n-1}]]$;
indeed, $G_n=G_{n-m}^\beta(t_m,\dots,t_{n-1})$ where $\beta_n=\alpha_{m+n}$.
\end{remark}

\section{Strong Automorphisms of $\T_{\log}$}\label{sec:Tlog}

\noindent
We have $\der(\T_{\log})\subseteq \T_{\log}$, and below we view $\T_{\log}$ as a differential subfield of $\T$.  More generally, fix a subring $\k\supseteq \Q$ of $\R$. Then $\der\big(\R[[\mathfrak{L}_n^\k]]\big)\subseteq  \R[[\mathfrak{L}_n^\k]]$, so 
$\T_{\log}^\k=\bigcup_n  \R[[\mathfrak{L}_n^\k]]$ is a differential subfield of $\T_{\log}$ (with constant field $\R$).
In this section we clarify the structure of the group~$\SigmaAut_\der(\T_{\log}^\k)$ of strong automorphisms of this differential field. First we
analyze  strongly additive differential field embeddings~$\T_{\log}^{\k}\to\T$; this will
also be useful for the proof of Theorem A in the next section.

\subsection*{Properties of strongly additive embeddings }
Let
$K:=\R[[\ell_0^{\k}\cdots\ell_n^{\k}]]$, and let
$\sigma\colon K\to\T$ be a strongly additive differential field embedding. Then $\sigma$ is the identity on $\R$, so  
 $\sigma$ is $\R$-linear. Since $\k$ is divisible, $K$ is real closed, so $\sigma$ is order-preserving, and thus a valued field embedding (between valued subfields of the Hahn field 
$\R[[G^{\operatorname{LE}}]]$). From the beginning of the previous section we obtain
$$\log(K^>)\cap K\  =\  \left(\bigoplus_{i=1}^n \k\ell_i\right) \oplus K_{\preceq}.$$
Define $\alpha\colon \log(K^>)\cap K\to\T$ by
$$ \alpha(f)\ :=\  \log\big(\sigma(\exp f)/\exp \sigma(f)\big)\  =\ \log\sigma(\exp f) - \sigma(f).$$
As in the proof of Lemma~\ref{lem:sigma alpha} we see that $\alpha$ is additive and takes values in $\R$.

\begin{lemma}\label{lem:alpha}
We have   $K_\preceq\subseteq\ker\alpha$ and so $\alpha(f)=\alpha(f_\succ)$ for all $f\in \log(K^{>})\cap K$.
\end{lemma}
\begin{proof}
For $c\in\R$ we have 
$\alpha(c)=\log\sigma(\ex^c)-\sigma(c)=\log \ex^c - c = 0$. 
For 
$\varepsilon\in K_\prec$ we have $\exp \varepsilon=\sum_n \frac{1}{n!}\varepsilon^n$ and
so
$\sigma(\exp \varepsilon)=\sum_n \frac{1}{n!}\sigma(\varepsilon)^n=\exp{\sigma(\varepsilon)}$ by strong additivity of $\sigma$,
thus $\alpha(\varepsilon)=0$.   
\end{proof}

\noindent
For $f\in K^>$  with $\log f\in K$ we have
$\sigma(\log f)=\log\sigma(f)-\alpha(\log f)$.  
In particular:
\begin{equation}\label{eq:sigma ell_i}
\sigma(\ell_{i+1})\ =\  \log \sigma(\ell_i)-\alpha(\ell_{i+1}) \qquad\text{for $i=0,\dots,n-1$.}
\end{equation}
From $\sigma(x)'=\sigma(x')=\sigma(1)=1$ we get $\sigma(x)\in x+\R$.
So by Lemma~\ref{lem:ellnstar} and~\eqref{eq:sigma ell_i}:

\begin{lemma} \label{lem:sigma elli}
$\sigma(\ell_i)\in\ell_i-\alpha(\ell_i)+\R[[\ell_0^\Z\cdots\ell_{i-1}^\Z]]_{\prec}$ for $i=1,\dots,n$.
\end{lemma}
 

\noindent
Let $f\in K^{>}$, $r\in \k$. Then $f^r\in K^{>}$ and $(f^r)^\dagger=rf^\dagger$ with $f^\dagger\in K$, and likewise~$(\sigma(f)^r)^\dagger=r\sigma(f)^\dagger$. It follows that 
$\big(\sigma(f^r)\big){}^\dagger = \big(\sigma(f)^r\big){}^\dagger$ and thus
$$\beta(f,r):=\sigma(f^r)/\sigma(f)^r\in\R^>.$$  Given any $r\in \k$ the map 
$f\mapsto \beta(f,r)\colon K^{>}\to \R^{>}$ is a morphism of multiplicative groups, and given any $f\in K^{>}$ the map 
$r\mapsto \beta(f,r) \colon \k \to \R^{>}$ is a morphism of the additive group $\k$ to the multiplicative group $\R^>$. 
If $f\in K^>$, $\log f\in K$, and $r\in \k$, then
$$\alpha(\log f^r)-r\alpha(\log f)=\log \sigma(f^r)-r\log \sigma(f),$$ hence
$$\beta(f,r)\  =\  \ex^{\alpha(r\log f)-r\alpha(\log f)}.$$
In particular 
\begin{equation}\label{eq:sigmaellir}
\sigma(\ell_i^r)\ =\ \ex^{\alpha(r\ell_{i+1})-r\alpha(\ell_{i+1})} \sigma(\ell_i)^r \qquad\text{for $i=0,\dots,n-1$ and $r\in\k$.}
\end{equation}
For $r=(r_0,\dots,r_{n})\in\k^{1+n}$ put $\chi(\ell^r):=\prod_{i=0}^n \beta(\ell_i,r_i)$. Then 
$$\chi\ \colon\ \mathfrak L_n^\k=\ell_0^{\k}\cdots\ell_n^{\k}\to\R^>$$ is a group morphism, and
$\sigma(\ell^r) = \chi(\ell^r)\sigma(\ell)^r$
for all $r\in\k^{1+n}$, where $\sigma(\ell)=\big(\sigma(\ell_0),\dots,\sigma(\ell_n)\big)$.

\begin{lemma}\label{lem:sigma is strong}
We have
$\sigma(K) =  K$, and $\sigma\colon K\to K$ is strong.
\end{lemma}
\begin{proof}
By Lemma~\ref{lem:sigma elli} we have
$\sigma(\ell_i)\in \R[[\ell_0^\Z\cdots\ell_i^\Z]]\subseteq K$ for $i=0,\dots,n$; 
therefore~$\sigma(K)\subseteq K$, 
by strong linearity of $\sigma$ and the remarks preceding the lemma.
Indeed, by Corollary~\ref{cor:subst into ell}, we have 
the strong field automorphism $f\mapsto f_{+g}$ of $K$,
where $g=(g_0,\dots,g_n)$ with $g_i:=\sigma(\ell_i)-\ell_i \preceq 1\prec \ell_i$ for $i=0,\dots,n$, and $\sigma$ is the composition of this automorphism with the strong automorphism $f\mapsto f_{\times\chi}$ of $K$ defined after
Corollary~\ref{cor:strongly linear inverse}: $\sigma(f)=f_{\times\chi,+g}$.
Hence $\sigma(K)=K$, and $\sigma$ is strong.
\end{proof}

\begin{cor}\label{corsad} 
Each strongly additive differential field embedding $\T_{\log}^\k\to\T$ has image $\T_{\log}^\k$,
and is a strong healthy automorphism of the differential field $\T_{\log}^{\k}$.
\end{cor}

\begin{cor}\label{cor:small pert of id}
The following are equivalent:
\begin{enumerate}
\item[\textup{(i)}] $\sigma(f)\sim f$ for all $f\in K^\times$;
\item[\textup{(ii)}] $\sigma(\ell_i^r)=\sigma(\ell_i)^r$ for $i=0,\dots,n$ and all $r\in\k$;
\item[\textup{(iii)}] $\sigma(f^r)=\sigma(f)^r$ for all $f\in K^>$ and $r\in\k$.
\end{enumerate}
\end{cor}
\begin{proof}
Assume (i). Let $i\in\{0,\dots,n\}$ and $r\in\k$; then $\sigma(\ell_i)\sim\ell_i$
and~$\sigma(\ell_i^r)\sim\ell_i^r$, so~$\sigma(\ell_i)^r\sim\ell_i^r\sim\sigma(\ell_i^r)$, hence
$1\sim\sigma(\ell_i^r)/\sigma(\ell_i)^r=\beta(\ell_i,r)\in\R^>$ and thus~$\sigma(\ell_i^r)=\sigma(\ell_i)^r$.
This proves (i)~$\Rightarrow$~(ii), and  (ii)~$\Rightarrow$~(iii) follows from the definitions
and strong linearity of $\sigma$. Now assume (iii).  Then $\beta$ and $\chi$ take the constant value~$1$. Hence the proof of Lemma~\ref{lem:sigma is strong} gives $\sigma(f)=f_{+g}$ where~$g=(g_0,\dots,g_n)$, $g_i=\sigma(\ell_i)-\ell_i$ for $i=0,\dots,n$, and $f_{+g}\sim f$ by the remark before Corollary~\ref{cor:subst into ell}.
This proves~(iii)~$\Rightarrow$~(i).
\end{proof}

\subsection*{Structure of $\SigmaAut_\der(\T_{\log}^\k)$}
Let $\sigma$ be a strong automorphism of the differential field $\T_{\log}^\k$.
We associate to $\sigma$ the additive map $\alpha\colon\log\big( (\T_{\log}^\k)^>\big) \to \R$
given by 
$$\alpha(f)\ :=\  \log\sigma(\exp f) - \sigma(f).$$
Its kernel contains $(\T_{\log}^\k)_\preceq$ by Lemma~\ref{lem:alpha}. By 
\eqref{eq:sigma ell_i} and \eqref{eq:sigmaellir}, for all $n$ and $r\in\k$:
$$\sigma(\ell_{n+1})\ =\  \log \sigma(\ell_n)-\alpha(\ell_{n+1}),\qquad
\sigma(\ell_n^r)\ =\ \ex^{\alpha(r\ell_{n+1})-r\alpha(\ell_{n+1})} \sigma(\ell_n)^r.
$$
Hence  
from $c_\sigma:=\sigma(x)-x\in\R$ and the restriction $\gamma_\sigma$ of $\alpha$ to the subgroup~$L_\k:=\bigoplus_{n\geq 1} \k\ell_n$ of $\log\big( (\T_{\log}^\k)^>\big)$ we recover $\sigma$. 
We thus have an injective map
\begin{equation}\label{eq:csigma,alphasigma}
\sigma\mapsto (c_\sigma, \gamma_{\sigma})\ \colon\  \SigmaAut_\der(\T_{\log}^\k) \to
\R\times \Hom(L_\k,\R),
\end{equation}
where $\Hom(L_\k,\R)$ denotes the group of additive maps $L_\k\to\R$.

\begin{prop}\label{prop:SigmaAut of Tlog}
The map \eqref{eq:csigma,alphasigma} is an isomorphism of groups.
\textup{(}In particular, the group $\SigmaAut_\der(\T_{\log}^\k)$ is abelian, torsion-free, and divisible.\textup{)}
\end{prop}

\noindent
Before we give the proof of this proposition we observe:

\medskip\noindent

\begin{lemma}~\label{lemobs}
Let $\tau\in\SigmaAut_\der(\T_{\log}^\k)$. Then 
$$\tau\left(\log\big( (\T_{\log}^\k)^>\big)\right)\ =\ \log\big( (\T_{\log}^\k)^>\big) \text{ and }\alpha\big( \tau(f)\big)\ =\ \alpha(f) \text{ for }
f\in \log\big((T_{\log}^\k)^>\big).$$
\end{lemma}
\begin{proof}
This follows from $\tau\big( (\T_{\log}^\k)_\preceq \big) = (\T_{\log}^\k)_\preceq$, the fact that by Lemma~\ref{lem:sigma elli} we have~$\tau(\ell_n)-\ell_n\preceq 1$ for all $n$, and
$\ker\alpha\supseteq (\T_{\log}^\k)_\preceq$.
\end{proof}

\noindent
For the proof of Proposition~\ref{prop:SigmaAut of Tlog},
we first deal with  surjectivity. Thus let $c\in\R$ and~$\gamma\in \Hom(L_\k,\R)$. Inductively define the elements $\ell_n^*$   of~$(\T_{\log}^\k)^{>\R}$  by
$$\ell_0^*\ :=\ \ell_0+c,\qquad \ell_{n+1}^*\ :=\ \log\ell_n^* -  \gamma(\ell_{n+1}).$$ 
Put $g_n:=\ell_n^*-\ell_n$, so $g_n\preceq 1 \prec \ell_n$ by   Lemma~\ref{lem:ellnstar}.
Fix $n$ and put $K_n:=\R[[\mathfrak L_n^k]]$ and~$g:=(g_0,\dots,g_n)$. 
For $r=(r_0,\dots,r_n)\in\k^{1+n}$, set 
$$\chi(\ell^r)\ :=\ \prod_{i=0}^n \ex^{\gamma(r_{i}\ell_{i+1})-r_{i}\gamma(\ell_{i+1})}.$$ 
The map $\chi\colon \mathfrak L_n^k\to\R^>$ is a group morphism.
Let $\sigma_n$ be the strong automorphism of the field~$K_n$ 
given by $\sigma_n(f)=f_{\times\chi,+g}$ for $f\in K_n$.
Let $i\in\{0,\dots,n\}$ and $r\in\k$. Then $\sigma_n(\ell_i^r)=\ex^{\gamma(r\ell_{i+1})-r\gamma(\ell_{i+1})}(\ell_i^*)^r$, in particular,
$\sigma_n(\ell_i)=\ell_i^*$. 
 We have
$$(\ell_i^r)'\  =\  r \frac{\ell_i^r}{\ell_0\cdots\ell_i}$$
and 
$$\big(\sigma_n(\ell_i^r)\big)^\dagger\ =\  \left( \ex^{\gamma(r\ell_{i+1})-r\gamma(\ell_{i+1})} (\ell_i^*)^r\right)^\dagger \
							=\  \big(  (\ell_i^*)^r \big)^\dagger\ 
							=\  r (\ell_i^*)^\dagger\ =\  r\frac{1}{\ell_0^*\cdots\ell_i^*}.
$$
Using first the identity for $(\ell_i^r)'$ and then the one for $\big(\sigma_n(\ell_i^r)\big)^\dagger$, we obtain
$$\sigma_n\big( (\ell_i^r)' \big)\ =\  r\frac{\sigma_n(\ell_i^r)}{\ell_0^*\cdots\ell_i^*}\  =\  \big(\sigma_n(\ell_i^r)\big)^\dagger\cdot\sigma_n(\ell_i^r)\  =\  
\big(\sigma_n(\ell_i^r)\big)'.$$
Thus $\sigma_n(\fm')=\big(\sigma_n(\fm)\big){}'$ for all $\fm\in\mathfrak L_n^k$,
hence by strong linearity of $\sigma_n$ and of the derivation of $K_n$, it follows that $\sigma_n(f')=\big(\sigma_n(f)\big){}'$ for all $f\in K_n$.
Clearly $\sigma_{n+1}$ extends $\sigma_n$. The common extension $\sigma$ of all~$\sigma_n$ to an operator on 
$\T_{\log}^\k$ is a strong differential field automorphism of $\T_{\log}^\k$ with $c_{\sigma}=c$.  For $n\ge 1$ and $r\in \k$ we have
\begin{align*}  \log \sigma\big(\!\exp(r\ell_n)\big)\ 	&=\ 	\log \sigma\big(\!\exp(r\log\ell_{n-1})\big)\
						=\	\log \sigma\big(\!\exp(\log \ell^r_{n-1})\big) \\
						&=\	\log \sigma(\ell^r_{n-1}) \ 
						=\	\gamma(r\ell_n)-r\gamma(\ell_n) + \log (\ell_{n-1}^*)^r \\
						&=\	\gamma(r\ell_n)-r\gamma(\ell_n) + r \log (\ell_{n-1}^*)\\  
						&=\	\gamma(r\ell_n)-r\gamma(\ell_n) + r \ell_n^* + r \gamma(\ell_n)\\
						&=\	\gamma(r\ell_n) + r\ell_n^*\ =\ \gamma(r\ell_n)+\sigma(r\ell_n),
\end{align*} 
so $\gamma_\sigma(r \ell_n) =  \log \sigma\big(\!\exp(r\ell_n)\big) - \sigma(r\ell_n) =  \gamma(r \ell_n)$. Thus $ \gamma_\sigma=\gamma$. 

\medskip
\noindent
Now let $\sigma,\tau\in \SigmaAut_\der(\T_{\log}^\k)$ have respective associated additive maps 
$$\alpha,\beta\ \colon\ \log\big( (\T_{\log}^\k)^>\big) \to \R.$$  Then for $f\in \log\big( (\T_{\log}^\k)^>\big)$ we have
$\tau(\exp f) = \ex^{\beta(f)}\exp \tau(f)$ and so
\begin{align*}
\sigma\big(\tau(\exp f)\big)\  &=\  \ex^{\beta(f)}\sigma\big(\exp \tau(f)\big)\\ &=\ \ex^{\beta(f) + \alpha(\tau(f))}\exp \sigma\big(\tau(f)\big) \\ &=\ 
\ex^{\beta(f) + \alpha(f)}\exp \sigma\big(\tau(f)\big) 
\end{align*}
using Lemma~\ref{lemobs}  for the last two equalities. So $(c_{\sigma\circ\tau}, \gamma_{\sigma\circ\tau})=(c_\sigma+c_\tau, \gamma_\sigma+\gamma_\tau)$. \qed

\section{Proof of Theorem A}\label{sec:proof of main thm}

\noindent
 We recall here that ``Theorem A'' refers to a result stated in the Introduction. 

\subsection*{Embeddings of $E{\downarrow_n}$ into $\T$}
Let $\sigma\colon E{\downarrow_n}\to \mathbb T$ be a strongly additive differential field embedding. Then $\sigma$ is $\R$-linear, order-preserving, and thus a valued field embedding (of valued subfields  of $\T[[G^{\operatorname{LE}}]]$).
Further analysis of $\sigma$ here will give a hint how to construct such
embeddings in the next subsection. Our main tool is the morphism  of additive groups
$\alpha\colon E{\downarrow_n}\to\R$ from Lemma~\ref{lem:sigma alpha}, given by
$$\alpha(f)\ :=\ \log\big(\sigma(\ex^f)/\ex^{\sigma(f)}\big)\quad  \text{(so $\sigma(\ex^f)=\ex^{\alpha(f)+\sigma(f)}$)}. $$
By that lemma, 
$$\ker\alpha\ =\ \big\{ f\in E{\downarrow}_n: \sigma(\exp f)=\exp{\sigma(f)} \big\} \quad\text{ (so $\R\subseteq \ker \alpha$),}$$
and for $f\in\exp(E{\downarrow_n})$, $r\in\R$,
\begin{equation}\label{eq:sigma log}
\sigma(\log f)=\log\sigma(f)-\alpha(\log f),\quad \sigma(f^r)=\ex^{\alpha(r\log f)-r\alpha(\log f)} \sigma(f)^r.
\end{equation}
As in the proof of Lemma~\ref{lem:alpha} we obtain:

\begin{lemma}\label{lem:extendsigmaT, 0}
  $(E{\downarrow_n})_\preceq\subseteq\ker\alpha$ and so $\alpha(f)=\alpha(f_\succ)$ for $f\in E{\downarrow_n}$.
\end{lemma}

\noindent
By Lemmas~\ref{lem:intersect with Tlog}  and~\ref{lem:sigma is strong}, $\sigma$ restricts to a strong automorphism of the differential subfield $E{\downarrow_n}\cap\T_{\log}=\R[[\ell_0^{\R}\cdots\ell_n^{\R}]]$ of $E{\downarrow_n}$.

\begin{lemma}\label{lemEdown}
The map $\sigma\colon E{\downarrow_n}\to\T$ is healthy and  $\sigma(E{\downarrow_n})\subseteq E{\downarrow_n}$.
\textup{(}Hence every strongly additive differential field embedding $\T\to\T$ is healthy.\textup{)}
\end{lemma}
\begin{proof}
By induction on~$m$  we show that $\sigma(E_{m,n})\subseteq E_{m+n,n}$.
For $m=0$ we have~$\sigma(E_{0,n})=\sigma(\R[[\ell_n^\R]]) \subseteq \R[[\ell_0^{\R}\cdots\ell_n^{\R}]]\subseteq E_{n,n}$ by the remark before  the lemma. Suppose $\sigma(E_{m,n})\subseteq E_{m+n,n}$   for a certain~$m$.
For this inclusion to hold with~${m+1}$ in place of $m$, it is enough by strong linearity of $\sigma$ that $\sigma(G_{m+1,n})\subseteq E_{m+n+1,n}$.
 Now $G_{m+1,n}=G_{m,n}\ex^{A_{m,n}}$, so using the inductive hypothesis we see that it is enough
that $\sigma(\ex^{A_{m,n}})\subseteq  E_{m+n+1,n}$.
For this, let $a\in A_{m,n}$; then $\sigma(a)\in E_{m+n,n}$ by inductive hypothesis, so
$\sigma(\ex^a)=\ex^{\alpha(a)+\sigma(a)}\in E_{m+n+1,n}$.
\end{proof}

\noindent
In the next subsection  we strengthen Lemma~\ref{lemEdown} and prove~$\sigma(E{\downarrow_n}) = E{\downarrow_n}$, and show
that the inverse of $\sigma\colon E{\downarrow_n}\to E{\downarrow_n}$ is also strongly linear.

\begin{lemma} \label{lem:sigmaelln=elln}
The following are equivalent:
\begin{enumerate}
\item[\rm(i)] $\sigma(\ell_n)=\ell_n$;
\item[\rm(ii)] $\sigma(\ell_i)=\ell_i$ for $i=0,\dots,n$;
\item[\rm(iii)] $\sigma(x)=x$ and $\alpha(\ell_i)=0$ for $i=1,\dots,n$;
\item[\rm(iv)] $\sigma(\R^> \cdot G_{m,n})= \R^>\cdot  G_{m,n}$ for all $m$;
\item[\rm(v)] $\sigma(A_{m,n})= A_{m,n}$ for all $m$;
\item[\rm(vi)] $\sigma\big( (E_{m,n})_\succ \big) = (E_{m,n})_\succ$ for all $m$;
\item[\rm(vii)] $\sigma\big( (E{\downarrow_n})_\succ \big) = (E{\downarrow_n})_\succ$.
\end{enumerate}
\end{lemma}
\begin{proof}
The equivalence of (i), (ii), (iii) follows from Lemma~\ref{lem:id on elln} and \eqref{eq:sigma ell_i}.
For the implication~(i)~$\Rightarrow$~(iv), first note
that for every $m$ we have (iv$_m$)~$\Rightarrow$~(v$_m$), where
\begin{enumerate}
\item[(iv$_m$)]  $\sigma(\R^>\cdot G_{m-1,n})= \R^> \cdot G_{m-1,n}$ and $\sigma(\R^> \cdot G_{m,n})= \R^> \cdot G_{m,n}$;
\item[(v$_m$)]  $\sigma(A_{m,n})= A_{m,n}$.
\end{enumerate}
This follows from strong linearity of $\sigma$ and the definition of $A_{m,n}$, by which  
$$A_{m,n} \ =\  \big\{ f\in E_{m,n}:\,  \supp f \succ G_{m-1,n} \big\}.$$
Now assume (i). We prove $\sigma(\R^> \cdot G_{m,n})= \R^> \cdot G_{m,n}$ by induction on $m$.
For   $r\in\R$ we have $(\ell_n^r)^\dagger=r\ell_n^\dagger$, so (1) gives $(\sigma(\ell_n^r))^\dagger=(\ell_n^r)^\dagger$, hence
$\sigma(\ell_n^r)=c\,\ell_n^r$ with $c\in\R^>$, and thus $\sigma(\R^> \cdot G_{0,n})= \R^> \cdot G_{0,n}$.
Assume inductively that $\sigma(\R^> \cdot G_{i,n})= \R^> \cdot G_{i,n}$ for $i=0,\dots,m$; in particular, (iv$_m$) holds. Let $\fm=\fn\ex^a\in G_{m+1,n}$,  $\fn\in G_{m,n}$, $a\in A_{m,n}$.
Then $\sigma(\fm)=\sigma(\fn)\sigma(\ex^a)$ where $\sigma(\fn)\in\R^>G_{m,n}$ by inductive hypothesis, 
and $\sigma(\ex^a)=c\ex^{\sigma(a)}$ where
$c=\ex^{\alpha(a)}\in\R^>$ and $\sigma(a)\in A_{m,n}$ by (iv$_m$)~$\Rightarrow$~(v$_m$); so we have $\sigma(\fm)\in \R^>\cdot G_{m+1,n}$.
Thus $\sigma(\R^> \cdot G_{m+1,n})\subseteq\R^>\cdot G_{m+1,n}$; the reverse inclusion is proved similarly.
This proves~(i)~$\Rightarrow$~(iv). The implication (iv)~$\Rightarrow$~(v) follows from~(iv$_m$)~$\Rightarrow$~(v$_m$),
(v)~$\Rightarrow$~(vi) from $(E_{m,n})_\succ=A_{0,n}\oplus\cdots\oplus A_{m,n}$, and~(vi)~$\Rightarrow$~(vii)
from $(E{\downarrow_n})_\succ = \bigcup_m (E_{m,n})_\succ$.
If (vii) holds, then  
$\sigma(\ell_n)-\ell_n\in \T_{\preceq} \cap (E{\downarrow_n})_\succ$ by Lem\-ma~\ref{lem:sigma elli}. 
Thus (vii)~$\Rightarrow$~(i) in view of  $\T_{\preceq}\cap (E{\downarrow_n})_\succ=\{0\}$.
\end{proof}

\subsection*{Proving Theorem A}
Let an additive map $\alpha\colon \T\to\R$ with $\T_{\preceq}\subseteq\ker\alpha$ and a real number~$c$ be given. 
Recursively  define the elements $\ell_n^*$   of~$\T_{\log}^{>\R}$  by
$$\ell_0^*\ :=\ \ell_0+c,\qquad \ell_{n+1}^*\ :=\ \log\ell_n^* -  \alpha(\ell_{n+1}).$$
Then $(\ell_{n}^*)'=\frac{1}{\ell_0^*\cdots\ell_{n-1}^*}$, and
$\ell_n^*\in \ell_n +  \R[[\ell_0^\Z\cdots\ell_{n-1}^\Z]]_\preceq$
by Lemma~\ref{lem:ellnstar}.

\begin{lemma}\label{lem:extendsigmaT, uniqueness}
Let $\sigma,\tau\colon E{\downarrow_n}\to\T$ be strongly linear field embeddings such that $\sigma(\ell_n^r)=\tau(\ell_n^r)$ for  all
$r\in\R$ and $\sigma(\ex^f)=\ex^{\alpha(f)+\sigma(f)}$, $\tau(\ex^f)=\ex^{\alpha(f)+\tau(f)}$  for all~${f\in E{\downarrow_n}}$. Then $\sigma=\tau$.
\end{lemma}
\begin{proof} Strong linearity of $\sigma$, $\tau$ gives $\sigma(f)=\tau(f)$ for all $f\in E_{0,n}=\R[[\ell_n^\R]]$. Suppose~$\sigma(f)=\tau(f)$ for all $f\in E_{m,n}$. Then for $\fm=\ell_n^r\ex^a$, where $r\in\R$, $a\in (E_{m,n})_\succ$:
$$\sigma(\fm)\ =\ \sigma(\ell_n^r)\ex^{\alpha(a)+\sigma(a)}\ =\ 
\tau(\ell_n^r)\ex^{\alpha(a)+\tau(a)}\ =\ 
\tau(\fm),$$ 
hence $\sigma(f)=\tau(f)$ for $f\in E_{m+1,n}$ by strong linearity of $\sigma$, $\tau$.
\end{proof}

\noindent
We now prove the ``at most one''  part of Theorem A. Let $\sigma$ be any strongly 
additive endomorphism of~$\T$ such that 
\begin{equation}\label{eq:sigma}
\sigma(x)\ =\ x+c\qquad\text{and}\qquad\sigma(\ex^f)\ =\ \ex^{\alpha(f)+\sigma(f)} \quad\text{for all $f\in\T$.}
\end{equation}
Then $\sigma(\ell_n)=\ell_n^*$ by induction on $n$
and  \eqref{eq:sigma ell_i},
hence   $\sigma(\ell_n^r)=\ex^{\alpha(r\ell_{n+1})-r\alpha(\ell_{n+1})}(\ell_n^*)^r$ for 
$r\in\R$ by \eqref{eq:sigmaellir}.  
Thus by the previous lemma, there is at most one strongly additive  endomorphism $\sigma$ of $\T$ satisfying~\eqref{eq:sigma}.
Conversely, for every strongly additive  endomorphism $\sigma$ of $\T$ there is a unique pair~$\alpha$,~$c$ as above such that~\eqref{eq:sigma} holds.

\medskip
\noindent
Next we turn to constructing a strongly additive automorphism of $\T$ as required.

\begin{lemma}\label{lem:extendsigmaT}
Let $\sigma\colon (E_{m,n})_\succ \to (E_{m,n})_\succ$ be a strictly increasing map.
Then there is a unique strongly linear map $\tilde\sigma\colon E_{m+1,n}\to E_{m+1,n}$
such that 
$$\tilde\sigma(\fm)=\ex^{\alpha(\log\fm)-r\alpha(\ell_{n+1})}\ell_n^r\ex^{\sigma(a)}\quad
\text{for $\fm=\ell_n^r\ex^a$, $r\in\R$, $a\in (E_{m,n})_\succ$.}$$
This map $\tilde\sigma$ is injective. If $\sigma$ is bijective, then so is $\tilde\sigma$ and we have~$\tilde\sigma\big( (E_{m+1,n})_\succ\big) = (E_{m+1,n})_\succ$. 
If $\sigma$ is additive, then $\tilde\sigma$ is a field embedding. 
 \end{lemma}
\begin{proof}
Since $E_{m+1,n}=\R[[G_{m+1,n}]]$ where $G_{m+1,n}=\ell_n^\R\ex^{(E_{m,n})_\succ}$ (see Section~\ref{sec:prelims, trans}),
it is clear that there can be at most one such map $\tilde\sigma$.
For $\fm=\ell_n^r \ex^a\in G_{m+1,n}$, where~$r\in\R$, $a\in (E_{m,n})_\succ$, define $\tilde\sigma(\fm):=\ex^{\alpha(\log\fm)-r\alpha(\ell_{n+1})}\ell_n^r\ex^{\sigma(a)}\in\R^>\cdot G_{m+1,n}$. Then
for $a,b\in (E_{m,n})_\succ$, $r,s\in\R$, and $\fm=\ell_n^r\cdot\ex^a,\fn=\ell_n^s\cdot\ex^b$, we have
\begin{align*}
\fm \prec \fn &\quad\Longleftrightarrow\quad \text{$a < b$, or $a=b$ and $r<s$} \\
&\quad\Longleftrightarrow\quad \text{$\sigma(a) < \sigma(b)$, or $\sigma(a)=\sigma(b)$ and $r<s$}
\quad\Longleftrightarrow\quad \tilde\sigma(\fm)\prec\tilde\sigma(\fn).
\end{align*}
Hence $\tilde\sigma\colon G_{m+1,n}\to E_{m+1,n}$ is well-based, so extends (uniquely) to a strongly linear map $E_{m+1,n}\to E_{m+1,n}$ by Proposition~\ref{prop:extend to strongly linear}.
The rest is routine.
\end{proof}

\noindent
The next lemma is obvious. 

\begin{lemma}\label{lem:extendsigmaT, 1}
Let $\sigma$
be a strongly linear field endomorphism of $E_{m,n}$ such that
$\sigma\big( (E_{m,n})_\succ \big) \subseteq (E_{m,n})_\succ$,
and let $\tilde\sigma$ be as in the previous lemma applied to the restriction of $\sigma$ to a map $(E_{m,n})_\succ\to (E_{m,n})_\succ$ in place of $\sigma$. Then $\tilde\sigma$ extends $\sigma$
iff 
$$\sigma(\fm)=\ex^{\alpha(\log\fm)-r\alpha(\ell_{n+1})}\ell_n^r\ex^{\sigma(a)}\quad\text{for all $\ \fm=\ell_n^r\ex^a$,
$r\in\R$, $a\in (E_{m-1,n})_\succ$.}$$
\end{lemma}

\begin{cor}\label{cor:extendsigmaT} There exists a unique strongly linear
field automorphism $\sigma$ of~$E{\downarrow_n}$ with 
$\sigma(\ell_n^r)=\ex^{\alpha(r\ell_{n+1})-r\alpha(\ell_{n+1})}\ell_n^r$ for 
$r\in\R$ and $\sigma(\ex^f)=\ex^{\alpha(f)+\sigma(f)}$ for  $f\in E{\downarrow_n}$.
\end{cor}
\begin{proof}
We inductively construct  strongly linear  auto\-mor\-phisms~$\sigma_m$ of~$E_{m,n}$ such that
for all $m$:
\begin{enumerate}
\item $\sigma_m\big( (E_{m,n})_\succ \big)=(E_{m,n})_\succ$;
\item $\sigma_m(\fm)=\ex^{\alpha(\log\fm)-r\alpha(\ell_{n+1})}\ell_n^r\ex^{\sigma_m(a)}$ for $\fm=\ell_n^r\ex^a$, $r\in\R$, $a\in (E_{m-1,n})_\succ$;  
\item if $m\geq 1$, then
$\sigma_m$ extends $\sigma_{m-1}$.
\end{enumerate}
We let $\sigma_0$ be the strongly linear automorphism of $E_{0,n}=\R[[\ell_n^\R]]$
sending $\ell_n^r$ (${r\in\R}$) to $\ex^{\alpha(r\ell_{n+1})-r\alpha(\ell_{n+1})}\ell_n^r$; then 
clearly (1)--(3) hold for $m=0$. Assume (1)--(3) for a certain $m$.  
Let $\sigma_{m+1}$ be the strongly linear automorphism~$\tilde\sigma$ of $E_{m+1,n}$ obtained
from Lemma~\ref{lem:extendsigmaT} applied to $\sigma_m|_{(E_{m,n})_\succ}$ 
in place of~$\sigma$. By Lemma~\ref{lem:extendsigmaT, 1}, items~(1)--(3) hold with $m+1$ instead of $m$. Clearly the common extension~$\sigma$ of the $\sigma_m$ to a map $E{\downarrow_n}\to E{\downarrow_n}$
is a strongly linear field automorphism of $E{\downarrow_n}$ with~$\sigma(\ell_n^r)=\ex^{\alpha(r\ell_{n+1})-r\alpha(\ell_{n+1})}\ell_n^r$ for 
$r\in\R$. Let $f\in E{\downarrow_n}$, and take $m$ such that~${f\in E_{m,n}}$.
Then~$a:=f_\succ\in (E_{m,n})_\succ$ gives $\sigma(\ex^{a})=\sigma_{m+1}(\ex^{a})= 
\ex^{\alpha(a)+\sigma(a)}$
by~(2) with~$m+1$ in place of $m$. For $b:=f_\preceq=f-a$ we have $\sigma(\ex^b)=\ex^{\sigma(b)}=\ex^{\alpha(b)+\sigma(b)}$ by strong linearity of $\sigma$ and $\ker\alpha\supseteq \T_\preceq$, so $\sigma(\ex^f)=\sigma(\ex^a)\sigma(\ex^b)=\ex^{\alpha(f)+\sigma(f)}$.
This shows the existence of $\sigma$ as required; uniqueness follows from Lem\-ma~\ref{lem:extendsigmaT, uniqueness}. 
\end{proof}

\noindent
In the following we denote the automorphism $\sigma$ of $E{\downarrow_n}$ from Corollary~\ref{cor:extendsigmaT} by~$\sigma_{\alpha,n}$. 
In general, $\sigma_{\alpha,n}$ is not an automorphism of
the {\it differential}\/ field $E{\downarrow_n}$. In order to mend this, we define $\sigma^*_n\colon E{\downarrow_n}\to E{\downarrow_n}$  by
$$\sigma^*_n(f)\  =\  (f{\uparrow_n})\circ \ell_n^*.$$
Note that $\sigma^*_n$ is a strongly linear automorphism of exponential ordered fields, by
 Corollary~\ref{cor:comp} and Lemma~\ref{lem:ellnstarinverse}, with $\sigma^*_n(\ell_n^r)=(\ell_n^*)^r$ for $r\in\R$.
Next we put 
$$\sigma_n\ :=\ \sigma_n^*\circ \sigma_{\alpha,n},$$
a strongly linear field automorphism of $E{\downarrow_n}$. We have 
\begin{equation}\label{eq:sigmanellnr}
\sigma_n(\ell_n^r)\	=\  \ex^{\alpha(r\ell_{n+1})-r\alpha(\ell_{n+1})} (\ell_n^*)^r\qquad\text{ for $r\in\R$,}
\end{equation}
 in particular $\sigma_n(\ell_n)=\ell_n^*$,
and
$$\sigma_n(\ex^f)\ =\ \ex^{\alpha(f)+\sigma_n(f)}\qquad\text{ for $f\in E{\downarrow_n}$.}$$  
Moreover:

\begin{lemma}\label{lem:extendsigmaT, 2}
$\sigma_n(\ell_i)=\ell_i^*$ for $i=0,\dots,n$.
\end{lemma}
\begin{proof}
By induction on $j=0,\dots,n$ we prove $\sigma_n(\ell_{n-j})=\ell_{n-j}^*$.
We already discussed the case $j=0$.
Suppose $j<n$ and we have   shown~$\sigma_n(\ell_{n-j})=\ell_{n-j}^*$. Then
\begin{align*}
\sigma_n(\ell_{n-j-1})\ &=\ \sigma_n\big(\!\exp(\ell_{n-j})\big)\\
&=\ \exp\big(\alpha(\ell_{n-j})+\sigma_n(\ell_{n-j})\big)\\
&=\ \exp\big(\alpha(\ell_{n-j})+\ell_{n-j}^*\big)\\
&=\ \exp\big(\alpha(\ell_{n-j})+\log\ell_{n-j-1}^*- \alpha(\ell_{n-j})\big)=\ell_{n-j-1}^*
\end{align*}
as required.
\end{proof}

\noindent
Combining Lemmas~\ref{lem:id on elln} and \ref{lem:extendsigmaT, 2} yields:

\begin{cor}\label{cor:id on elln}  
The following are equivalent:
\begin{enumerate}
\item[\rm(i)] $\sigma_n(\ell_n)=\ell_n$;
\item[\rm(ii)] $\sigma_n(f)=f$ for all $f\in \R[[\ell_0^\Q\cdots\ell_n^\Q]]$;
\item[\rm(iii)] $c=0$ and $\alpha(\ell_i)=0$ for $i=1,\dots,n$.
\end{enumerate}
\end{cor}

\noindent
Next we show that $\sigma_n$ is an automorphism of the {\em differential\/}  field $E{\downarrow_n}$:

\begin{lemma}\label{lem:extendsigmaT, 3}
For all $f\in E{\downarrow_n}$ we have $\sigma_n(f')=\sigma_n(f)'$.
\end{lemma}
\begin{proof}
Let $r\in\R$; then \eqref{eq:sigmanellnr} and
 arguing as in the proof of Proposition~\ref{prop:SigmaAut of Tlog} shows~$\sigma_n\big( (\ell_n^r)' \big) = \sigma_n(\ell_n^r)'$. 
By strong linearity of $\sigma_n$ and the derivation of $E{\downarrow_n}$, this gives $\sigma_n(f')=\sigma_n(f)'$ for all $f\in E_{0,n}=\R[[\ell_n^\R]]$.
Assume for a certain $m$ that for all $f\in E_{m,n}$ we have $\sigma_n(f')=\sigma_n(f)'$.  
Then for all $a\in (E_{m,n})_\succ$:
\begin{align*}
\sigma_n(\ex^a)'\  &=\  \big(\ex^{\alpha(a)+\sigma_n(a)}\big)' \\
&=\ \ex^{\alpha(a)+\sigma_n(a)}\cdot \sigma_n(a)' \\
&=\ \ex^{\alpha(a)+\sigma_n(a)}\cdot\sigma_n(a') 
 =\ \sigma_n(\ex^{a}a')\ =\ \sigma_n \big( (\ex^a)'\big).
\end{align*}
Using $\sigma_n\big( (\ell_n^r)' \big) = \sigma_n(\ell_n^r)'$ for $r\in\R$ and
strong linearity as before, we get $\sigma_n(f')=\sigma_n(f)'$ for all $f\in E_{m+1,n}$.
\end{proof}

\begin{cor}
The automorphism 
$\sigma_{n+1}$ of $E{\downarrow_{n+1}}$  extends  $\sigma_n$.
\end{cor}
\begin{proof} Use $\ell_n^r=\ex^{r\ell_{n+1}}$, 
Lemma~\ref{lem:extendsigmaT, 2},  ~\eqref{eq:sigmanellnr}, and $\ell_{n+1}^*=\log \ell_n^* -\alpha(\ell_{n+1})$ to obtain 
$$\sigma_{n+1}(\ell_n^r)\ =\ \ex^{\alpha(r\ell_{n+1})-r\alpha(\ell_{n+1})} (\ell_n^*)^r\  =\  \sigma_n(\ell_n^r)\qquad (r\in \R).$$ 
 Now appeal to Lemma~\ref{lem:extendsigmaT, uniqueness}.  
\end{proof}

\noindent
The common extension  of the $\sigma_n$ to a map $\sigma\colon\T\to \T$
is a strongly additive  automorphism of the differential field $\T$ satisfying~\eqref{eq:sigma}. This finishes the proof of  Theorem A. \qed

\medskip
\noindent
Let $\sigma_{\alpha,c}$ be the
strongly additive automorphism $\sigma$ of~$\T$ satisfying~\eqref{eq:sigma}. 
With~$0$ denoting the zero map $\T\to\R$ we have $\sigma_{0,0}=\id_{\T}$. If $\beta\colon\T\to\R$ is also an additive map with $\T_{\preceq}\subseteq\ker\beta$ and $d\in\R$, then 
$\beta\circ \sigma_{\alpha,c}+\alpha$ is an additive map $\T\to\R$ with~$\T_\preceq$ in its kernel 
(since $\sigma_{\alpha,c}(\T_\preceq)=\T_\preceq$)
and
\begin{equation}\label{eq:composition}
\sigma_{\beta,d}\circ\sigma_{\alpha,c}\ =\ \sigma_{\gamma,c+d}\qquad\text{ where $\gamma:=\beta\circ\sigma_{\alpha,c}+\alpha$.}
\end{equation}
Hence if $\beta=-\alpha\circ\sigma_{\alpha,c}^{-1}$ then $\sigma_{\beta,-c}\circ\sigma_{\alpha,c}=\id_{\T}$, that is, 
\begin{equation}\label{eq:inversion}
\sigma_{\alpha,c}^{-1}\ =\ \sigma_{\beta,-c} \qquad\text{ where $\beta:=-\alpha\circ\sigma_{\alpha,c}^{-1}$.}
\end{equation}
In particular:

\begin{cor}
The inverse of  a strongly additive automorphism of~$\T$  is also strongly additive.
\end{cor}

\noindent
That is, each strongly additive automorphism of~$\T$ is strong.
 Recall that the  strong automorphisms of~$\T$ form 
 a sub\-group~$\SigmaAut_\der(\T)$ of the group~$\Aut_\der(\T)$ of differential field automorphisms of   $\T$
under composition,   with identity element $1=\id_{\T}$. 
By remarks after the proof of Lemma~\ref{lem:extendsigmaT, uniqueness}, every
strongly additive endomorphism of $\T$ is of the form $\sigma_{\alpha,c}$, and is thus an automorphism of $\T$. 

\begin{cor}\label{cor:extend from Euparrown}
Every strongly additive differential field 
embedding $E{\downarrow_n}\to\T$ extends to a strong automorphism of $\T$ and maps $E{\downarrow_n}$ onto  $E{\downarrow_n}$.
\end{cor}
{\sloppy
\begin{proof}
Let $\sigma\colon E{\downarrow_n}\to\T$ be a strongly additive  differential field embedding.
Let~$\alpha\colon E{\downarrow_n}\to\R$ be as at the beginning of this section and $c:=\sigma(x)-x\in\R$. We extend~$\alpha$ to an additive map $\T\to\R$ (also called $\alpha$) vanishing on $\T_\preceq$, such that 
$$\alpha(r\ell_{n+1})\ :=\ \log\big(\sigma(\ell_n^r)/\sigma(\ell_n)^r\big)\qquad (r\in \R).$$
Let $r\in \R$. Then $\sigma_{\alpha,c}(\ell_n^r)=\sigma_{\alpha,c}(\ex^{r\ell_{n+1}})=
\ex^{\alpha(r\ell_{n+1})}\cdot \ex^{\sigma_{\alpha,c}(r\ell_{n+1})}$. Moreover, by Lem\-ma~\ref{lem:extendsigmaT, 2}, 
$$\sigma_{\alpha,c}(r\ell_{n+1})\ =\ r\sigma_{\alpha,c}(\ell_{n+1})=r\ell_{n+1}^*=r\big(\log \ell_n^* - \alpha(\ell_{n+1})\big),$$ and thus
$$\sigma_{\alpha,c}(\ell_n^r)=\ex^{\alpha(r\ell_{n+1})-r\alpha(\ell_{n+1})}\ell_n^*r=\sigma(\ell_n^r),$$ using \eqref{eq:sigmanellnr} for the last equality. Hence~$\sigma_{\alpha,c}$ extends~$\sigma$ by Lemma~\ref{lem:extendsigmaT, uniqueness}.
In the course of the proof of Theorem A we saw that $\sigma_{\alpha,c}(E{\downarrow_n})=E{\downarrow_n}$ for all~$n$.
\end{proof}}

\noindent
Thus $\SigmaAut_{\der}(\T|E{\downarrow_n})$ is a normal subgroup of $\SigmaAut_{\der}(\T)$, by Lemma~\ref{lemgalois}. 

\medskip
\noindent
 For $n\ge m$, each $\sigma\in \SigmaAut_\der(E{\downarrow}_n)$
restricts to a strongly linear automorphism~$\pi^n_{m}(\sigma)$ of the differential field $E{\downarrow}_m$, by Corollary~\ref{cor:extend from Euparrown}.
The restriction maps $$\pi^n_{m}\ \colon\  
\SigmaAut_\der(E{\downarrow}_n)\to \SigmaAut_\der(E{\downarrow}_m)$$ are group morphisms,
and $(\pi^n_{m})_{n\ge m}$ is an inverse system of group morphisms. Clearly, the induced map 
from $ \SigmaAut_\der(\T)$ into the inverse limit of this system is a group isomorphism.

\medskip
\noindent
Combining Corollaries~\ref{cor:extend from Euparrown} and ~\ref{cor:id on elln}  with Lemma~\ref{lem:sigmaelln=elln}  yields:

 \begin{cor} \label{cor:sigmaelln=elln}  The following conditions on $\sigma=\sigma_{\alpha,c}$ are equivalent: \begin{enumerate}
 \item[\rm(i)] $\sigma(\ell_n)=\ell_n$ for all $n$;
 \item[\rm(ii)] $\sigma$ is the identity on $\T_{\log}^\Q$; 
 \item[\rm(iii)] $c=0$ and $\alpha(\ell_n)=0$ for all $n\ge 1$; 
 \item[\rm(iv)] $\sigma(\R^>\cdot G^{\operatorname{LE}})=\R^>\cdot G^{\operatorname{LE}}$;
\item[\rm(v)] $\sigma(\T_\succ)=\T_\succ$.
\end{enumerate}
\end{cor}

\noindent
(For example (i)~$\Rightarrow$~(iv) follows from (i)~$\Rightarrow$~(iv) in Lemma~\ref{lem:sigmaelln=elln}, (iv)~$\Rightarrow$~(v) holds by strong linearity of $\sigma$, and (v)~$\Rightarrow$~(i) follows from Corollary~\ref{cor:extend from Euparrown} and (vii)~$\Rightarrow$~(i) in Lem\-ma~\ref{lem:sigmaelln=elln}.) 
From (i)~$\Rightarrow$~(iv) of Corollary~\ref{cor:sigmaelln=elln} we obtain:

\begin{cor}\label{cortr} Let $\mathcal{G}$ be a subset of $\SigmaAut_\der(\T)$ such that $\sigma(\ell_n)=\ell_n$ for all $\sigma\in \mathcal{G}$ and all $n$. Then its fixed field $\T^{\mathcal{G}}$ is truncation closed in $\R[[G^{\operatorname{LE}}]]$.
\end{cor} 

\noindent
Note that for $\mathcal{G}$ as in Corollary~\ref{cortr}, every $f\in \T$ has by Corollary~\ref{corbestapp} a best approximation in 
$\T^{\mathcal{G}}$; in other words, there are no divergent pc-sequences in $\T^{\mathcal{G}}$ with a pseudolimit in $\T$. 
Next we show that there are no non-trivial
``small perturbations of the identity'' in $\SigmaAut_\der(\T)$, thus putting
the automorphisms studied here outside the purview of \cite{BKKPS}:

\begin{cor}
Let $\sigma\in\SigmaAut_\der(\T)$ and $\sigma(\fm) \asymp \fm$ for all $\fm\in G^{\operatorname{LE}}$. Then $\sigma=\id_{\T}$.
\end{cor}
\begin{proof}
Take $\alpha$, $c$ such that $\sigma=\sigma_{\alpha,c}$. Then
$$\ex^{x^2}\ \asymp\ \sigma(\ex^{x^2})\ =\ \ex^{\alpha(x^2)+\sigma(x^2)}\ =\ \ex^{\alpha(x^2)+c^2}\cdot \ex^{2cx}\ex^{x^2}$$
and $\ex^{2cx}\in G^{\operatorname{LE}}$, so $c=0$.
By induction on $n$ we show that $\alpha(\ell_m)=0$ for~$m=1,\dots,n$. This is vacuously true for $n=0$.
Suppose it holds for a certain $n$.  Using~$c=0$ and $\sigma(\ell_m)=\ell_m^*$  this yields $\sigma(\ell_m)=\ell_m$ for $m=0,\dots,n$. Then 
$$\sigma(\ell_{n+1})\ =\ \ell_{n+1}^*\ =\ \ell_{n+1}-\alpha(\ell_{n+1}),$$
so $\sigma\big(\ell_{n+1}^2\big)= \alpha(\ell_{n+1})^2- 2\alpha(\ell_{n+1})\ell_{n+1}+ \ell_{n+1}^2$ and thus
$$\ex^{\ell_{n+1}^2}\ \asymp\ \sigma(\ex^{\ell_{n+1}^2})\ =\ \ex^{\alpha(\ell_{n+1}^2)+\alpha(\ell_{n+1})^2}\cdot \ex^{-2\alpha(\ell_{n+1})\ell_{n+1}}\ex^{\ell_{n+1}^2},$$
so $\alpha(\ell_{n+1})=0$. This finishes the induction. So $\sigma(\T_{\succ})= \T_{\succ}$ by Corollary~\ref{cor:sigmaelln=elln}.
Let $a\in \T_{\succ}$. Then $\fm:=\ex^a\in G^{\operatorname{LE}}$ and  $\sigma(\fm)=\ex^{\alpha(a)} \ex^{\sigma(a)}$ with $\sigma(a)\in \T_{\succ}$, so~$\ex^{\sigma(a)}\in G^{\operatorname{LE}}$, which in view of $\sigma(\fm)\asymp \fm$ gives $\sigma(a)=a$. Thus $\sigma$ is the identity on $\T_{\succ}$, in particular, $\sigma(\fn)=\fn$ for all $\fn\succ 1$ in $G^{\operatorname{LE}}$, and hence $\sigma(\fn)=\fn$ for all~$\fn\in G^{\operatorname{LE}}$.
Strong linearity of $\sigma$ then gives $\sigma=\id_\T$.
\end{proof}

\section{The Group $\SigmaAut_\der(\T)$}\label{sec:group}

\noindent
In this section we show how $\SigmaAut_\der(\T)$ is an internal semidirect product $\mathcal{E}\rtimes \mathcal{T}$
of a normal subgroup $\mathcal{E}$ with an abelian subgroup $\mathcal{T}$. 

\subsection*{The exponential torus of $\T$}
Let $\Hom(\T,\R)$ be the $\R$-module of morphisms from the additive group $\T$ to the additive group $\R$,
and $$D\ :=\ \big\{\alpha\in \Hom(\T,\R):\, \ker \alpha \supseteq \T_{\preceq}\big\},$$  a submodule of $\Hom(\T,\R)$.
Then $\alpha\mapsto \alpha|_{\T_\succ}$ is an isomorphism from $D$ onto the $\R$-module $\Hom(\T_\succ,\R)$ of additive maps $\T_\succ\to\R$.
For $\alpha\in D$, set $$\sigma_\alpha\ :=\ \sigma_{\alpha,0}\in \SigmaAut_\der(\T),$$ so 
$\sigma_\alpha(x)=x$ and $\sigma_\alpha(\ex^f)=\ex^{\alpha(f)+\sigma_\alpha(f)}$ for all $f\in\T$. 
By rules \eqref{eq:composition} and~\eqref{eq:inversion} for composition and inversion in $\SigmaAut_\der(\T)$,   
$$\mathcal E\ :=\ \big\{\sigma\in \SigmaAut_{\der}(\T):\ \sigma(x)=x\big\}\ =\ \{\sigma_\alpha:\alpha\in D\}$$ is a subgroup of $\SigmaAut_\der(\T)$, which we call the {\it exponential torus}\/ of $\T$.
Denote the inverse of the bijection $\alpha\mapsto\sigma_\alpha\colon D\to\mathcal E$ by $\sigma\mapsto \alpha_\sigma$. Then 
\begin{equation}\label{eq:cocycle}
\alpha_{\sigma\circ\tau}\  =\  \alpha_\sigma\circ\tau+\alpha_\tau\quad\text{for $\sigma,\tau\in\mathcal E$.}
\end{equation}
So for $\sigma,\tau\in\mathcal E$ we have 
$$\sigma\circ\tau\ =\ \tau\circ\sigma\qquad\Longleftrightarrow\qquad \alpha_\sigma\circ\tau+\alpha_\tau\ =\  \alpha_\tau\circ\sigma+\alpha_\sigma.$$
The group~$\mathcal E$ is  not abelian: take  $r\ne s$ in $\Q^{>}$, and $\alpha,\beta\in D$ such that
$\alpha(x)=\log r$, $\beta(x)=\log s$, and $\alpha(\ex^x)=\beta(\ex^x)=1$.
Then with $\sigma=\sigma_\alpha$, $\tau=\sigma_\beta$,
$$(\alpha_\sigma\circ\tau+\alpha_\tau)(\ex^x)\ =\ 
\alpha(\ex^{\beta(x)+\tau(x)})+\beta(\ex^x)\ =\ \alpha(s\ex^x)+1\ =\ s+1
$$
and similarly $(\alpha_\tau\circ\sigma+\alpha_\sigma)(\ex^x)=r+1$. Hence $\sigma\circ\tau \neq \tau\circ\sigma$.

\medskip\noindent
Thus the bijection $\sigma\mapsto \alpha_\sigma\colon \mathcal{E}\to D$ is not a group morphism; rather, by \eqref{eq:cocycle}  this map is a
$1$-cocycle of the opposite group $\mathcal E^{\operatorname{opp}}$ of $\mathcal E$ in $D$ with respect to the natural right action $(\alpha,\sigma)\mapsto \alpha\circ\sigma$ of
$\mathcal E$ on $D$. (See \cite[Chapter~VI, \S{}10]{Lang}.)
Let $\operatorname{Aff}(D)$ be the group of automorphisms of the affine space $D$ over $\R$.
For $\sigma\in\mathcal E$ we have the $\R$-linear bijection $\alpha\mapsto \alpha\circ \sigma\colon D\to D$ with inverse $\beta\mapsto \beta\circ \sigma^{-1}$, and so we obtain~$\hat\sigma\in\operatorname{Aff}(D)$ given by $\alpha\mapsto \alpha\circ\sigma+\alpha_\sigma$.

\begin{lemma}
The map $\sigma\mapsto\hat\sigma\colon \mathcal E^{\operatorname{opp}}\to\operatorname{Aff}(D)$ is a group embedding. Moreover, the \textup{(}left\textup{)} action $(\sigma,\alpha)\mapsto \hat\sigma(\alpha)$ of
$\mathcal E^{\operatorname{opp}}$ on $D$ is regular: for all~$\alpha,\beta\in D$ there is a unique $\sigma\in\mathcal E$ with $\hat\sigma(\alpha)=\beta$.
\end{lemma}
\begin{proof}
Let $\alpha$, $\beta$ range over $D$ and $\sigma$, $\tau$   over $\mathcal E$.
Then 
\begin{align*}
\hat{\sigma\circ\tau}(\alpha)\ 	&=\ \alpha\circ(\sigma\circ\tau)+\alpha_{\sigma\circ\tau}\\
								&=\  (\alpha\circ \sigma)\circ\tau + \alpha_\sigma\circ\tau+\alpha_\tau \\
								&=\   (\alpha\circ\sigma+\alpha_\sigma )\circ\tau + \alpha_\tau \\
								&=\ \hat\tau(\alpha\circ\sigma+\alpha_\sigma)=(\hat\tau\circ\hat\sigma)(\alpha).
\end{align*}
We have $\hat\sigma(0)=\alpha_\sigma$, hence
$\sigma=\id_{\T}$ iff $\hat\sigma=\id_D$. This proves the first statement.
For given $\alpha$, $\beta$ take $\sigma$, $\tau$ with $\alpha_\sigma=\alpha$, $\alpha_\tau=\beta$; then 
$\hat\sigma(0)=\alpha$, $\hat\tau(0)=\beta$, so
with $\rho:=\sigma^{-1}\circ\tau$ we have 
$\hat\rho(\alpha)=\hat \tau(\hat\sigma^{-1}(\alpha))=\hat\tau(0)=\beta$.  Thus the 
action of~$\mathcal E^{\operatorname{opp}}$ on $D$ is transitive. Now suppose $\hat\sigma(\alpha)=\alpha$, and
take $\tau$ with $\alpha_\tau=\alpha$; then~$\alpha_{\tau\circ\sigma}=\hat{\tau\circ\sigma}(0)=\hat\sigma(\hat\tau(0))=\hat\sigma(\alpha)=\alpha=\alpha_\tau$ and so $\tau\circ\sigma=\tau$, which yields~$\sigma=\id_{\T}$.
This proves regularity.
\end{proof}

\begin{cor}\label{cor:tfree}
The group $\mathcal E$ is torsion-free.  
\end{cor}
\begin{proof}
Every $\varphi\in\operatorname{Aff}(D)$   with $\varphi^n=\id_D$ ($n\ge 1$) has  fixed point $\frac{1}{n}\sum_{i=0}^{n-1} \varphi^i(0)$. Now use the previous lemma.
\end{proof}

\noindent
The group $\mathcal E$ is not simple: the following lemma says that
for $a\in (E_0)^{\ne}_\succ$,   the map~$\sigma\mapsto \alpha_\sigma(a)\colon\mathcal E \to \mathbb R$
is a group morphism whose kernel is a nontrivial normal proper subgroup of $\mathcal E$. 

\begin{lemma}
Let $\sigma,\tau\in\mathcal E$. Then for $a\in E_0$ and $f\in E_1$ we have
$$\alpha_{\sigma\circ\tau}(a)\ =\ \alpha_\sigma(a)+\alpha_\tau(a), \qquad
(\sigma\circ\tau)(f)\ =\ (\tau\circ\sigma)(f).$$
\end{lemma}
\begin{proof} By Lemma~\ref{lem:sigma alpha}(2) we have $\sigma(x^r)=x^r$ for $r\in \R$, so $\sigma$ is the identity on~$E_0=\R[[x^\R]]$;
likewise for $\tau$. Using \eqref{eq:cocycle} this yields the equality on the left.  Since $\sigma$, $\tau$ are the identity on $E_0$ 
we also obtain for $a\in (E_0)_{\succ}$: 
$$(\sigma\circ \tau)(\ex^a)\ =\ \ex^{\alpha_{\sigma\circ \tau}(a)}\cdot \ex^a,\qquad (\tau\circ \sigma)(\ex^a)\ =\ \ex^{\alpha_{\tau\circ \sigma}(a)}\cdot \ex^a,$$
so $\sigma\circ \tau$ and $\tau\circ \sigma$ agree on $\ex^{(E_0)_{\succ}}$. Using strong linearity and $G_1=x^{\R}\ex^{(E_0)_{\succ}}$ this yields the equality on the right. 
\end{proof}

\begin{lemma}
The cardinality of $\mathcal E$ is $2^{2^{\aleph_0}}$.
\end{lemma}
\begin{proof}
This equals the cardinality of $\Hom(\T_\succ,\R)$. Each additive map $\T_\succ\to\R$ is $\Q$-linear, and 
so choosing a basis $\mathcal B$ of $\T_\succ$ as $\Q$-linear space we have an obvious isomorphism
$\Hom(\T_\succ,\R)\cong \R^{\mathcal B}$. Now use $|\T_\succ|=|\T|=2^{\aleph_0}$
by \cite[\S{}2.10]{DMM2}.
\end{proof}

\subsection*{The translation group of $\T$} 
Every real number $c$ yields the strong  auto\-mor\-phism~$f(x)\mapsto f(x+c)$ of the differential field $\T$, to be
denoted by $\sigma_c$, so $\sigma_c=\sigma_{0,c}$ in the notation of Section~\ref{sec:proof of main thm}. 
The map $c\mapsto\sigma_c$ is an embedding of the additive group~$\R$
into the group~$\SigmaAut_\der(\T)$, so
we obtain an abelian subgroup
$$\mathcal T\ :=\ \{ \sigma_c:c\in\R \}$$
of~$\SigmaAut_\der(\T)$. We call $\mathcal T$ the {\it translation group}\/ of  $\T$.
Next we determine the fixed field of $\mathcal T$.

\begin{lemma}\label{lem:fixed field of M}
Let $f\in\T$, $f\notin \R$, $g\in\T^{>\R}$ and $f\circ g=f$. Then $g=x$.
\end{lemma}
\begin{proof}
Replacing $f$ by $-f$, subtracting a constant, and taking a multiplicative inverse, if necessary, we reduce to the case that $f>\R$. Then $g\in\T^{>\R}$ and $f\circ g=f$ gives
$g=f^{\inv}\circ(f\circ g)=f^{\inv}\circ f=x$.
\end{proof}

\begin{cor}\label{cor:R definably closed}
For each $\sigma\in\mathcal T$ with $\sigma\neq 1$ we have $\T^{\sigma}=\R$. 
\end{cor}

\noindent
Thus $\R$ is definably closed in $\T$. Therefore the constant field $C$ of any differential field $K$ elementarily equivalent to $\T$ is definably closed in $K$.

\subsection*{Structure of $\SigmaAut_\der(\T)$}
To describe how $\mathcal E$ and $\mathcal T$ fit together in $\SigmaAut_\der(\T)$, let~${\alpha\in D}$, $c\in\R$.
Then $\sigma_{\alpha,c}=\sigma_c\sigma_{\alpha}\in \mathcal T\mathcal E=\mathcal E\mathcal T$ by \eqref{eq:composition}.
Thus $\SigmaAut_\der(\T)=\mathcal E\mathcal T$ with $\mathcal E\cap\mathcal T=\{1\}$, and moreover: 

\begin{lemma}\label{lem:T normal} For $\alpha_c:=\alpha \circ\sigma_c\in D$ we have
 $\sigma_c^{-1}\sigma_\alpha\sigma_c=\sigma_{\alpha_c}$.
\textup{(}Hence~$\mathcal E$ is a normal subgroup of $\SigmaAut_\der(\T)$.\textup{)} 
\end{lemma}
\begin{proof}
We have
$(\sigma_c^{-1}\sigma_\alpha\sigma_c)(x)=x$, and
for $f\in\T$ we have
$$(\sigma_c^{-1}\sigma_\alpha\sigma_c)(\ex^f) = (\sigma_c^{-1}\sigma_\alpha) (\ex^{\sigma_c(f)}) =
\sigma_c^{-1}(\ex^{\alpha_c(f)+(\sigma_\alpha\sigma_c)(f)})=\ex^{\alpha_c(f)+(\sigma_c^{-1}\sigma_\alpha\sigma_c)(f)}$$
and thus $\sigma_c^{-1}\sigma_\alpha\sigma_c=\sigma_{\alpha_c}$ by the uniqueness part of the main 
theorem. 
\end{proof}

\begin{cor}\label{cor:semidirect product} The group~$\SigmaAut_\der(\T)$ is the internal semidirect product $\mathcal{E}\rtimes \mathcal{T}$ of its subgroup $\mathcal T$ by its normal subgroup~$\mathcal E$ \textup{(}according to Bourbaki's terminology~\cite[Chapter~I, \S{}6.1]{Bourbaki}\textup{)}.
\end{cor}

\noindent
From Corollary~\ref{cor:semidirect product} and  torsion-freeness of $\mathcal E$ and $\mathcal T$ (Corollary~\ref{cor:tfree}) we obtain:

\begin{cor}\label{cor:tfreeS}
The group $\SigmaAut_\der(\T)$ is torsion-free. 
\end{cor}





\subsection*{The field of exponents}
Given an automorphism $\sigma$ of $\T$ we say that a real number~$r$ is an {\it exponent}\/ of~$\sigma$
if  $\sigma(f^r)=\big(\sigma(f)\big){}^r$ for all $f\in\T^>$.
The set of exponents of $\sigma$ forms a subfield of $\R$, the {\it field of exponents}\/ of $\sigma$.
Given $\mathcal G\subseteq \mathcal E$, the field of exponents of $\mathcal G$
is defined to be the intersection of the fields of exponents of all $\sigma\in\mathcal G$.
Thus the field of exponents of~$\mathcal T$ is $\R$. The field of exponents of $\mathcal E$ is $\Q$: for $r\in \R\setminus \Q$ and $\alpha\in D$ with $\alpha(x)=0$, $\alpha(rx)\ne 0$ we have 
$$\sigma_{\alpha}(\ex^x)\ =\ \ex^x, \qquad \sigma_{\alpha}\big((\ex^x)^r\big)\ =\ \sigma_{\alpha}(\ex^{rx})\ =\ \ex^{\alpha(rx)}\ex^{rx}\ \ne\  \sigma_{\alpha}(\ex^x)^r.$$
If $\k$ is the field of exponents of $\mathcal G$ and $r\in \k$, then 
the fixed field $\T^{\mathcal G}$ of~$\mathcal G$ is closed under $r$-th powers:
if $f\in\T^{\mathcal G}$ and $f>0$, then $f^r\in\T^{\mathcal G}$.
Given $\alpha\in D$,  a {\em scalar of~$\alpha$\/} is an $r\in \R$ such that
$\alpha(rf)=r\alpha(f)$ for all  $f\in\T$. Thus~$r$ is an exponent of $\sigma\in\mathcal E$ iff it is a scalar of $\alpha_\sigma$, by \eqref{eq:sigma log}.
The scalars of any $\alpha\in D$ form a  subfield $\k$ of $\R$, and the map $\alpha\colon\T\to\R$ is
$\k$-linear. 

\subsection*{Power-preserving automorphisms}
Let $\k$ be a subfield of $\R$. We say that an automorphism  of $\T$ is {\it $\k$-power preserving}\/
if its field of exponents contains~$\k$.
The $\k$-power preserving $\sigma\in\SigmaAut_\der(\T)$  form a   subgroup $\SigmaAut_\der^\k(\T)$ of $\SigmaAut_\der(\T)$ with~$\mathcal T\subseteq \SigmaAut_\der^\k(\T)$.
  Setting $\mathcal E^\k:=\SigmaAut_\der^\k(\T)\cap\mathcal E$ we have
$$\mathcal E^\k\ =\ \{\sigma_\alpha:\alpha\in D^\k\big\},\ \text{ where }\  D^\k := \{ \alpha\in D: \text{$\alpha$ is $\k$-linear} \}.$$
The group $\SigmaAut_\der^\k(\T)$ is the internal semidirect product $\mathcal E^\k\rtimes \mathcal T$ of its subgroup~$\mathcal T$ by its normal subgroup $\mathcal E^\k$. 
We have $$\SigmaAut_\der^\Q(\T)\ =\ \SigmaAut_\der(\T),\quad
\mathcal E^\Q\ =\ \mathcal E,\quad D^\Q\ =\ D,$$
and if $\k$ is a subfield of the subfield $\mathbf l$ of $\R$, then 
 $$\SigmaAut_\der^{\k}(\T)\ \supseteq\  \SigmaAut_\der^\l(\T),\quad
\mathcal E^\k\ \supseteq\ \mathcal E^\l,\quad D^\k\ \supseteq\ D^\l.$$
For $F\subseteq\T$ we  set
$$\SigmaAut_\der^\k(\T|F)\ :=\ \SigmaAut_\der(\T|F)\cap \SigmaAut_\der^\k(\T), \qquad \mathcal E^\k|F\ :=\ \SigmaAut_\der(\T|F)\cap \mathcal E^\k.$$

\begin{lemma}\label{lemt1} Let $\sigma\in \mathcal{E}^\R$. Then $\sigma(E_{m,1}) = E_{m,1}$ for all $m$.
\end{lemma}
\begin{proof} Since $\mathcal{E}^\R$ is a subgroup of $\SigmaAut_\der(\T)$, it suffices to prove $\sigma(E_{m,1}) \subseteq  E_{m,1}$. We proceed by induction  on $m$. We have $E_{0,1}=\R[[\ell_1^{\R}]]$, and $\sigma(\ell_0)=\ell_0$, so  with
$\alpha\in D^{\R}$ such that $\sigma=\sigma_{\alpha}$ we obtain
 $\sigma(\ell_1)=\ell_1-\alpha(\ell_1)\in E_{0,1}$, and thus $\sigma(E_{0,1})\subseteq E_{0,1}$. 
 
 Suppose  $\sigma(E_{m,1})\subseteq E_{m,1}$ for a certain $m$. We have $E_{m+1,1}=\R[[G_{m+1,1}]]$, so to obtain $\sigma(E_{m+1,1})\subseteq E_{m+1,1}$ it is enough to show $\sigma(G_{m+1,1})\subseteq E_{m+1,1}$. Now $G_{m+1,1}=G_{m,1}\ex^{A_{m,1}}$, so
$\sigma(G_{m+1,1})=\sigma(G_{m,1})\sigma(\ex^{A_{m,1}})$. We can assume  inductively that $\sigma(G_{m,1})\subseteq E_{m,1}$ and that for $a\in A_{m,1}$ we have
 $\sigma(a)\in E_{m,1}$, so~$\sigma(\ex^a)=\ex^{\alpha(a)}\ex^{\sigma(a)}\in E_{m+1,1}$.
\end{proof}

\subsection*{Action by composition on the exponential torus} 
For $\alpha\in\Hom(\T,\R)$ and ${g\in\T^{>\R}}$, let $g\alpha$ be the group morphism $\T\to\R$ given by $g\alpha(f)=\alpha(f\circ g)$ for $g\in\T$. 
The map $(g,\alpha)\mapsto g\alpha$ is a left action of the group~$\T^{>\R}$ (under composition)
on the set~$\Hom(\T,\R)$: $x\alpha=\alpha$ and $(g\circ h)\alpha=g(h\alpha)$ for $\alpha\in\Hom(\T,\R)$, $g,h\in\T^{>\R}$.

Let $g\in\T^{>\R}$. Then $\ker g \alpha=(\ker\alpha)\circ g^{\inv}$ for  $\alpha\in \Hom(\T,\R)$, and
 $\alpha\mapsto g\alpha$ is an automorphism of the group $\Hom(\T,\R)$
with inverse $\beta\mapsto g^{\inv}\beta$, and for any  subfield $\k$ of $\R$ it maps the subgroup $D^\k$ of $\Hom(\T,\R)$ onto itself.
Note also that  for $\alpha, \beta\in \Hom(\T,\R)$ we have  $g(\alpha+\beta)=g\alpha+g\beta$.

\begin{lemma}\label{faithful} 
The action of $\T^{>\R}$ on $D^\R$ is faithful, that is, if $g\in\T^{>\R}$ and $g\alpha=\alpha$ for all $\alpha\in  D^{\R}$, then
$g=x$.
\end{lemma}
\begin{proof}
Let $g\in \T^{>\R}$, $g\neq x$. 
It suffices to find $f\in\T_\succ$ such that $f\circ g\notin\R f+\T_{\preceq}$,  since that would give
an $\alpha\in D^\R$ with $\alpha(f\circ g) \neq \alpha(f)$.
Let $0\ne f\in\T_{\succ}$.  Then~$f\succ 1$, so $f\circ g\succ 1$, hence if $f\circ g\in\R f+\T_{\preceq}$, then $(f\circ g)^\dagger=f^\dagger+(1+\varepsilon)^\dagger$ with $\varepsilon\prec 1$ in $\T$. Also
$(f\circ g)^\dagger=(f^\dagger\circ g)\cdot g'$. 
Thus if $g'\not\sim 1$, then for $f:=\ex^x$ we have $f\circ g\notin \R f + \T_{\preceq}$,  and  if $g'\sim 1$, then 
$f:=\ex_n$ for high enough $n$ gives $f^\dagger\circ g\not\sim f^\dagger\succ 1$ by Lemma~\ref{lemest} (since $g\ne x$), and thus
$f\circ g\notin \R f + \T_{\preceq}$. 
\end{proof}

\noindent
This yields a left action $(g,\sigma_\alpha)\mapsto\sigma_{g\alpha}$ of the group~$\T^{>\R}$
on the set $\mathcal E$, under which~$\mathcal E^\k$ (for any subfield $\k$ of $\R$) is invariant.  
The image of $(g,\sigma)\in \T^{>\R}\times \mathcal E$ in~$\mathcal E$ under this action will be denoted by
$g\sigma$. 
For $\alpha\in D$, $f,g\in\T$, $g>\R$, we have 
$$\sigma_{g\alpha}(\ex^f)\ =\ \ex^{\alpha(f\circ g)+\sigma_{g\alpha}(f)}.$$

\begin{lemma}\label{lem:T under comp}
Let $g\in\T^{>\R}$, $\alpha\in D$. Then $\sigma_\alpha(g)\in\T^{>\R}$, and for  $f\in\T$ we have
$$\sigma_\alpha(f\circ g)\ =\ \sigma_{g\alpha}(f)\circ\sigma_\alpha(g).$$
\end{lemma}
\begin{proof}
Clearly $\sigma_\alpha(g)\in\T^{>\R}$.
By induction on $m$ we show that the desired equality holds for $f\in E_{m,n}$.
For $m=0$ it is enough to show this for $f=\ell_n^r$ ($r\in\R$) since~$E_{0,n}=\R[[\ell_n^\R]]$ and the relevant maps are strongly linear.
Now by \eqref{eq:sigma log}, 
\begin{align*}
\sigma_\alpha(\ell_n^r\circ g)\ =\  \sigma_\alpha\big((\ell_n\circ g)^r\big) &=\ 
\ex^{\alpha(r\log(\ell_n\circ g))-r\alpha(\log(\ell_n\circ g))}\big(\sigma_\alpha(\ell_n\circ g)\big)^r \\
&=\ \ex^{\alpha(r(\ell_{n+1}\circ g))-r\alpha(\ell_{n+1}\circ g)}\big(\sigma_\alpha(\ell_n\circ g)\big)^r \\
&=\ \ex^{g\alpha(r\ell_{n+1})-rg\alpha(\ell_{n+1})}\big(\sigma_\alpha(\ell_n\circ g)\big)^r
\end{align*}
and
\begin{align*}
\sigma_{g\alpha}(\ell_n^r)\circ\sigma_\alpha(g)\ &=\ 
\ex^{g\alpha(r\log\ell_n)-rg\alpha(\log\ell_n)}\big(\sigma_{g\alpha}(\ell_n)\big)^r\circ\sigma_\alpha(g) \\
&=\ \ex^{g\alpha(r\ell_{n+1})-rg\alpha(\ell_{n+1})}\big(\sigma_{g\alpha}(\ell_n)\circ\sigma_\alpha(g) \big)^r,
\end{align*}
so
we reduce to the case $f=\ell_n$, which we prove by induction on $n$. For $n=0$ we have $\sigma_{g\alpha}(\ell_0) = \ell_0 = x$ and so
$\sigma_\alpha(\ell_0\circ g) = \sigma_\alpha(g) = \sigma_{g\alpha}(\ell_0) \circ \sigma_\alpha(g)$.
Suppose the desired identity holds for a certain $n$. Then by \eqref{eq:sigma log},
\begin{align*}
\sigma_\alpha(\ell_{n+1}\circ g)\  &=\  \sigma_\alpha\big(\!\log(\ell_n\circ g)\big) \\
						&=\ \log \sigma_\alpha(\ell_n\circ g) -   \alpha\big(\log(\ell_n\circ g)\big) \\
					&=\ \log\big(\sigma_{g\alpha}(\ell_n)\circ\sigma_\alpha(g)\big) -  \alpha(\ell_{n+1}\circ g) \\
					&=\  \big(\!\log\sigma_{g\alpha}(\ell_n)\big) \circ \sigma_\alpha(g)  -  g\alpha(\ell_{n+1}) \\
			&=\  \big(\sigma_{g\alpha}(\ell_{n+1}) +   g\alpha(\ell_{n+1})\big) \circ \sigma_\alpha(g)  -  g\alpha(\ell_{n+1}) \\
								&=\  \sigma_{g\alpha}(\ell_{n+1})\circ\sigma_\alpha(g).
\end{align*}
This proves the case $m=0$. Suppose $m$ is such that the desired equality holds for all $f\in E_{m,n}$. 
Let $f\in E_{m+1,n}$, and let $a$ range over $(E_{m,n})_\succ$. Then $f=\sum_{a} f_a\ex^a$ where $f_a\in\R[[\ell_n^\R]]=E_{0,n}$, and so
\begin{align*}
\sigma_\alpha(f\circ g)\	&=\  \sum_a \sigma_\alpha(f_a\circ g) \sigma_\alpha(\ex^{a\circ g}) \\
			&=\  \sum_a \big(\sigma_{g\alpha}(f_a)\circ\sigma_\alpha(g)\big)\, \ex^{\alpha(a\circ g)+\sigma_\alpha(a\circ g)} \\
	&=\  \sum_a \big(\sigma_{g\alpha}(f_a)\circ\sigma_\alpha(g)\big)\, \ex^{g\alpha(a)+\sigma_{g\alpha}(a)\circ\sigma_\alpha(g)} \\
	&=\  \left(\sum_a \sigma_{g\alpha}(f_a) \,\ex^{g\alpha(a)+\sigma_{g\alpha}(a)}\right) \circ \sigma_\alpha(g) = \sigma_{g\alpha}(f) \circ \sigma_\alpha(g),
\end{align*}
where we used the case $m=0$ for the second equality, and the inductive hypothesis for the third equality.
\end{proof}

\begin{cor}\label{cor:T under comp, 1}
Let $g\in \T^{>\R}$ and set $\mathcal{E}|g:=\big\{\sigma\in \mathcal{E}:\, \sigma(g)=g\big\}$. Then the map 
$$\sigma\mapsto g\sigma\ \colon\  \mathcal{E}|g\to\mathcal E$$ is a group embedding.
\end{cor}
\begin{proof}
Let $\alpha,\beta\in D$ with $\sigma_\beta(g)=g$ and $f\in\T$. Then
$\sigma_\beta(f\circ g) =  \sigma_{g\beta}(f)\circ g$ by the previous lemma, hence
\begin{align*}
g(\alpha\circ\sigma_\beta+\beta)(f) &=\  
\alpha\big(\sigma_\beta(f\circ g)\big)+\beta(f\circ g)\\ &=\ \alpha\big(\sigma_{g\beta}(f)\circ g\big)+\beta(f\circ g)\ 
=\ (g\alpha\circ\sigma_{g\beta})(f)+g\beta(f)\\
&=\ (g\alpha\circ \sigma_{g\beta} + g\beta)(f).
\end{align*}
Thus, using \eqref{eq:cocycle}, we have $g(\sigma_{\alpha}\sigma_\beta)=\sigma_{g\alpha}\sigma_{g\beta}$. 
\end{proof}

\noindent
Taking $g:=\ell_n$ gives the automorphism $\alpha\mapsto\alpha{\downarrow_n}:=\ell_n\alpha$ of~$\Hom(\T,\R)$, with 
inverse~$\beta\mapsto\beta{\uparrow^n}:=\ex_n\beta$.
So $\sigma\mapsto \sigma{\downarrow_n}:=\ell_n\sigma$
is a permutation of~$\mathcal E$, with inverse~$\tau\mapsto \tau{\uparrow^n}:=\ex_n\tau$. For any subfield $\k$ of $\R$
we have $\mathcal E^\k{\downarrow_n}=\mathcal E^\k=\mathcal E^\k{\uparrow^n}$.
Also, $\ker (\alpha{\downarrow_n}) = (\ker\alpha){\uparrow^n}$ and
$\ker (\alpha{\uparrow^n}) = (\ker\alpha){\downarrow_n}$.
By Lemma~\ref{lem:T under comp} we have:

\begin{cor}\label{cor:T under comp, 2}
Let $\sigma\in \mathcal E$ and $f\in\T$. Then
$$\sigma(f{\uparrow^n})\  =\  \sigma{\uparrow^n}(f)\circ \sigma(\ex_n),\quad
\sigma(f{\downarrow_n})\  =\ \sigma{\downarrow_n}(f)\circ \sigma(\ell_n).
$$
\end{cor}

\subsection*{Translations are very healthy}
By  Lemma~\ref{lemEdown}  each $\sigma\in\SigmaAut_\der(\T)$ is healthy.
For use in the next section we now strengthen this for $\sigma\in\mathcal T$ (see also \cite{Sch}): 

\begin{prop}\label{prop:very healthy}
Each $\sigma\in\mathcal T$ is very healthy. 
\end{prop}

\noindent
The proof rests on a Taylor formula from \cite[(4)$_m$ in 6.8]{DMM2}, where $K_m$ is our $E_m$ and $L_{M,n}$ is our $E_{n,M}$:

\begin{lemma}\label{lem:Taylor}
Let $f\in E_m$ and $\varepsilon\in \T$, and suppose
$m=0$ and $\varepsilon\prec x$, or
$m\ge 1$ and $\abs{\varepsilon}< E_{m-1}^>$.
Then
the family $\big( \frac{f^{(i)}}{i!} \varepsilon^i \big)_{i\in\N}$ is summable in $\T$
and
$$f(x+\varepsilon) = \sum_{i=0}^\infty \frac{f^{(i)}}{i!} \varepsilon^i,$$
so if in addition $\varepsilon\in E_{m}$, then
$f(x+\varepsilon)\in E_{m}$.  
\end{lemma}

\noindent
Proposition~\ref{prop:very healthy} is an immediate consequence of the next lemma:

\begin{lemma} \label{lemvhe}
Let $c\in\R$ and $f\in E_{m,n}$ where $m\ge n$. Then $f(x+c)\in E_{m,n}$.
\end{lemma}
\begin{proof} Let $\alpha=(\alpha_i)_{i\in \N}$ be the sequence given by $\alpha_0=c$ and $\alpha_i:=0$ for $i\ge 1$.
Let $G_i=G_i^\alpha\in\R[[t_0,\dots,t_{i-1}]]$ be the power series defined in Lemma~\ref{lem:ellnstar}.
By the remark following that lemma, 
$G_{n+1}\in t_0\R[[t_0,\dots,t_{n}]]$. With $\ell_n^*$ defined just before that lemma, induction on $n$ gives  $\ell_n^*=\ell_n(x+c)$, and then the lemma itself yields
$$\ell_n(x+c)=\ell_n+\varepsilon\qquad\text{ where $\varepsilon:=G_n(1/\ell_0,\dots,1/\ell_{n-1})$.}$$
We first treat the case $m=n$, where the subcase $n=0$ is taken care of by Lemma~\ref{lem:Taylor} for $m=0$. Let $n\ge 1$, $f\in E_{n,n}$.
Put $g:=f{\uparrow^n}\in E_n$ and~$h:=f(x+c){\uparrow^n}$.
Then $f(x+c)=g\circ\ell_n\circ(x+c)=g(\ell_n+\varepsilon)$ and
 so $h=g(x+\delta)$ where~$\delta:=G_n(1/{\ex_n},\dots,1/{\ex_1})\in E_n$ satisfies  
 $\delta\preceq 1/{\ex_n}$,
so $|\delta| < E_{n-1}^>$.
Then Lem\-ma~\ref{lem:Taylor} gives $h\in E_n$ and so $f(x+c)=h{\downarrow_n}\in E_{n,n}$.

Suppose the lemma holds for a certain $m\ge n\ge 1$. Let $\fm\in G_{m+1,n}$. From ${G_{m+1}=G_m\ex^{A_m}}$ we obtain $G_{m+1,n}=G_{m,n}\ex^{A_{m,n}}$, so $\fm=\fn\ex^a$ where~${\fn\in G_{m,n}}$ and $a\in A_{m,n}$;  then $\fn\circ(x+c), a\circ (x+c)\in E_{m,n}$ by the inductive hypothesis, and so $\fm\circ (x+c)=(\fn\circ(x+c))\ex^{a\circ(x+c)}\in E_{m+1,n}$. By strong linearity   this yields~$f(x+c)\in E_{m+1,n}$ for all $f\in E_{m+1,n}$.
\end{proof}

\section{The Center of $\SigmaAut_\der(\T)$}\label{sec:center}

\noindent
In this section we show that $\SigmaAut_\der(\T)$ has trivial center. (With Corollary~\ref{cor:tfree}, this completes the proof of Theorem~B from the Introduction.) We precede this with a decomposition of
$\T$ into $\mathcal T$-invariant subspaces, to be used in the subsequent study of a certain linear difference
equation on $\T$.

\subsection*{Decomposing $\T$ into $\mathcal T$-invariant subspaces}
We decompose $\T$ as an internal direct sum of the strong truncation closed
$\R$-linear subspaces $ \T_{\flatter}$, $\T_{\comp}$, $\T_{\steeper}$: 
$$\T\  =\  \T_{\flatter} \oplus \T_{\comp} \oplus \T_{\steeper}$$
where for $\square$ any of the symbols $\flatter$, $\comp$, $\steeper$:
$$\T_{\square}\ :=\ \{ f\in\T: \supp f \,\square\, \ex^x \}.$$
Thus $\T_{\flatter}=\T^\flat$, and $\T_{\comp}$, $\T_{\steeper}$ are $\T^\flat$-linear subspaces of $\T$.
Proposition~\ref{prop:5.11} and the display just before Lemma~\ref{lem:5.11} give $\mathcal T$-invariance of $\T^\flat$:

\begin{lemma}\label{flatdec1} For $f\in \T^\flat$ and $\sigma\in\mathcal T$ we have $\sigma(f)\in\T^\flat$, with $\sigma(f)\sim f$ if $f\ne 0$.
\end{lemma} 

\noindent
For $\fm\in G^{\operatorname{LE}}$ and $\fm\comp \ex^x$ we have $c\in\R^\times$ and~$\varepsilon\in\T_\succ$ with~$\log\fm = cx+\varepsilon$
where~$\varepsilon\prec x$, and then $\fm=\fn\ex^{cx}$ with~$\fn:=\exp\varepsilon\in G^\flat$. 
Now $$\T^\flat\ =\ \bigcup_{m,n}\R[[G^\flat_{m,n}]]\ \subseteq\ \R[[ G^\flat]], \qquad \ex^{\R^\times\cdot x}\ \subseteq\ G^\sharp,$$ and inside $\R[[ G^{\operatorname{LE}}]]=\R[[G^\flat]][[G^\sharp]]$ we have 
\begin{equation}\label{eq:Tcomp} \T_{\comp}\ =\  \bigcup_{m,n}\R[[G^\flat_{m,n}]][[\ex^{\R^\times\cdot x}]]. \end{equation} 
(The inclusion of   $\T_{\comp}$ in the right side requires a little argument: for $f\in \T_{\comp}$, take~$m$,~$n$ such that
$\ex^{\R x} \cup \supp f \subseteq G_{m,n}$; then for $\fm=\fn\ex^{cx}\in \supp f$ with $\fn\in G^\flat$ and $c\in \R^\times$ we have
$\fn\in G_{m,n}\cap G^\flat=G^\flat_{m,n}$, so $f\in \R[[G^\flat_{m,n}]][[\ex^{\R^\times\cdot x}]]$.) We can now prove the $\mathcal{T}$-invariance of $\T_{\comp}$:

\begin{lemma}\label{flatdec2} Let $\sigma\in \mathcal{T}$ and $g\in \T_{\comp}$. Then $\sigma(g)\in \T_{\comp}$ and $\sigma(g)\asymp g$.
\end{lemma}
\begin{proof} By Lemma~\ref{flatdec1} we have  $\sigma(\T^\flat)\subseteq \T^\flat$. By Lemma~\ref{lemvhe}, 
 $\sigma(E_{m,n})\subseteq E_{m,n}$ for~$m\ge n$. Also $E_{m,n}=\R[[G_{m,n}]]$, so
 $\T^\flat\cap E_{m,n}=\R[[G^\flat_{m,n}]]$, and hence for~${m\ge n}$ we have 
 $\sigma(\R[[G^\flat_{m,n}]])\subseteq \R[[G^\flat_{m,n}]]$.  We have $\sigma=\sigma_d$ with $d\in \R$. Then for~${c\in\R^\times}$ and $f\in \R[[G^\flat_{m,n}]]$ with $m\ge n$, we have $\sigma( f\ex^{cx} ) =  \ex^{cd}\sigma(f)\ex^{cx}\in \T^\flat \ex^{cx}$, with~$\sigma(f\ex^{cx})\sim \ex^{cd} (f\ex^{cx})$ if $f\ne 0$. Now use \eqref{eq:Tcomp} and the strong linearity of $\sigma$. 
\end{proof} 

\noindent
Next we decompose $\T_{\steeper}$ as an internal direct sum 
\begin{equation}\label{eq:Tsteeper}
\T_{\steeper}\  =\  \T_{\steeper,\prec} \oplus \T_{\steeper,\succ}
\end{equation}
of the following strong truncation closed $\R$-linear subspaces of $\T$: 
$$\T_{\steeper,\prec}\ :=\ \T_{\steeper}\cap\T_{\prec} \ =\  \{ f\in\T:  
\log\supp f \succ x, \ \log \supp f<0
\} $$ 
and
$$\T_{\steeper,\succ}\ :=\ \T_{\steeper}\cap\T_{\succ}\  =\  \{ f\in\T:  
\log\supp f \succ x, \ \log \supp f>0
\}.$$ 
Inside $\R[[G^{\operatorname{LE}}]]$ we have $\T_{\steeper,\prec}  =  \T\cap \R[[G_{\steeper,\prec}]]$ and $\T_{\steeper,\succ}=  \T\cap \R[[G_{\steeper,\succ}]]$  where
\begin{align*} G_{\steeper,\prec}\ &:=\ \{\fm\in G^{\operatorname{LE}}:\ \log \fm \succ x,\ \log \fm <0\},\\
G_{\steeper,\succ}\ &:=\  \{\fm\in G^{\operatorname{LE}}:\ \log \fm \succ x,\ \log \fm >0\}.
\end{align*} 
Next we focus on $\T_{\steeper,\prec}$. Let $\sigma\in \mathcal{T}$. For  $\fm\in G^{\operatorname{LE}}$,
and $b:=\sigma(\log\fm)$ we have
\begin{multline*}
\sigma(\fm)\  =\  \exp b\  =\  c\,\fn\, h \\ \text{where $\fn:=\exp(b_\succ)\in G^{\operatorname{LE}}$, 
  $c:=\exp(b_1)\in\R^\times$,  $h:=\exp(b_\prec)\in 1+\T_\prec$.}
  \end{multline*} So 
$\supp \sigma(\fm)\subseteq \fn\supp h$.
If also $\fm\in  \T_{\steeper}$, that is, $x\prec \log\fm$, then
  $x\sim\sigma(x)\prec b$ and hence $b\sim b_\succ$. Here is the $\mathcal{T}$-invariance of $\T_{\steeper,\prec}$:

\begin{lemma}\label{lem:steeper, prec}
Let $f\in (\T_{\steeper,\prec})^{\ne}$ and $\sigma=\sigma_c\in\mathcal T$, $c\in\R$. 
Then $\sigma(f)\in \T_{\steeper,\prec}$, with~$\sigma(f)\prec f$ if $c>0$
and $\sigma(f) \succ f$ if   $c<0$.
\end{lemma}
\begin{proof}
Let $\fm\in G^{\operatorname{LE}}\cap \T_{\steeper}$, $b:=\sigma(\log\fm)$. 
  Consider an element $\fn\fv$ of $\supp\sigma(\fm)$, where $\fn=\ex^{b_\succ}$, $\fv\in\supp\ex^{b_\prec}$.
Then~$\fv\preceq 1$ and therefore $\log \fv\le 0$, hence~$\log(\fn\fv)=b_\succ+\log\fv\le b_\succ$.
Now if in addition~${\fm\prec 1}$, then~$\log \fm <0$, hence~$b_\succ<0$,  
  so~$\log(\fn\fv)<0$ and $\log(\fn\fv)\succeq b_\succ\succ x$, and
thus~$\sigma(\fm)\in \T_{\steeper,\prec}$.
By strong linearity of $\sigma$ this yields~$\sigma(f)\in \T_{\steeper,\prec}$.
For the rest use Corollary~\ref{cor:steep fm}.
\end{proof}

\noindent
The case of $\T_{\steeper,\succ}$ requires more effort. We begin with two observations:

\begin{lemma}\label{lem:exp bound}
Let $b\in \T_\prec$ and $\fw\in G^{\operatorname{LE}}$ be such that $\log\fm\preceq\fw$ for all $\fm\in\supp b$. Then $\log\fv\preceq\fw$ for all
$\fv\in\supp\ex^b$.
\end{lemma}
\begin{proof}
By $b\prec 1$ we have $\ex^b=\sum_n \frac{b^n}{n!}$.
Hence, if $\fv\in\supp\ex^b$, then we can take $n$ and 
$\fm_1,\dots,\fm_n\in\supp b$ with $\fv=\fm_1\cdots\fm_n$, where
$\log\fm_1,\dots,\log\fm_n\preceq\fw$, and thus also~$\log \fv=\log \fm_1+\cdots+\log\fm_n\preceq \fw$.
\end{proof}

\begin{lemma}\label{lem:logbd}
Let $f\in\T_\succ$, $\sigma\in\mathcal T$, and $g=\sigma(f)$. Then for all $\mathfrak u\in\supp(g_\prec)$ we have $\log \mathfrak u\preceq\log x$.
\end{lemma}
{\sloppy\begin{proof}
For $r\in\R$ we have $\sigma(\ell_n^r)\in \R[[\ell_0^{\R}\cdots\ell_n^{\R}]]$, by Lemma~\ref{lem:sigma is strong}, and
if~${r_0,\dots,r_n\in\R}$, then
$\log ( \ell_0^{r_0}\cdots\ell_n^{r_n} ) =r_0\ell_1+\cdots+r_n\ell_{n+1}\preceq \ell_1= \log x$. Thus~$\log \mathfrak u\preceq\log x$ for all~${\fm\in \ell_n^{\R}}$ and $\fu\in\supp \sigma(\fm)$.
Using strong linearity of $\sigma$, this yields the base case~$m=0$ of an induction on $m$ to prove the conclusion of the lemma for~${f\in (E_{m,n})_\succ}$. 
Now suppose  the claim holds  for a certain $m$.
To obtain it for~$m+1$ in place of $m$, it suffices by strong linearity of~$\sigma$ to treat the case~$f\in G_{m+1,n}^{\succ 1}$. 
We have~$G_{m+1,n}=
\ell_n^{\R}\exp\big((E_{m,n})_\succ\big)$, so~$f=\ell_n^r\ex^a$ where~$r\in\R$, $a\in (E_{m,n})_\succ$.
Then $g=\sigma(\ell_n^r)\sigma(\ex^a)$. Let $\fu\in \supp(g_{\prec})$. Then $\fu=\fn\fv$ with~$\fn\in \supp \sigma(\ell_n^r)$ and $\fv\in \supp \sigma(\ex^a)$. Now $\fu\prec 1$, so $\fv\prec \fn^{-1}$. If $\fv\succeq 1$, then
$0\le \log \fv < \log \fn^{-1}\preceq  \log x$ in view of $\fn, \fn^{-1}\in \ell_0^\R\cdots \ell_n^\R$ and the first part of the proof, and so
$\log \fu=\log \fn + \log \fv\preceq \log x$. So it remains to handle the case $\fv\prec 1$, where it again suffices to show
$\log \fv\preceq \log x$, since we already have $\log \fn\preceq \log x$ by the first part of the proof.  Now $f\succ 1$ and $\ell_n^{\R}$ is convex in~$G_{m+1,n}$, so $\ex^a\succ 1$.  So we have a reduction to the case $f=\ex^a$, $g=\sigma(\ex^a)=\ex^b$ with $b=\sigma(a)$. Let $\fu\in \supp( g_\prec)$; it suffices to show that then $\log \fu\preceq \log x$.

 The inductive hypothesis gives $\log \mathfrak \fm\preceq\log x$ for all~$\fm\in\supp(b_\prec)$
 and thus~$\log \mathfrak v\preceq\log x$ for all~$\fv\in\supp \ex^{b_\prec}$, by Lemma~\ref{lem:exp bound}. 
 Now $\fu=\fn\,\fv\prec 1$ where~$\fn=\ex^{b_\succ}$ and~$\fv\in\supp(\ex^{b_\prec})$.
 Also $a>0$, so $b>0$, hence $\R< \fn < \fv^{-1}$,  so~$\R<\log\fn < \log \fv^{-1}=-\log\fv$,  so $\log\fn\preceq\log\fv\preceq\log x$, and thus $\log \fu \preceq\log x$. 
\end{proof}}

\noindent
We can now prove the $\mathcal{T}$-invariance of $\T_{\steeper,\succ}$:

\begin{lemma}\label{lem:steeper, succ}
Let $f\in (\T_{\steeper,\succ})^{\ne}$ and $\sigma=\sigma_c\in\mathcal T$, $c\in\R$. Then $\sigma(f)\in \T_{\steeper,\succ}$, with~$\sigma(f)\succ f$ if $c>0$ and $\sigma(f) \prec f$ if $c<0$.
\end{lemma}
\begin{proof}
Let $\fm\in G^{\operatorname{LE}}\cap \T_{\steeper,\succ}$; then $\log\fm\succ x$ and $\log\fm>0$.
Put~$b:=\sigma(\log\fm)$, so~$b \succ x$, $b>0$, and thus~$b\sim b_\succ$.  Also $\sigma(\fm)=\ex^b$.  
Let $\fu\in \supp \sigma(\fm)$. Then $\fu=\fn\fv$ where~$\fn=\ex^{b_\succ}$, $\fv\in\supp \ex^{b_\prec}$.  Now $\log \fm\in \T_{\succ}$, so applying Lemma~\ref{lem:logbd} to $\log \fm$ in the role of $f$ we obtain $\log \mathfrak{a} \preceq \log x$ for all
$\mathfrak{a}\in \supp b_{\prec}$. Then Lemma~\ref{lem:exp bound} applied to $b_{\prec}$ in the role of $b$ yields
$\log \fv\preceq \log x\prec x$. Hence~$\log (\fn\fv)=b_\succ+\log\fv\sim b_\succ$, so~${\log \fu\succ x}$ and~$\log \fu>0$.
Therefore~$\sigma(\fm)\in \T_{\steeper,\succ}$.  By strong linearity of $\sigma$ this yields~$\sigma(f)\in \T_{\steeper,\succ}$.
For the rest use Corollary~\ref{cor:steep fm}.
\end{proof}

\subsection*{Solving a difference equation}		
Consider the  strongly linear map 
$$\Delta\ \colon\ \T\to\T,\qquad \Delta(f)\ :=\ f(x+1)-f(x).$$
Elaborating on \cite[Section~5]{vdHmeta}, we shall prove:

\begin{prop}\label{prop:Delta} $\Delta$ is surjective. 
\end{prop}

\noindent
Towards this we exhibit as in the previous subsection ~$\T$ as an internal direct sum of its
$\R$-linear subspaces $\T_{\flatter}, \T_{\comp}$, and $\T_{\steeper}$, 
$$\T\  =\  \T_{\flatter} \oplus \T_{\comp} \oplus \T_{\steeper},$$
each of which is $\Delta$-invariant.
Proposition~\ref{prop:Delta flat} yields a strongly linear right inverse  of~$\Delta|_{\T_{\flatter}}=\Delta^\flat\colon \T^\flat\to \T^\flat$.  In the  lemmas below we treat $\Delta|_{\T_{\comp}}$ and $\Delta|_{\T_{\steeper}}$ likewise. 

\begin{lemma}\label{lemlricomp} 
$\Delta|_{\T_{\comp}}\colon \T_{\comp}\to \T_{\comp}$ has a linear right inverse $\Sigma_{\comp}$.
\end{lemma}
\begin{proof}
Recall how \eqref{eq:Tcomp} exhibits $\T_{\comp}$ as a union. Let $m\ge n$, $f\in \R[[G^\flat_{m,n}]]$, and~${c\in \R^\times}$.
Then by the proof of Lemma~\ref{flatdec2}, 
$$\Delta (f\ex^{cx})\  =\  \big(\!\ex^cf(x+1) -f(x)\big)\cdot\ex^{cx}\  =\ \Delta_c(f)\ex^{cx}\in \R[[G^\flat_{m,n}]]\ex^{cx}$$
where $\Delta_c:=\big(\!\ex^c \exp(t)-1\big)(\der)$, and $\der$ denotes the derivation of~$\T^\flat$.  With~$S_{\ex^c}\in\R[[t]]$   as in Remark~\ref{rem:S},   
 $\Sigma_c:=S_{\ex^c}(\der)$ is a strongly linear right inverse of $\Delta_c\colon \T^\flat\to \T^\flat$, and $\Sigma_c$ maps
 $\R[[G^\flat_{m,n}]]$ into $\R[[G^\flat_{m,n}]]$ for $m\ge n$.   

Let now $f\in \T_{\comp}$. Then for some $m\ge n$ we have
$f=\sum_{c\in \R^\times} f_c \ex^{cx}$ with all~$f_c\in \R[[\G^\flat_{m,n}]]$, and then $\Delta(f)=\sum_{c\in \R^\times}\Delta_c(f_c)\ex^{cx}$. For $m\ge n$ and such $f$ the family $\big(\Sigma_c(f_c)\ex^{cx}\big){}_{c\in \R^\times}$ in $\R[[G^\flat_{m,n}\ex^{\R x}]]$ is summable, and this yields a linear map
$$f=\sum_{c\in \R^\times}  f_c\ex^{cx}\ \mapsto\   \sum_{c\in \R^\times} \Sigma_c(f_c)\ex^{cx}\ :\  \R[[G^\flat_{m,n}]][[\ex^{\R^\times x}]]\to\R [[G^\flat_{m,n}]][[\ex^{\R^\times x}]].$$
These maps for $m\ge n$ have a common extension to a linear map
$\Sigma_{\comp} \colon  \T_{\comp}\to\T_{\comp}$,
which is a right inverse of $\Delta|_{\T_{\comp}}$.
\end{proof}

\begin{lemma} The restriction of $\Delta$ to a map from $\T_{\steeper}$ to itself is surjective. 
\end{lemma}

\begin{proof}
By \eqref{eq:Tsteeper}, $\T_{\steeper}$  
is the internal direct sum of its $\R$-linear subspaces $\T_{\steeper,\prec}$ and~$\T_{\steeper,\succ}$, each of which is  
$\mathcal T$-invariant 
by Lemmas~\ref{lem:steeper, prec} and~\ref{lem:steeper, succ}.
 Set $\sigma=\sigma_1$. Then~$\sigma(f)\prec f$ for~$f\in\T^{\ne}_{\steeper,\prec}$,
and $\sigma^{-1}(f)\prec f$ for~${f\in\T^{\ne}_{\steeper,\succ}}$, and
by Proposition~\ref{prop:very healthy}, $\sigma$ and $\sigma^{-1}$ map $E_{m,n}$ into itself for $m\ge n$. 
Now  $$ \T_{\steeper, \prec}\ =\ \bigcup_{m,n} \big(E_{m,n}\cap \T_{\steeper,\prec}\big),\quad E_{m,n}\cap \T_{\steeper,\prec}\ =\  \R[[G_{m,n}\cap G_{\steeper, \prec}]],$$
and Corollary~\ref{cor:cnPhin} gives for $m\ge n$ a strongly linear map 
$$f\mapsto \sum_{i=0}^\infty \sigma^i(f)\ :\  \R[[G_{m,n}\cap G_{\steeper, \prec}]]\to \R[[G_{m,n}\cap G_{\steeper, \prec}]]$$ such that
 $\Delta\big({-\sum_i \sigma^i(f)}\big)=f$ for each $f\in  E_{m,n}\cap \T_{\steeper,\prec}$. Thus the restriction of~$\Delta$ to a map from $\T_{\steeper,\prec}$ to itself is surjective. 
Likewise we have for~${m\ge n}$ a strongly linear map~$f\mapsto \sum_{i=0}^\infty \sigma^{-(i+1)}(f)$ from
$\R[[G_{m,n}\cap G_{\steeper,\succ}]]$ to itself such that~$\Delta\left(\sum_{i=0}^\infty \sigma^{-(i+1)}(f)\right)=f$, so the restriction of $\Delta$ to a map from $\T_{\steeper,\succ}$ to itself is also surjective. 
\end{proof}

\noindent
This concludes the proof of Proposition~\ref{prop:Delta}. \qed

\medskip\noindent
{\it Remark}. Proposition~\ref{prop:Delta} suffices for use in the next subsection, but it is worth stating a strengthening: {\em For each $g\in \T$ there is a unique $f\in \T$ with~${\Delta(f)=g}$ and~$f_1=0$. The map $\Sigma \colon \T\to \T$ that assigns to each $g\in \T$ this  $f$ is strongly linear}. The first of these claims holds because $\ker \Delta=\R$ by Lemma~\ref{lem:fixed field of M}. For the strong linearity claim we note first that the map $\Sigma$ is obviously linear. Moreover, it agrees on $\T^\flat$ with $h\mapsto \Sigma^\flat(h)-\Sigma^{\flat}(h)_1$,
where $\Sigma^\flat$ is as in the proof of Proposition~\ref{prop:Delta flat}.  On~$\T_{\comp}$ it agrees with $\Sigma_{\comp}$ from the proof of Lemma~\ref{lemlricomp}, on $\T_{\steeper, \prec}$ with $-\sum_i \sigma_1^i$, and on $\T_{\steeper, \succ}$ with $\sum_i \sigma_1^{-(i+1)}$.  
A close look at these partial right inverses to $\Delta$ then yields the strong linearity of $\Sigma$.

\subsection*{The center of $\SigmaAut_\der(\T)$}
Set $\mathcal G:=\SigmaAut_\der(\T)=\mathcal E\ltimes\mathcal T$. 
For  $\alpha\in D$, $c\in\R$,  
$$\sigma_{\alpha_c}\ =\ \sigma_c^{-1}\sigma_\alpha\sigma_c\quad\text{ where $\alpha_c=\alpha\circ\sigma_c\in D$.}$$
We now aim to determine the center $Z(\mathcal G)$ of $\mathcal G$. Let
\begin{align*}
C_{\mathcal E}(\mathcal T)\	:&=\ \{\sigma\in\mathcal E:\ \text{$\sigma\tau=\tau\sigma$ for all $\tau\in\mathcal T$}\}\\
					&=\ \{ \sigma_\alpha:\ \text{$\alpha\in D$, $\alpha_c=\alpha$ for all $c\in\R$}\}
\end{align*}
be the centralizer of $\mathcal T$ in $\mathcal E$. 

\begin{lemma}
$Z(\mathcal G)\subseteq C_{\mathcal E}(\mathcal T)\mathcal{T}$.
\end{lemma}
\begin{proof}
Let $\sigma\in\mathcal E$, $\tau\in\mathcal T$.
For $\sigma\tau$ to lie in $Z(\mathcal G)$, it must commute with every~$\rho\in\mathcal G$:
$(\sigma\tau)\rho=\rho(\sigma\tau)$. 
Since $\mathcal T$ is abelian, this simplifies for $\rho\in\mathcal T$ to $\sigma\rho=\rho\sigma$, that is,
$\sigma\in C_{\mathcal E}(\mathcal T)$.
\end{proof}

\noindent
The next lemma gives $C_{\mathcal E}(\mathcal T)=\{\id_{\T}\}$:

\begin{lemma}\label{lemcomsigmaone}
Suppose $\sigma\in\mathcal E$ commutes with $\sigma_1\in\mathcal T$. Then $\sigma=\id_{\T}$.
\end{lemma}
\begin{proof}
Set $\alpha:=\alpha_\sigma$. Then $\sigma_{\alpha}=\sigma$, so $\sigma_{\alpha_1}=\sigma_1^{-1}\sigma\sigma_1=\sigma=\sigma_{\alpha}$, hence
$$ \alpha\ =\ \alpha_1\ :=\ \alpha\circ \sigma_1.$$ 
Therefore, $f(x+1)-f(x)\in\ker\alpha$
for all $f\in\T$, so $\ker\alpha=\T$ by Proposition~\ref{prop:Delta}, and thus $\sigma=\id_{\T}$.
\end{proof}

\begin{lemma}
$Z(\mathcal G)\cap\mathcal T=\{\id_{\T}\}$.
\end{lemma}
\begin{proof}
Let $c\in\R$ and $\sigma_c\in Z(\mathcal G)$. Then for all $\alpha\in D$   we have
$\sigma_{\alpha_c}=\sigma_c^{-1}\sigma_\alpha\sigma_c=\sigma_\alpha$ and so $\alpha_c=\alpha$, that is, $f(x+c)-f(x)\in\ker\alpha$
for all $f\in\T$. Taking $f=\ex^x$ we obtain, for each $\alpha\in D$: $(\ex^c-1)\ex^x\in\ker\alpha$. This forces $c=0$: otherwise, we can choose $\alpha\in D$
such that $(\ex^c-1)\ex^x\notin\ker\alpha$.
\end{proof}

\noindent
Combining the last three lemmas now yields:

\begin{cor}
$Z(\mathcal G)=\{\id_{\T}\}$.
\end{cor}

\noindent
Let us also justify here a remark from the Introduction:

\begin{cor} $\mathcal{T}$ is not a normal subgroup of $\SigmaAut_{\der}(\T)$.
\end{cor}
\begin{proof} Let $\sigma\in \mathcal{E}$, $\sigma\ne \id_{\T}$. We have $(\sigma \sigma_1\sigma^{-1})(x)=x+1$, so
if $\sigma\sigma_1\sigma^{-1}\in \mathcal{T}$, then~$\sigma$ commutes with $\sigma_1$, contradicting Lemma~\ref{lemcomsigmaone}. Thus $\sigma \sigma_1\sigma^{-1}\notin \mathcal{T}$.
\end{proof}

\section{Some Filtrations of $\mathcal E$}\label{sec:filtrations}

\noindent
Here we consider various subgroups of $\mathcal E$ and determine their fixed fields,
and use this  to prove Corollaries~C and~D from the Introduction.
Throughout this section~$\k$ is a subfield of $\R$.

\subsection*{The groups $\mathcal E_n$ and $\mathcal L_n$}
We define the  subgroups
$$\mathcal E_n\ :=\ \{\sigma\in \mathcal{E}:\ \sigma(\ex_n)=\ex_n\}, \qquad \mathcal L_n\ :=\ \{\sigma\in \mathcal{E}:\ \sigma(\ell_n)=\ell_n\}.$$
of $\mathcal{E}$. Thus $\mathcal E_0 =\mathcal L_0 =\mathcal E$. We also have $\mathcal{L}_n=\SigmaAut_\der (\T|\ell_n)$ by Lemma~\ref{lem:sigmaelln=elln}. 
From Corollaries~\ref{cor:T under comp, 1} and \ref{cor:T under comp, 2} we obtain:

\begin{cor}\label{cor:TR and updownshift}
The maps
$$\sigma\mapsto \sigma{\uparrow^n}\ \colon\ \mathcal E_n\to\mathcal E, \qquad  \sigma\mapsto \sigma{\downarrow_n}\ \colon\ \mathcal L_n\to\mathcal E $$
are group embeddings. Moreover, for $f\in\T$,
$$\sigma{\uparrow^n}(f)\ =\ \sigma(f{\uparrow^n}){\downarrow_n}\  \ \text{\ for $\sigma\in\mathcal E_n$}, \qquad 
\sigma{\downarrow_n}(f)\ =\ \sigma(f{\downarrow_n}){\uparrow^n}\ \text{\ for $\sigma\in\mathcal L_n$.}$$
\end{cor}

\noindent
It follows that we have mutually inverse group isomorphisms 
$$\sigma\mapsto \sigma{\uparrow^n}\ :\ \mathcal E_n\to \mathcal L_n, \qquad \tau\mapsto \tau{\downarrow_n}\ :\ \mathcal L_n\to \mathcal E_n.$$
Combining Corollaries~\ref{cor:id on elln} and~\ref{cor:TR and updownshift} yields:

\begin{cor}\label{cor:elln and en}
Let $\sigma\in\mathcal E$, $\alpha=\alpha_\sigma$. Then the following are equivalent:
\begin{enumerate}
\item[\rm(i)] $\sigma\in \mathcal{L}_n$;
\item[\rm(ii)] $\sigma(\ell_i)=\ell_i$  for  $i=0,\dots,n$;
\item[\rm(iii)] $\alpha(\ell_i)=0$ for $i=1,\dots,n$;
\item[\rm(i$'$)] $\sigma{\downarrow_n}(\ex_n)=\ex_n$;
\item[\rm(ii$'$)] $\sigma{\downarrow_n}(\ex_i)=\ex_i$ for $i=0,\dots,n$;
\item[\rm(iii$'$)] $\alpha{\downarrow_n}(\ex_i)= 0$  for $i=0,\dots,n-1$.
\end{enumerate}
\end{cor}
 
\noindent
We put $\mathcal E_n^\k:=\mathcal E_n\cap \SigmaAut_\der^\k(\T)$ and $\mathcal L_n^\k:=\mathcal L_n\cap \SigmaAut_\der^\k(\T)$. So $\mathcal E_n^\Q=\mathcal E_n$, and if $\k$ is a subfield of the subfield $\l$ of $\R$, then
$\mathcal E_n^\k\supseteq\mathcal E_n^\l$, and likewise with $\mathcal L$ in place of $\mathcal E$. We have
$\big(\mathcal E_n^{\k}\big){\uparrow^n}= \mathcal L_n^{\k}$. 
By Corollary~\ref{cor:elln and en}, 
$$\mathcal E_n^\k\	 =\  \SigmaAut_\der^\k\!\big(\T|\R[[\ex_n^\k\cdots\ex_0^\k]]\big)\  
=\ \big\{\sigma_\alpha: \alpha\in D^\k,\ \alpha(\ex_{n-1})=\cdots=\alpha(\ex_0)=0 \big\}$$ 
and
$$\mathcal L_n^\k\  =\   \SigmaAut_\der^\k\!\big(\T|\R[[\ell_0^\k\cdots\ell_n^\k]]\big)\   =\  
\big\{\sigma_\alpha: \alpha\in D^\k,\ \alpha(\ell_1)=\cdots=\alpha(\ell_n)=0 \big\}.$$
Note that $\R[[\ex_n^\k\cdots\ex_0^\k]]$ and $\R[[\ell_0^\k\cdots\ell_n^\k]]$ are differential subfields of $\T$, and that
$\mathcal L_n^\k$ is a normal subgroup of $\SigmaAut_\der^{\k}(\T)$ by
Lemma~\ref{lem:sigma is strong}.  Clearly,
$$
\mathcal E^\k = \mathcal E_0^\k \supseteq \cdots  \supseteq \mathcal E_n^\k \supseteq \mathcal E_{n+1}^\k \supseteq \cdots, \qquad
\mathcal E^\k = \mathcal L_0^\k \supseteq \cdots \supseteq \mathcal L_n^\k \supseteq \mathcal L_{n+1}^\k \supseteq \cdots.$$ 

\begin{lemma}\label{lem:fixed fields under updown shift}
If $\mathcal G\subseteq \mathcal E_n$, then $\T^{\mathcal G}{\downarrow_n}=\T^{\mathcal G{\uparrow^n}}$. If $\mathcal H\subseteq \mathcal L_n$, then $\T^{\mathcal H}{\uparrow^n}=\T^{\mathcal H{\downarrow_n}}$. \textup{(}Hence $\T^{\mathcal E_n^\k}{\downarrow_n}=\T^{\mathcal L_n^\k}$.\textup{)}
\end{lemma}

\begin{proof}
Let $f\in\T$ and $\mathcal G\subseteq\mathcal E_n$. Then by Corollary~\ref{cor:TR and updownshift}:
\begin{align*}
f\in\T^{\mathcal G}{\downarrow_n}	&\quad\Longleftrightarrow\quad f{\uparrow^n}\in\T^{\mathcal G} \\
										&\quad\Longleftrightarrow\quad \text{$\sigma(f{\uparrow^n})=f{\uparrow^n}$ for all $\sigma\in\mathcal G$}\\
										&\quad\Longleftrightarrow\quad \text{$\sigma{\uparrow^n}(f)=f$ for all $\sigma\in\mathcal G$} \quad\Longleftrightarrow \quad
f\in										\T^{\mathcal G{\uparrow^n}}.
\end{align*}
Thus $\T^{\mathcal G}{\downarrow_n}=\T^{\mathcal G{\uparrow^n}}$, and $\T^{\mathcal H}{\uparrow^n}=\T^{\mathcal H{\downarrow_n}}$ for $\mathcal H\subseteq\mathcal L_n$ follows by applying this to~$\mathcal G=\mathcal H{\downarrow}_n$. 
\end{proof}

\noindent
We also set
$$\mathcal E_{\infty}^\k\ :=\ \bigcap_n\mathcal E_n^\k\ =\ \big\{\sigma\in\mathcal E^{\k}:\  \text{$\sigma(\ex_n)=\ex_n$ for all $n$}\big\}\ \supseteq\  \SigmaAut_\der^{\k}(\T|E)$$
and
$$\mathcal L_{\infty}^\k\ :=\  \bigcap_n \mathcal L_n^\k\  =\  \big\{ \sigma\in\SigmaAut_\der^\k(\T):\  \text{$\sigma(\ell_n)=\ell_n$ for all $n$} \big\}\  =\  \SigmaAut_\der^\k(\T|\T_{\log}^\k).$$
Note that $\mathcal L_{\infty}^\k$ is a normal subgroup of $\SigmaAut_{\der}^{\k}(\T)$ (being an intersection of such).
We put $\mathcal E_{\infty}:=\mathcal E_{\infty}^\Q=\bigcap_n \mathcal E_n$ and $\mathcal L_{\infty}:=\mathcal L_{\infty}^\Q=\bigcap_n \mathcal L_n$. Recall from Corollary~\ref{cor:sigmaelln=elln}  that for $\sigma\in\SigmaAut_\der(\T)$: 
 $$\sigma\in \mathcal L_{\infty} \quad\Longleftrightarrow\quad \sigma(\T_\succ)=\T_\succ\quad \Longleftrightarrow\quad \sigma(\R^>\cdot G^{\operatorname{LE}})=\R^>\cdot G^{\operatorname{LE}}.$$

\subsection*{The groups $\mathcal D_n$} Next a condition on $\sigma\in \mathcal E$ to act as the identity on $E_{m,n}$:

\begin{lemma}\label{lem:Emn in Tsigma}
Let $\alpha\in D$ and $\sigma=\sigma_\alpha$. Then 
$$E_{0,n}\subseteq\T^\sigma \ \Longleftrightarrow\  \text{$\alpha(\ell_1)=\cdots=\alpha(\ell_n)=0$ and 
$\alpha(r\ell_{n+1})=r\alpha(\ell_{n+1})$ for all $r\in\R$,}$$
and
$$E_{m,n}\subseteq\T^\sigma \ \Longleftrightarrow\  \text{$E_{0,n}\subseteq\T^\sigma$
and $(E_{m-1,n})_\succ\subseteq \ker\alpha$.}$$
\end{lemma}
\begin{proof}
For $m=0$ we have $E_{0,n}=\R[[\ell_n^\R]]$ and so
by strong linearity of $\sigma$,
\begin{align*}
E_{0,n}\subseteq\T^\sigma	&\ \Longleftrightarrow\  \text{$\sigma(\ell_n^r)=\ell_n^r$ for all $r\in\R$} \\
							&\ \Longleftrightarrow\  \text{$\ex^{\alpha(r\ell_{n+1})-r\alpha(\ell_{n+1})}\big(\sigma(\ell_n)\big)^r=\ell_n^r$ for all $r\in\R$},\ \text{by}\  \eqref{eq:sigmanellnr} \\
							&\ \Longleftrightarrow\  \text{$\sigma(\ell_n)=\ell_n$ and
							 $\alpha(r\ell_{n+1})=r\alpha(\ell_{n+1})$
							for all $r\in\R$}\\
							&\ \Longleftrightarrow\  \text{$\alpha(\ell_1)=\cdots=\alpha(\ell_n)=0$ and $\alpha(r\ell_{n+1})=r\alpha(\ell_{n+1})$
							for all $r\in\R$,}
\end{align*}
using Corollary~\ref{cor:id on elln} for the last equivalence.
This proves the first equivalence of the lemma. We show the second equivalence by
induction on~$m$, the case $m=0$ being trivial.
By strong linearity of~$\sigma$   we have
\begin{align*}
E_{m+1,n}\subseteq\T^\sigma	&\ \Longleftrightarrow\  \text{$\sigma(\fm)=\fm$ for all $\fm\in G_{m+1,n}$} \\
							&\ \Longleftrightarrow\  \text{$E_{0,n}\subseteq\T^\sigma$  and $\alpha(a)+\sigma(a)=a$ for  all $a\in (E_{m,n})_\succ$.}						
\end{align*}
If $E_{m+1,n}\subseteq\T^\sigma$, then $(E_{m,n})_\succ\subseteq E_{m+1,n}\subseteq\T^\sigma$,
hence the above equivalences give~$E_{0,n}\subseteq\T^\sigma$ and $(E_{m,n})_\succ\subseteq\ker\alpha$.
Conversely, if this condition holds, then the inclusion $(E_{m-1,n})_\succ\subseteq (E_{m,n})_\succ$ and the inductive hypothesis
yield $(E_{m,n})_\succ\subseteq E_{m,n}\subseteq\T^\sigma$, and so $E_{m+1,n}\subseteq\T^\sigma$, again 
by the equivalences above.
\end{proof}

\noindent
Set
$$D_{m,n}\ :=\ \big\{\alpha\in D^\R:\,  \ker \alpha \supseteq (E_{m-1,n})_\succ+\Z\ell_1+\cdots+\Z\ell_n\big\}.$$ Then 
 $D_{m,n}\supseteq D_{m+1,n}$. Moreover, $\ell_1,\dots,\ell_n\in (E_{n-1,n})_\succ$, so if $m\geq n$, then
$$D_{m,n}\ =\ \big\{\alpha\in D^\R:\, \ker \alpha \supseteq  (E_{m-1,n})_\succ \big\}.$$
We also set $$\mathcal D_{m,n}\ :=\ \{\sigma_\alpha:\alpha\in D_{m,n}\},$$ so 
$\mathcal D_{m,n}\supseteq \mathcal D_{m+1,n}$.
By Lemma~\ref{lem:Emn in Tsigma} we have $\mathcal D_{m,n}=\mathcal{D}^\R|E_{m,n}$. 
Corollary~\ref{cor:id on elln}, (iii)~$\Rightarrow$~(ii), yields
$\R[[\ell_0^\R\cdots\ell_n^\R]]\subseteq \T^{\mathcal D_{m,n}}$.
Next, put
\begin{align*}D_n\ &:=\ \bigcap_m D_{m,n}\ =\ \big\{\alpha\in D^\R:\ \ker \alpha \supseteq (E{\downarrow_n})_\succ \big\},\\
\mathcal D_n\  &:=\ \bigcap_m\mathcal D_{m,n}\ =\ \{\sigma_\alpha:\ \alpha\in D_n\}\ =\ \mathcal{E}^{\R}|E_{\downarrow_n}\ =\  \SigmaAut_\der^\R(\T|E{\downarrow}_n).
\end{align*}
Thus the subgroup $\mathcal{D}_n$ of $\SigmaAut_\der(\T)$ is contained in $\mathcal E_{\infty}^\R\cap\mathcal L_n^\R$.
Note that   $\mathcal D_n\supseteq \mathcal D_{n+1}$, with $\bigcap_n \mathcal D_n=\{1\}$.  
By Corollary~\ref{cor:extend from Euparrown}, $\mathcal D_n$ is a normal subgroup of $\SigmaAut_\der^{\R}(\T)$.
 Corollary~\ref{cor:TR and updownshift} gives the following:

\begin{cor}
Let $f\in\T$. Then
$$\sigma(f{\uparrow^n})\  =\  \sigma{\uparrow^n}(f){\uparrow^n}\quad\text{for $\sigma\in\mathcal D_0$}, \qquad
\sigma(f{\downarrow_n})\  =\  \sigma{\downarrow_n}(f){\downarrow_n}\quad\text{for $\sigma\in\mathcal D_n$.}$$
\end{cor}

\noindent
Since$(E{\downarrow})_\succ = (E_\succ){\downarrow}$ we have
for $\alpha\in D^\R$: 
$\ker(\alpha)\supseteq E_\succ\ \Leftrightarrow\ \ker(\alpha{\uparrow})\supseteq (E{\downarrow})_\succ$.
Thus~$D_0{\uparrow}= D_1$. In general, 
$D_m{\uparrow^n} = D_{m+n}$, and hence $D_{m+n}{\downarrow_n}=D_m$.
As a consequence, $\mathcal D_m{\uparrow^n} = \mathcal D_{m+n}$ and  $\mathcal D_{m+n}{\downarrow_n}=\mathcal D_m$.
Thus by Lemma~\ref{lem:fixed fields under updown shift}:

\begin{cor}\label{cor:fixed fields under updown shift}
If $\mathcal G\subseteq\mathcal D_0$, then $\T^{\mathcal G}{\downarrow_n}=\T^{\mathcal G{\uparrow^n}}$.  If $\mathcal H\subseteq\mathcal D_n$, then $\T^{\mathcal H}{\uparrow^n}=\T^{\mathcal H{\downarrow_n}}$.
Hence
$\T^{\mathcal D_{m}}{\downarrow_n}=\T^{\mathcal D_{m+n}}$ and
$\T^{\mathcal D_{m+n}}{\uparrow^n}=\T^{\mathcal D_{m}}$.
\end{cor}

\subsection*{Some definably closed differential subfields of $\T$}
We now show that the fixed field of $\mathcal D_n$ is $E{\downarrow_n}$: Corollary~\ref{cor:fixed field of T_n}. First some auxiliary results.

\begin{lemma}\label{lemneq}
Let $f\in E_{m+1,n}\setminus E_{m,n}$, $m\geq n$. Then $\sigma(f)\neq f$ for some $\sigma\in\mathcal D_{m,n}$.
\end{lemma}
\begin{proof}
Let $a$ range over $(E_{m,n})_\succ$. Then
$f=\sum_a f_a\ex^a$ with all $f_a\in\R[[\ell_n^\R]]=E_{0,n}$.
Since $f\notin E_{m,n}$, we have $a_0\in (E_{m,n})_\succ\setminus (E_{m-1,n})_\succ$ with $f_{a_0}\neq 0$.
Take $\alpha\in D^\R$ with $(E_{m-1,n})_\succ\subseteq\ker\alpha$
and $a_0\notin \ker\alpha$.  Set $\sigma:=\sigma_\alpha\in\mathcal D_{m,n}$. Then
$E_{m,n}\subseteq\T^\sigma$, so
\[\qquad\sigma(f)\  =\  \sum_{a} \sigma(f_a) \ex^{\alpha(a)+\sigma(a)}\ =\ \sum_{a} f_a\ex^{\alpha(a)}\ex^a\ \neq\ \sum_{a} f_a\ex^a\ =\ f.\qquad \qedhere\]
\end{proof}

\begin{lemma}\label{mgeqnlemma}
Suppose $m\geq n$. Then
$E{\downarrow_n}\cap \T^{\mathcal D_{m,n}}=E_{m,n}$.
\end{lemma}
\begin{proof}
Clearly $E{\downarrow_n}\cap \T^{\mathcal D_{m,n}}\supseteq E_{m,n}$. By induction on $k\in \N$ we show that
if~$f\in E_{m+k,n}\setminus E_{m,n}$, then $\sigma(f)\neq f$ for some $\sigma\in\mathcal D_{m,n}$.
This holds trivially if $k=0$, so suppose $k>0$ and 
$f\in E_{m+k,n}\setminus E_{m,n}$. If $f\in E_{m+k-1,n}$, then  we are done by inductive hypothesis.
So suppose $f\notin E_{m+k-1,n}$. Then Lemma~\ref{lemneq} gives~${\sigma\in\mathcal D_{m+k-1,n}}$ with $\sigma(f)\neq f$; now use $\mathcal D_{m+k-1,n}\subseteq\mathcal D_{m,n}$.
\end{proof}

\begin{lemma}\label{lem:intersect n+1 to n}
$E{\downarrow_{n+1}}\cap\T^{\mathcal D_n}=E{\downarrow_n}$.
\end{lemma}
\begin{proof}  
If $E{\downarrow}\cap\T^{\mathcal D_n}{\uparrow^n}=E$, then
$E{\downarrow_{n+1}}\cap\T^{\mathcal D_n}=E{\downarrow_n}$. Also
 $\T^{\mathcal D_n}{\uparrow^n}=\T^{\mathcal D_0}$ by Corollary~\ref{cor:fixed fields under updown shift}. Hence it suffices to do the case $n=0$. Now $\mathcal D_0=\mathcal{E}^{\R}|E$, so $\T^{\mathcal D_0}\supseteq E$. 
Suppose $f\in E_{m,1}\setminus E$; we claim that $\sigma(f)\neq f$  for some $\sigma\in\mathcal D_0$. We prove this claim by induction on $m$.
Let $m=0$.  With $r$ ranging over $\R$ we have $f=\sum_r f_r\ell_1^r$ with all $f_r\in\R$ and $f_r\neq 0$ for some $r\neq 0$. Choose $\alpha\in  D_0$ with $\alpha(\ell_1)\neq 0$. Then~$\sigma_\alpha(f)=\sum_r f_r(\ell_1-c)^r$ where $c=\alpha(\ell_1)\in\R^\times$, so $\sigma_\alpha(f) \neq f$. (Proof: Reduce to the case that $r_0:=\max \{r: f_r\ne 0\}\ne 0$. Then
expand $\sum_r f_r(\ell_1-c)^r$ into a series of powers of $\ell_1$ and check that the coefficient of $\ell_1^{r_0-1}$ differs from
$f_{r_0-1}$.) 

Suppose our claim holds for all $f\in E_{m,1}\setminus E$ for a certain $m$. Let~${f\in E_{m+1,1}\setminus E}$.
Let $a$ range over $A_{m,1}$.
We have $E_{m+1,1}=E_{m,1}[[\ex^{A_{m,1}}]]$, 
so $f=\sum_{a} f_a\ex^a$ with all~$f_a\in E_{m,1}$.
Suppose first that we have $a_0\in A_{m,1}$ with $f_{a_0}\neq 0$ and~${a_0\notin E_\succ}$.
Then $a_0\notin E_\succ+(E_{m-1,1})_\succ$ by Corollary~\ref{cor+:distributive}, so we have $\alpha\in D^\R$ with
$$E_\succ+(E_{m-1,1})_\succ\subseteq \ker\alpha, \qquad \alpha(a_0)\neq 0.$$
Then $\sigma_\alpha\in \mathcal D_{m,1}$, hence $E_{m,1}\subseteq \T^{\sigma_\alpha}$
and so~$\sigma_\alpha(f_a)=f_a$ and $\sigma_\alpha(a)=a$
for all~$a$. This yields
$\sigma_\alpha(f)=\sum_{a} f_a\ex^{\alpha(a)}\ex^a\neq \sum_{a} f_a\ex^a=f$.
Next assume $a\in E_\succ$ for all~$a$ with $f_a\neq 0$.
Then $\ex^a\in E$ for all~$a$ with~$f_a\neq 0$, so we can take
$a_0\in A_{m,1}$ with~$f_{a_0}\notin E$. By inductive hypothesis we obtain
$\sigma\in\mathcal D_0$ with $\sigma(f_{a_0})\neq f_{a_0}$.  
Hence~$\sigma(f)=\sum_a \sigma(f_a)\ex^a\neq f$ by Lem\-ma~\ref{lemt1}.
\end{proof}

\begin{cor}\label{cor:fixed field of T_n}
$\T^{\mathcal D_n}=E{\downarrow_n}$. In particular, $E{\downarrow_n}$ is definably closed in $\T$.
\end{cor}
\begin{proof}
By Corollary~\ref{cor:fixed fields under updown shift} it is enough to prove this for $n=0$.
By induction on $k\in \N$ we show $E{\downarrow_{k}}\cap\T^{\mathcal D_0}=E$. The case $k=0$ being trivial,
 suppose $E{\downarrow_{k}}\cap\T^{\mathcal D_0}=E$ for a certain $k$. Then
 \begin{align*}
 E{\downarrow_{k+1}}\cap\T^{\mathcal D_1} &\ =\  
 E{\downarrow_{k+1}}\cap\T^{\mathcal D_0}{\downarrow}\ =\ (E{\downarrow_{k}}\cap\T^{\mathcal D_0}){\downarrow}\ =\ E{\downarrow}\\
\intertext{and so by $\T^{\mathcal D_0}\subseteq \T^{\mathcal D_1}$ and Lemma~\ref{lem:intersect n+1 to n},}
 E{\downarrow_{k+1}}\cap\T^{\mathcal D_0}&\  =\  
 (E{\downarrow_{k+1}}\cap\T^{\mathcal D_1}) \cap \T^{\mathcal D_0}\  =\ 
 E{\downarrow} \cap \T^{\mathcal D_0}\  =\  E.\qedhere\end{align*}
\end{proof}

\noindent
In view of $\mathcal{D}_{m,n}\supseteq \mathcal{D}_n$ and Lemma~\ref{mgeqnlemma} this yields: 

\begin{cor}\label{cor:fixed field of Tmn}
If $m\geq n$, then $\T^{\mathcal D_{m,n}}=E_{m,n}$. 
\end{cor}

\noindent
Now we determine the fixed fields of $\mathcal E_n^\k$ and $\mathcal L_n^\k$.
We  have $\mathcal D_0\subseteq\mathcal E_{\infty}^\R\subseteq\mathcal E_{\infty}^\k$, so
by Corollary~\ref{cor:fixed field of T_n} we obtain
$E=\T^{\mathcal D_0}  \supseteq \T^{\mathcal E_{\infty}^\k}\supseteq \T^{\mathcal E_n^\k}$. Next we show:

\begin{lemma}\label{emtenk}\mbox{}
\begin{enumerate}
\item[\textup{(i)}]
If $m\leq n$, then $E_{m}\cap\T^{\mathcal E_n^\k} = \R[[\ex_m^\k\cdots\ex_0^\k]]$; and
\item[\textup{(ii)}]
 if $m\geq n$, then $E_{m}\cap\T^{\mathcal E_n^\k} = \R[[\ex_n^\k\cdots\ex_0^\k]]$.
 \end{enumerate}
\end{lemma}
\begin{proof}
Let $f\in E_{m}=\R[[G_{m}]]$.  Recalling that $E_{-1}=\R$, so $(E_{-1})_{\succ}=\{0\}$, we have~$G_{m}=x^\R \ex^{(E_{m-1})_\succ}$, so $f=\sum_{r,a} f_{r,a}x^r\ex^a$, all $f_{r,a}\in\R$,
with~$r$ ranging over $\R$ and $a$ over $(E_{m-1})_\succ$. 
For $\sigma\in\mathcal E^\k$ we
have $\sigma\big( (E_{m-1})_\succ \big) = (E_{m-1})_\succ$ 
by Lemma~\ref{lem:sigmaelln=elln} for $n=0$, hence $\sigma(x^r\ex^a)=\ex^{\alpha_{\sigma}(r\ell_1)-r\alpha_{\sigma}(\ell_1)}\ex^{\alpha_{\sigma}(a)}x^r\ex^{\sigma(a)}$ and thus
$$\sigma(f)_{r,\sigma(a)}\ =\ f_{r,a}\ex^{\alpha_{\sigma}(r\ell_1)-r\alpha_{\sigma}(\ell_1)+\alpha_\sigma(a)}.$$
Suppose $(r_0,a_0)\in \R\times (E_{m-1})_\succ$ and $f_{r_0,a_0}\neq 0$. The subset
$S:=\{a:f_{r_0,a}\neq 0\}$ of~$(E_{m-1})_\succ$ is countable, by \cite[2.10]{DMM2}. 
Hence if $a_0\notin\k\ex_{n-1}+\cdots+\k\ex_0$, then
we can take $\alpha\in D^\k$ such that $\alpha(\ex_{n-1})=\cdots=\alpha(\ex_0)=0$ and 
$$\ex^{\alpha(r_0\ell_1)-r_0\alpha(\ell_1)+\alpha(a_0)}\notin \{f_{r_0,a}/f_{r_0,a_0}:\ a\in S\}.$$
Then $\sigma=\sigma_\alpha\in\mathcal E_n^\k$ and 
$\sigma(f)_{r_0,\sigma(a_0)}\neq   f_{r_0,\sigma(a_0)}$, hence $\sigma(f)\neq f$.
So if $f\in \T^{\mathcal E_n^\k}$, then $f\in \R[[\ex_n^\k\cdots\ex_1^\k\ex_0^\R]]$.
Now let $f\in \R[[\ex_n^\k\cdots\ex_1^\k\ex_0^\R]]$, so $f=\sum_r f_rx^r$ where all~$f_r\in\R[[\ex_n^\k\cdots\ex_1^\k]]$.
Then for $\sigma\in\mathcal E_n^\k$ and $\alpha=\alpha_\sigma$ we have
$$\sigma(f)\ =\ \sum_r f_r \ex^{\alpha(r\ell_1)-r\alpha(\ell_1)} x^r$$
and hence $\sigma(f)=f$ iff $\alpha(r\ell_1)=r\alpha(\ell_1)$ for all $r$ with $f_r\neq 0$.
Now suppose $r_0\in\R\setminus\k$ and $f_{r_0}\neq 0$.
Then we can take $\alpha\in D^\k$ such that
$$\k\ell_1+\k\ex_{n-1}+\cdots+\k\ex_0\ \subseteq\ \ker\alpha\quad\text{ and }\quad
\alpha(r_0\ell_1)\neq 0,$$
hence $\sigma=\sigma_\alpha\in\mathcal E_n^\k$ and $\sigma(f)\neq f$.
Thus if $f\in\T^{\mathcal E_n^\k}$, then $f\in \R[[\ex_n^\k\cdots\ex_0^\k]]$.
This shows $E_{m}\cap\T^{\mathcal E_n^\k} = E_{m}\cap\R[[\ex_n^\k\cdots\ex_0^\k]]$,
and the latter equals $\R[[\ex_m^\k\cdots\ex_0^\k]]$ if $m\leq n$ (by Corollary~\ref{cor:Em it exp}) and
$\R[[\ex_n^\k\cdots\ex_0^\k]]$ if $m\geq n$.
\end{proof}

\noindent
By Lemma~\ref{emtenk} and the remark preceding it we have $\T^{\mathcal E_n^\k} = \R[[\ex_n^\k\cdots\ex_0^\k]]$. Therefore~$\R[[\ex_n^\k\cdots\ex_0^\k]]$ is definably closed in~$\T$. 
In view of Lemma~\ref{lem:fixed fields under updown shift}  this yields:

\begin{cor}\label{cor:TLnk} 
$\T^{\mathcal L_n^\k}=\R[[\ell_0^\k\cdots\ell_n^\k]]$,
so $\R[[\ell_0^\k\cdots\ell_n^\k]]$ is definably closed in $\T$.
\end{cor}

\noindent
By this corollary, 
$\T_{\log}^\k=\bigcup_n \R[[\ell_0^\k\cdots\ell_n^\k]]$ is definably closed in $\T$ as well.
The next lemma strengthens this:

\begin{lemma}\label{lem:TLk}
We have
$\T^{\mathcal L_{\infty}^{\k}}\cap E{\downarrow_n}=\R[[\ell_0^\k\cdots\ell_n^\k]]$, hence
$\T^{\mathcal L_{\infty}^{\k}}=\T_{\log}^\k$.
\end{lemma}
\begin{proof}
Let $f\in \T^{\mathcal L_{\infty}^{\k}}\cap E{\downarrow_n}$ and   $\sigma\in\mathcal L_n^\k$, so $\alpha:=\alpha_\sigma\in D^\k$ and 
$\ker\alpha\supseteq\bigoplus_{m=1}^n\k\ell_m$.
Extend $\alpha|_{E{\downarrow_n}}$ to $\hat\alpha\in D^\k$ such that $\ker\hat\alpha\supseteq \bigoplus_{m>n} \k\ell_m$.
Then $\hat\sigma:=\sigma_{\hat\alpha}\in\mathcal L_{\infty}^{\k}$ and~$\sigma|_{E{\downarrow_n}}=\hat\sigma|_{E{\downarrow_n}}$,
thus $\sigma(f)=\hat\sigma(f)=f$. By the previous corollary this yields~$\T^{\mathcal L_{\infty}^{\k}}\cap E{\downarrow_n}\subseteq \T^{\mathcal L_n^\k}=\R[[\ell_0^\k\cdots\ell_n^\k]]$, and $\R[[\ell_0^\k\cdots\ell_n^\k]]\subseteq \T^{\mathcal L_{\infty}^{\k}}$ is clear.
\end{proof}

\noindent
Let $\SigmaAut_\der^\k(\T^\k_{\log})$ be the subgroup of $\SigmaAut_\der(\T^\k_{\log})$ consisting of the $\sigma\in \SigmaAut_\der(\T^\k_{\log})$ such that
$ \sigma(f^r)=\sigma(f)^r$ for all $f>0$ in $\T^\k_{\log}$ and all $r\in \k$. To $\sigma\in \SigmaAut_\der(\T^\k_{\log})$ we associated in Section~\ref{sec:Tlog} an additive map $\alpha\colon \log\big((\T^\k_{\log})^{>}\big) \to \R$, and here we note that~$\log\big((\T^\k_{\log})^{>}\big)$ is a $\k$-linear subspace of $\T^\k_{\log}$, and that $\sigma\in \SigmaAut_\der^\k(\T^\k_{\log})$ iff~$\alpha$ is $\k$-linear.  If $\alpha$ is $\k$-linear, it extends to an element of $D^\k$. Using also Corollary~\ref{corsad} we see that each $\sigma\in \SigmaAut_\der^\k(\T^\k_{\log})$ extends to a strong $\k$-power preserving automorphism
of $\T$. This yields a short exact sequence of groups: 
$$1 \to \mathcal L_{\infty}^\k \xrightarrow{\ \iota\ } \SigmaAut^\k_\der(\T) \xrightarrow{\ \rho\ } \SigmaAut_\der^\k(\T^\k_{\log})\to 1$$
where $\iota$ is the natural inclusion and $\rho(\sigma):=\sigma|_{\T^\k_{\log}}$ for each $\sigma\in \SigmaAut^\k_\der(\T)$.
Does this sequence split (so that  $\SigmaAut^\k_\der(\T)  = \mathcal L_{\infty}^\k\ltimes\SigmaAut_\der^\k(\T^\k_{\log})$)?

\begin{remarks}
Corollary~C in the Introduction now follows for $\R[[\ex_n^\R\cdots\ex_0^\R]]$ from remarks after Lem\-ma~\ref{emtenk}, for $\R[[\ell_0^\R\cdots\ell_n^\R]]$ from
Corollary~\ref{cor:TLnk}, for $\T_{\exp}$ from Corollary~\ref{cor:fixed field of T_n}, 
and for $\T_{\log}$ and  $\T^{\Q}_{\log}$ from Lem\-ma~\ref{lem:TLk}.
The proofs show that,  except for $\T^{\Q}_{\log}$, these differential subfields of $\T$ are
even definably closed in the expansion of the differential field $\T$ by the constant power map
$(r,f)\mapsto f^r\colon \R\times \T^> \to \T^>$, whose model theory is studied in \cite{Kaplan}.
\end{remarks}

\subsection*{Fixed fields of subgroups of $\mathcal L_{\infty}$} Recall that 
$$\mathcal L_{\infty}\ =\ \big\{\sigma\in \SigmaAut_{\der}(\T):\,\text{$\sigma(\ell_n)=\ell_n$  for all $n$}\big\}\ \subseteq\  \mathcal{E}.$$
In this subsection we investigate the structure of the fixed fields of
subsets of  $\mathcal{L}_{\infty}$. 
Let~$\mathcal{G} \subseteq \mathcal L_{\infty}$ (so $\T^{\mathcal{G}}\supseteq\T^{\mathcal L_{\infty}}=\T_{\log}^\Q$).
By Corollary~\ref{cor:sigmaelln=elln},  
$\sigma(\R^> \cdot G^{\operatorname{LE}})=\R^>\cdot G^{\operatorname{LE}}$ for~$\sigma\in\mathcal L_{\infty}$. It follows that
$\fd (\T^{\mathcal G})^\times = G^{\operatorname{LE}}\cap \T^{\mathcal G}$, and thus
$$(\T^{\mathcal G})^\times\  =\  \R^\times\cdot (G^{\operatorname{LE}}\cap \T^{\mathcal G})\cdot \big(1+(\T^{\mathcal G})^\prec\big).$$
By Corollary~\ref{cortr}, $\T^{\mathcal G}$ is truncation closed, so $\supp f\subseteq  \T^{\mathcal G}$ for all $f\in \T^{\mathcal{G}}$.

\begin{lemma}\label{lem:fixsigma, 2}
Let $\sigma\in\mathcal L_{\infty}$, $\alpha=\alpha_\sigma$. Then 
$$ G^{\operatorname{LE}}\cap \T^{\sigma}\  =\  \big\{\fm\in 
G^{\operatorname{LE}}:  \log\fm\in\T^{\sigma}\cap\ker\alpha\big\},\qquad
\log\!\big( (\T^{\sigma})^> \big)\  =\  \T^\sigma\cap\ker\alpha.$$
\end{lemma}
\begin{proof}
Let $\fm\in G^{\operatorname{LE}}$. Then 
$\sigma(\log\fm)=\log\sigma(\fm)-\alpha(\log\fm)$ by Lemma~\ref{lem:sigma alpha}(ii).
In view of $\log G^{\operatorname{LE}} = \T_\succ=\sigma(\T_\succ)$ it follows that
$\sigma(\fm)=\fm$ iff  $\sigma(\log\fm)=\log\fm$ and~$\alpha(\log\fm)=0$.
This proves the first equality.
Let $c\in\R^>$, $\fm\in G^{\operatorname{LE}}$, $\varepsilon\in\T^\prec$, and~$f=c\fm(1+\varepsilon)$. Then 
$\log f=\log c+\log \fm+\log(1+\varepsilon)$ 
where $\log c\in\R$, $\log\fm\in\T_\succ$, and $\log(1+\varepsilon)\in\T^\prec$.
We have
\begin{align*}
f\in\T^\sigma	&\quad\Longleftrightarrow\quad \fm\in\T^\sigma \text{ and } \varepsilon\in\T^\sigma \\
				&\quad\Longleftrightarrow\quad \log\fm\in\T^\sigma\cap\ker\alpha
				 \text{ and } \log(1+\varepsilon)\in\T^\sigma,
\end{align*}
and the second equality follows.
\end{proof}

\noindent
For $\fm\in G^{\operatorname{E}}{\downarrow_n}$ we have  $\fm=\ell_n^r\ex^a$ for a unique pair $(r,a)\in \R\times (E{\downarrow_n})_{\succ}$, and we set~$\fm_n:=\ex^a$. So for $\fm\in G_{m+1,n}$ and $\sigma\in \mathcal{L}_{\infty}$ we have $\log \fm_n\in E_{m,n}$, and~$\sigma(\log \fm)=\log \fm$ iff $\sigma(\log \fm_n)=\log \fm_n$.

By recursion on $m$ we can now define for $f\in E_{m,n}$ a set $\supp_{m,n} f\subseteq G_{m,n}$ as follows:
if $m=0$, so $E_{0,n}=\R[[\ell_n^\R]]$, then $\supp_{m,n} f:=\supp f\subseteq \ell_n^{\R}$, and for~$f\in E_{m+1,n}=\R[[G_{m+1,n}]]$ we set
$$\supp_{m+1,n} f\ :=\ \supp f\ \cup \bigcup_{\fm\in \supp f} \supp_{m,n} \log\fm_n.$$
The sets $\supp_{m,n}$ for $f\in E_{m,n}$ are countable  by \cite[2.10]{DMM2}. 
Example:  $f=\ell_n\in E_{0,n}$ gives $\supp_{0,n} f = \{\ell_{n}\}$, but 
 $\supp_{1,n+1} f = \supp_{1,n+1} \ex^{\ell_{n+1}} = \{ \ell_n, \ell_{n+1} \}$.

\begin{lemma}\label{lem:criterion for sigma(f)=f}
Let $\sigma\in\mathcal L_{\infty}$, $\alpha=\alpha_\sigma$. Then for all $f\in E_{m,n}$:
$$\sigma(f)=f \quad\Longleftrightarrow\quad 
\log\supp_{m,n} f \subseteq \ker\alpha.$$
\end{lemma}
\begin{proof} By induction on $m$. Let $f\in\ E_{0,n}$, so with $r$ ranging over $\R$
we have~$f=\sum_r f_r\ell_n^r$, all $f_r\in\R$. Then $\sigma(f)=\sum_r f_r \ex^{\alpha(r\ell_{n+1})}\ell_n^r$,
so 
\begin{align*}
\sigma(f) = f	&\quad \Longleftrightarrow\quad \text{$\alpha(r\ell_{n+1})=0$ for all $r$ with $f_r\neq 0$} \\
			&\quad\Longleftrightarrow\quad \log\supp_{0,n} f=\{r\ell_{n+1}:f_r\neq 0\}\subseteq\ker\alpha, 
\end{align*}
as required. 
Suppose $m$ is such that the claimed equivalence holds for all $n$ and~${f\in E_{m,n}}$.
Let $f\in E_{m+1,n}$, and
let $\fm$ range over~$\supp f\subseteq G_{m+1,n}$.
Then 
$$f\ =\ \sum_{\fm} f_\fm\fm,\qquad \sigma(f)\ =\ \sum_{\fm} f_\fm\ex^{\alpha(\log\fm)}\ex^{\sigma(\log\fm)}.$$
Assume $\sigma(f)=f$. Then $f\in \T^{\sigma}$, so for all $\fm$ we have $\fm\in\T^{\sigma}$,  hence
$\log\fm\in\ker\alpha$ and $\sigma(\log\fm)=\log\fm$
by Lemma~\ref{lem:fixsigma, 2}. 
Using the  inductive hypothesis  we obtain~$\log\supp_{m+1,n} f \subseteq \ker\alpha$.
Conversely, if this inclusion holds, then~${\alpha(\log\fm)=0}$ as well as $\sigma(\log\fm)=\log\fm$ by inductive hypothesis, and hence~${\sigma(f)=f}$.
\end{proof}

\noindent
For $F\subseteq\T$ we put  $\mathcal L_{\infty}|F:=\SigmaAut_\der(\T|F)\cap\mathcal L_{\infty}$ and $\mathcal E|F:=\SigmaAut_\der(\T|F)\cap\mathcal E$.
As explained in \cite[Section~5]{ADHdim}, the following proves Corollary~D from the Introduction.

\begin{cor}\label{cor:move exps}
Let $F\subseteq \T$  be countable. Then for some subgroup~$R$ of $\R$
with~$|R|=2^{\aleph_0}$ there is for every $r\in R^{\neq}$ a $\sigma\in\mathcal L_{\infty}|F$ such that $\sigma(\ex^{rx})\neq \ex^{rx}$.
\end{cor}
\begin{proof}
For each $f\in F$, choose $m(f), n(f)\in \N$ such that $f\in E_{m(f),n(f)}$, and put
$$S\ :=\ \bigcup_{f\in F} \log \supp_{m(f),n(f)} f\  \subseteq\  \T_\succ.$$
Then $S$ is countable, hence so is the $\Q$-linear subspace $M$ of $\T_\succ$ generated by ${S\cup\{\ell_n: n\ge 1\}}$. 
Choose a $\Q$-linear subspace $R$ of $\R$ such that $\R  = R \oplus ({\R\cap x^{-1}M})$. Then $|R|=2^{\aleph_0}$.
Given $r\in R^{\neq}$, choose a $\Q$-linear map $\alpha\colon\T\to\R$ such that $M\oplus \T^{\preceq}\subseteq\ker\alpha$
and $\alpha(rx)\neq 0$. Then  $\sigma:=\sigma_\alpha\in \mathcal{L}_{\infty}$ gives $\sigma(f)=f$ for all~$f\in F$, by Lemma~\ref{lem:criterion for sigma(f)=f}, and
$\sigma(\ex^{rx})=\ex^{\alpha(rx)+\sigma(rx)}=\ex^{\alpha(rx)}\cdot \ex^{rx}\neq\ex^{rx}$.
\end{proof}

\noindent
For use in a later paper we include a modified form of this corollary:

\begin{cor}
Let $F\subseteq \T$ and $G\subseteq\T^{>\R}$  be countable. Then for some subgroup~$R$ of $\R$
with~$|R|=2^{\aleph_0}$ there is for all $g\in G$, $r\in R^{\neq}$ a $\tau\in\mathcal E|F$ such that~$\tau(\ex^{rg})\neq \ex^{rg}$.
\end{cor}
\begin{proof}
Arrange   $x\in F$, and take $R$ as in Corollary~\ref{cor:move exps} with  the countable subset
$$F\circ G^{\inv}\ :=\ \{f\circ g^{\inv}:\,f\in F, g\in G\}$$ of $\T$ in place of $F$.
Let $g\in G$, $r\in R^{\ne}$, and take 
$\sigma\in\mathcal E|(F\circ G^{\inv})$ with~$\sigma(\ex^{rx})\ne \ex^{rx}$. Put $\tau:=g^{\inv}\sigma\in\mathcal E$.
By Lemma~\ref{lem:T under comp} we have for $f\in\T$, 
$$\sigma(f\circ g^{\inv})\ =\ \tau(f)\circ\sigma(g^{\inv})\ =\ \tau(f)\circ g^{\inv}.$$
In particular, $\ex^{rx}\ne\sigma(\ex^{rx})=\sigma(\ex^{rg}\circ g^{\inv})=\tau(\ex^{rg})\circ g^{\inv}$
and thus $\ex^{rg}\ne\tau(\ex^{rg})$.
\end{proof}

\noindent
Suppose now that we are given a strong differential subfield $K$ of
$\T$ such that
$$K^\times\ =\  \R^\times\cdot G\cdot  (1+K^\prec), $$
where $G$ is a divisible subgroup of $G^{\operatorname{LE}}$ such that $x\in G$ and $\log G\subseteq K$. Therefore~$G=\fd(K^\times)\subseteq K^\times$, and so $K=\T\cap\R[[G]]$ by Lemma~\ref{lem:fixsigma, 1}. 
From $\log G^{\operatorname{LE}} = \T_\succ$ we obtain
$$\log( K^>)\  =\  \log G \oplus \R \oplus K^\prec\  \subseteq\  K.$$
Also,  $x\in G$ and $\log G\subseteq K$ yields $\ell_n\in G$ for all $n$, so $\T_{\log}^\Q\subseteq K$.
Set
$$\mathcal G\ :=\ \{\sigma\in\mathcal E:\,  \ker\alpha_\sigma \supseteq  \log G\}.$$
Then $\mathcal G\subseteq\mathcal L_{\infty}$. We prove below that $K=\T^{\mathcal G}$. First:

\begin{lemma}\label{lem:K subset TG}
$K\subseteq\T^{\mathcal G}$.
\end{lemma}
\begin{proof}
Let $\sigma\in\mathcal G$;
by induction on $m$ we prove $K\cap E_{m,n}\subseteq\T^{\sigma}$.
Let $f\in K\cap E_{m,n}$.
Since $K=\T\cap\R[[G]]$    we have
$\log\supp f\subseteq\log G\subseteq \ker\alpha_\sigma$. If $m=0$, this already gives
$f\in \T^\sigma$, by Lemma~\ref{lem:criterion for sigma(f)=f}.
Suppose $m$ is such that $K\cap E_{m,n}\subseteq\T^{\sigma}$, and let~$f\in K\cap E_{m+1,n}$. Then  for $\fm\in \supp f$
we have $$\log \fm\in K, \qquad \log \fm\in \R\ell_{n+1}+\log \fm_n,$$ so $\log \fm_n\in K\cap E_{m,n}\subseteq \T^{\sigma}$. 
Hence by Lemma~\ref{lem:criterion for sigma(f)=f}
$$\log\supp_{m,n} \log\fm_n\ \subseteq\ \ker\alpha_\sigma\quad \text{ for all }\fm\in \supp f. $$
Hence $\log \supp_{m+1,n} f \subseteq \ker\alpha_\sigma$,
so $f\in\T^{\sigma}$.
\end{proof}

\begin{cor}
$\mathcal G$ is a subgroup of $\mathcal L_{\infty}$.
\end{cor}
\begin{proof}
Let $\sigma,\tau\in\mathcal G$. From  $\alpha_{\sigma\circ\tau}=\alpha_\sigma\circ\tau+\alpha_\tau$,
and  $\log G\subseteq K\subseteq\T^\tau$ we obtain~$\ker\alpha_{\sigma\circ\tau}\supseteq\log G$, so $\sigma\circ\tau\in\mathcal G$.
Next, let $\tau\in \mathcal{G}$ and $\sigma=\tau^{-1}$. Then $\alpha_\sigma=-\alpha_\tau\circ\sigma$, so 
in view of   $\log G\subseteq K\subseteq\T^\tau=\T^{\sigma}$ we obtain $\ker\alpha_\sigma\supseteq\log G$.
\end{proof}

\noindent
Recall the terminology ``closed'' from the discussion before Lemma~\ref{lemgalois}.

\begin{cor}
The subgroup
$\mathcal G$ is  of $\SigmaAut_\der(\T)$ is closed \textup{(}and hence the differential subfield $\T^{\mathcal G}$ of $\T$ is closed\textup{)}.
\end{cor}
\begin{proof}
We have $\mathcal G\subseteq \operatorname{c}(\mathcal G)=\SigmaAut_\der(\T|\T^{\mathcal G})$; it suffices to show $\operatorname{c}(\mathcal G)\subseteq\mathcal G$.
For this, let $\sigma\in\operatorname{c}(\mathcal G)$. Then $\T_{\log}^{\Q}\subseteq K\subseteq \T^{\mathcal{G}}$ gives
 $\sigma\in\mathcal L_{\infty} \subseteq\mathcal E$. Moreover, if
$g\in \log G$, then~$\ex^g\in G\subseteq K\subseteq\T^{\mathcal G}$ and~$g\in K\subseteq\T^{\mathcal G}$,  so
$\ex^g=\sigma(\ex^g)=\ex^{\alpha_\sigma(g)+\sigma(g)}=\ex^{\alpha_\sigma(g)}\cdot\ex^g$
and thus $\alpha_\sigma(g)=0$. This yields $\sigma\in \mathcal G$.
\end{proof}

\begin{lemma}
$\T^{\mathcal G}\subseteq K$.
\end{lemma}
\begin{proof}
By induction on $m$ we show $\T^{\mathcal G}\cap E_{m,n}\subseteq K$. Now $\T^{\mathcal G}\cap E_{m,n}$ is a strong $\R$-linear subspace of $\T$, so
by Lemma~\ref{lem:fixsigma, 1} it is enough to show that
$\T^{\mathcal G}\cap G_{m,n}\subseteq K$. We have $G_{0,n}=\ell_n^\R$. Suppose $\fm=\ell_n^r\notin K$, $r\in \R$. Then
$r\ell_{n+1}=\log \ell_n^r\notin\log G$, so we have $\alpha\in D$ such that $\ker\alpha\supseteq \log G$ and $\alpha(r\ell_{n+1})\neq 0$. Hence $\sigma_{\alpha}\in\mathcal G$ and~$\sigma_{\alpha}(\fm)\neq\fm$, so $\fm\notin \T^{\mathcal G}$. This argument yields $\T^{\mathcal G}\cap G_{0,n}\subseteq K$.

Assume $m$ is such that $\T^{\mathcal G}\cap G_{m,n}\subseteq K$. 
We have $G_{m+1,n}=\ell_n^\R \ex^{(E_{m,n})_\succ}$.
Let~$\fm=\ell_n^r\ex^a\in G_{m+1,n}$, where $r\in\R$ and $a\in (E_{m,n})_\succ$.
Then for $\sigma\in\mathcal L_{\infty}$, $\alpha=\alpha_\sigma$, we have~$\sigma(\fm)=\ex^{\alpha(r\ell_{n+1})+\alpha(a)}\ell_n^r\ex^{\sigma(a)}$
where $\sigma(a)\in (E_{m,n})_\succ$ by Lemma~\ref{lem:sigmaelln=elln}, hence 
$$\sigma(\fm)=\fm\quad\Longleftrightarrow\quad \alpha(r\ell_{n+1}+a)=0 \text{ and } \sigma(a)=a.$$
Suppose $\fm\in\T^{\mathcal G}$; our aim is to show $\fm\in K$. If
$\log \fm=r\ell_{n+1}+a\in \log G$, then~$\fm\in G\subseteq K$ and we are done. So assume $r\ell_{n+1}+a\notin \log G$.
We do have~${a\in \T^{\mathcal G}}$, so 
$a\in K\cap \log G^{\operatorname{LE}}=K\cap \T_\succ$ by the inductive assumption. Take $\alpha\in D$ with $\ker \alpha\supseteq \log G$ and $\alpha(r\ell_{n+1}+a)\ne 0$. Then $\sigma_{\alpha}\in \mathcal{G}$, but $\sigma(\fm)\ne \fm$, a contradiction.
\end{proof}

\begin{cor}
The fixed fields of subgroups of~$\mathcal L_{\infty}$ are exactly the 
strong differential subfields $F$ of $\T$ 
such that
$F^\times = \R^\times\cdot G\cdot  (1+F^\prec)$ for some divisible subgroup~$G$ of~$G^{\operatorname{LE}}$ with $x\in G$ and $\log G\subseteq F$.
\end{cor}
 
\noindent
For any $\mathcal G\subseteq \mathcal L_{\infty}$, the fixed field $F:=\mathbb T^{\mathcal G}$ is newtonian, by Corollary~\ref{cor:newt}.

\subsection*{Application to $\T^\flat$}
An instructive example is provided by the differential field~$\T^\flat$ of flat transseries.
Recall from Section~\ref{sec:prelims, trans} that the flat part   of $G:=G^{\operatorname{LE}}$ is its   subgroup~$G^\flat=\{\fm\in G:\,\fm\flatter \ex^x\}=\{\fm\in G:\, 1\prec \supp \log \fm \prec x\}$,
and~$\T^\flat=\T\cap\R[[G^\flat]]$ is a strong differential subfield of $\T$ with
$$\big(\T^\flat\big)^\times\  =\  \R^\times \cdot G^\flat \cdot (1+(\T^\flat)^\prec)
\quad\text{and}\quad \log\big((\T^\flat)^>\big)\ \subseteq\ \T^\flat.$$
Since the subgroup $G^\flat$ of $G$ is convex, it is closed under powers and hence divisible, and  $x\in G^\flat$.
Hence $\T^\flat$ is the fixed field of a closed subgroup~$\mathcal G^\flat$ of $\SigmaAut_\der(\T)$, namely
$$\mathcal G^\flat\ :=\   \{\sigma\in\mathcal E:\  \ker \alpha_\sigma\supseteq\log G^\flat \}\ \subseteq\  \mathcal{L}_{\infty}.$$
Thus $\mathcal G^\flat=\SigmaAut_{\der}(\T|\T^\flat)$ and $\T^\flat$ is definably closed in $\T$.
Recall that
$G_{m,n}^\flat:=G_{m,n}\cap G^\flat$, and
$$G_{0,n}^\flat\ =\ G_{0,n}\ =\ \ell_n^{\R},\qquad G_{m+1,n}^\flat\ =\ \{\fm\in G_{m+1,n} : 1\prec\supp \log \fm \prec x\}.$$
Moreover, with $E^\flat_{m,n}:=E_{m,n}\cap\T^\flat=\R[[G_{m,n}^\flat]]$
we have  $\T^\flat=\bigcup_{m,n} E^\flat_{m,n}$.
 
\begin{prop}\label{proptflat} Let $\sigma\in \SigmaAut_\der(\T)$. Then $\sigma(\T^\flat)\subseteq \T^\flat$. 
\end{prop}
\begin{proof}
By Proposition~\ref{prop:5.11} this holds for $\sigma\in\mathcal T$. 
Let $\sigma\in\mathcal E$, $\alpha:=\alpha_\sigma$; we claim that~${\sigma(\T^\flat)\subseteq\T^\flat}$.
We fix $n$ and prove by induction on $m$ that $\sigma(E^\flat_{m,n})\subseteq\T^\flat$.
Let~$\fn=\ell_n^r$ ($r\in\R$).  By Lemma~\ref{lem:sigma elli} we have
$$\sigma(r\ell_{n+1})\ =\ r\sigma(\ell_{n+1})\ =\ r\ell_{n+1}-r\alpha(\ell_{n+1})+\varepsilon\quad\text{ where $\varepsilon\in\R[[\ell_0^\Z\cdots\ell_n^\Z]]_{\prec}\subseteq \T^\flat_\prec$.}$$
Hence
$$\sigma(\fn)\ 	=\	\sigma(\ex^{r\ell_{n+1}})=\ex^{\alpha(r\ell_{n+1})+\sigma(r\ell_{n+1})}\ 
			=\	\ex^{\alpha(r\ell_{n+1})-r\alpha(\ell_{n+1})} \fn \sum_n \frac{\varepsilon^n}{n!}\in \T^\flat. $$
By strong linearity of $\tau$ this yields $\tau(E_{0,n})=\tau(E_{0,n}^\flat)\subseteq\T^\flat$.
Next, let $\fn\in G_{m+1,n}^\flat$, so~$\fn=\ell_n^r\ex^a$ where $r\in \R$ and $a\in (E^\flat_{m,n})_\succ$. Inductively,
 $\sigma(a)\in\T^\flat$,
so $\sigma(a)=\varphi+c+\delta$ where~$\varphi\in \T^\flat_\succ$,  $c\in\R$, $\delta\in \T^\flat_\prec$.
We have $a\prec x$, so $\sigma(a)\prec x$, hence
$\varphi\prec x$, and thus $\ex^\varphi\in G^\flat$. Therefore,
$$\sigma(\ex^a)\ =\ \ex^{\alpha(a)+\sigma(a)}\  =\  \ex^{\alpha(a)+c} \ex^{\varphi} \sum_n \frac{\delta^n}{n!}\in\T^\flat.$$
Together with strong linearly of $\sigma$ and $\sigma(\ell_n^{\R})\subseteq\T^\flat$ this yields
$\sigma(E_{m+1,n}^\flat)\subseteq\T^\flat$.
\end{proof}

\begin{cor}
The group morphism
\begin{equation}\label{eq:restrictflat}
\sigma\mapsto\sigma|_{\T^\flat}\  \colon\  \SigmaAut_\der(\T)\to\SigmaAut_\der(\T^\flat)
\end{equation}
has kernel $\mathcal G^\flat$ $($so $\mathcal G^\flat$ is a normal subgroup of $\SigmaAut_\der(\T))$ and is surjective.
\end{cor}
\begin{proof}
By Lemma~\ref{lemgalois}(i),  the group morphism \eqref{eq:restrictflat}  
has kernel $\mathcal G^\flat$. 
It remains to show that it is surjective.
Let $\sigma\in\SigmaAut_\der(\T^\flat)$;
we shall find $\tau\in\SigmaAut_\der(\T)$ such that $\sigma=\tau|_{\T^\flat}$.
Replacing~$\sigma$ by~$\sigma \circ(\sigma_{-c})|_{\T^\flat}$ where $c:=\sigma(x)-x\in\R$ we first arrange $\sigma(x)=x$.
Now consider the $\Q$-linear subspace $\log\big( (\T^\flat)^>\big)$ of $\T^\flat$. We have~$\log\big((\T^\flat)^>\big)\cap \T_{\preceq}=\T^\flat_\preceq:=\T^\flat\cap\T_\preceq$. We introduce the $\Q$-linear function
$$\alpha\ \colon\ \log\big((\T^\flat)^>\big)\to\R, \qquad
\alpha(f)\  :=\  \log \sigma(\ex^f)-\sigma(f),$$
 which satisfies $\ker\alpha\supseteq\T^\flat_\preceq$. 
Choose a $\Q$-linear subspace~$V$ of $\T_{\preceq}$ with $\T_{\preceq}=V\oplus\T^\flat_{\preceq}$;
then $V\cap \log\big((\T^\flat)^>\big)= \{0\}$.
Next choose a $\Q$-linear subspace $W$ of~$\T$ with~$\T=\left(V\oplus\log\big( (\T^\flat)^>\big)\right)\oplus W$, and
extend $\alpha$ to the $\Q$-linear map~$\T\to\R$, also denoted by~$\alpha$, with $\alpha|_{V\oplus W}=0$. Then $\alpha\in D$,
and we claim that with $\tau:=\sigma_\alpha\in\mathcal E$ we have~$\sigma=\tau|_{\T^\flat}$.
We have $\T^{\R}_{\log}\subseteq\T^\flat$, and by Corollary~\ref{corsad}, 
$\sigma$ and $\tau$ both restrict to strong automorphisms of $\T^{\R}_{\log}$. 
Since $L_{\R}=\bigoplus_{n\ge 1}\R\ell_n\subseteq\log\big( (\T^\flat)^>\big)$, 
Section~\ref{sec:Tlog} shows that $\sigma$ and $\tau$ agree on $\T^{\R}_{\log}$.
Fixing $n$, induction on $m$ as in the proof of Proposition~\ref{proptflat} now  shows that $\sigma$ and $\tau$ agree on $E_{m,n}^\flat$.
\end{proof}
 
\bibliographystyle{amsplain}

\end{document}